\documentclass[12pt, a4paper,]{article}

\usepackage[margin=1in]{geometry}
\usepackage{amsmath,amssymb,amsfonts}
\usepackage{amsthm}
\usepackage{booktabs}
\usepackage{multirow,multicol}
\usepackage{diffcoeff}
\usepackage{mathtools}
\usepackage{hyperref}
\usepackage{cleveref}
\usepackage{algorithm}
\usepackage{algpseudocodex}
\algrenewcommand\algorithmicrequire{\textbf{Input:}}
\algrenewcommand\algorithmicensure{\textbf{Output:}}
\usepackage{marvosym}
\usepackage{graphics}
\usepackage{graphicx}
\usepackage{blindtext}
\usepackage{wrapfig}
\usepackage{url}
\usepackage{comment}
\usepackage{tikz}
\usepackage{tkz-base,tkz-euclide}
\usepackage{tikz-3dplot}
\usepackage{pgfplots}
\pgfplotsset{compat=1.18}
\usetikzlibrary{calc,intersections,through,arrows.meta,external,shapes.geometric}
\usepackage{stmaryrd}
\SetSymbolFont{stmry}{bold}{U}{stmry}{m}{n}
\usepackage{nicefrac}
\usepackage{listings}
\usepackage{xcolor}
\usepackage{geometry}
\usepackage{tabularx}
\usepackage{rotating}
\usepackage{enumitem}
\usepackage{lmodern}
\usepackage{placeins}
\usepackage{orcidlink}
\usepackage{bm}

\usepackage[T1]{fontenc}
\usepackage[english]{babel}

\usetikzlibrary{angles,quotes}

\usepackage[%
backend=biber, 
sorting=nyt, 
sortcites, 
safeinputenc, 
citestyle=numeric, %
bibstyle=ext-numeric, %
hyperref=true, 
maxbibnames=99, 
maxcitenames=2, 
url=false, 
isbn=false,
doi=true, 
giveninits=true,
date=year
]%
{biblatex}

\DeclareFieldFormat{titlecase:title}{\MakeSentenceCase*{#1}}

\renewbibmacro{in:}{%
	\ifentrytype{article}{}{\printtext{\bibstring{in}\intitlepunct}}}

\usepackage{lineno}
\usepackage{etoolbox}

\definecolor{mutedRose}{HTML}{CC6677}
\definecolor{mutedIndigo}{HTML}{332288}
\definecolor{mutedSand}{HTML}{DDCC77}
\definecolor{mutedGreen}{HTML}{117733}
\definecolor{mutedCyan}{HTML}{88CCEE}
\definecolor{mutedWine}{HTML}{882255}
\definecolor{mutedTeal}{HTML}{44AA99}
\definecolor{mutedOlive}{HTML}{999933}
\definecolor{mutedPurple}{HTML}{AA4499}
\definecolor{mutedGrey}{HTML}{DDDDDD}

\definecolor{brightblue}{HTML}{4477AA}
\definecolor{brightcyan}{HTML}{66CCEE}
\definecolor{brightgreen}{HTML}{228833}
\definecolor{brightyellow}{HTML}{CCBB44}
\definecolor{brightred}{HTML}{EE6677}
\definecolor{brightpurple}{HTML}{AA3377}
\definecolor{brightgrey}{HTML}{BBBBBB}

\theoremstyle{thmstyleone}
\newtheorem{theorem}{Theorem}
\newtheorem{proposition}[theorem]{Proposition}
\newtheorem{lemma}[theorem]{Lemma}
\newtheorem{corollary}[theorem]{Corollary}

\theoremstyle{thmstyletwo}
\newtheorem{example}{Example}

\theoremstyle{thmstylethree}
\newtheorem{definition}{Definition}
\newtheorem{assumption}{Assumption}

\newcommand{\T}{\mathsf{T}}

\begin{document}

\title{Approximation Algorithms for Right-Hand-Side Parametric Linear Programming Problems}

\author{Levin Nemesch\footnote{RPTU Kaiserslautern-Landau, Department of Mathematics, Paul-Ehrlich-Str. 31, 67663 Kaiserslautern, Germany, email:\{l.nemesch,stefan.ruzika\}@math.rptu.de}, Stefan Ruzika\footnotemark[\value{footnote}]~,\\
     Clemens Thielen\footnote{Technical University of Munich, Campus Straubing for Biotechnology and Sustainability, Professorship of Optimization and Sustainable Decision Making, Am Essigberg~3, 94315~Straubing, Germany, email:\{clemens.thielen,alina.wittmann\}@tum.de}, Alina Wittmann\textsuperscript{\Letter}\footnotemark[\value{footnote}]
}

\date{\today}

\maketitle

\begin{abstract}
In parametric linear optimization problems, the objective function or the feasible set depends linearly on one or several real-valued parameters. The goal is to provide an optimal solution set, i.e., a set containing an optimal solution for each non-parametric problem induced by a particular choice of parameter values. In this paper, we consider right-hand-side single-parametric linear programming problems (RHS-PLPPs) that are usually very difficult to solve exactly. An optimal solution set of minimum cardinality can be infinite, while a representation by optimal bases may still require super-polynomially many different bases. Hence, we develop approximation methods for RHS-PLPPs that compute approximation sets of polynomially bounded cardinality, each containing, for every parameter value, a feasible solution whose objective value is at most a prescribed multiplicative factor times the optimum. Our first approach constructs an approximation set using a greedy staircase procedure and yields an FPTAS under general sufficient conditions. We show that these conditions hold for some important classes of RHS-PLPPs. For monotone instances, in which the feasible sets are nested as the parameter varies, our algorithm yields an approximation set of minimum cardinality. Beyond monotonicity, we establish near-minimum cardinality guarantees under more general sufficient conditions.
Finally, we introduce a dual approximation algorithm that returns solutions with objective value no greater than the optimum while allowing a controlled multiplicative violation of the parameter-dependent right-hand sides.

\medskip


\medskip


\medskip

\noindent
\textbf{Funding:} \,This research was funded by the Deutsche Forschungsgemeinschaft (DFG, German Research Foundation) -- Project number 508981269, and GRK 2982, 516090167 ``Mathematics of Interdisciplinary Multiobjective Optimization''.

\medskip
\end{abstract}

    \section{Introduction}

    Parametric optimization provides a natural framework for studying optimization problems under changing input data. In a parametric optimization problem, the objective function, the feasible set, or both depend on one or several unknown parameters. Instead of solving a single optimization problem for fixed data, the aim is to determine a set of solutions that contains an optimal solution for every possible choice of parameter values. Parametric optimization has found applications across a range of complex areas, including biomedical systems and power electronics, as well as fleet planning, aeronautics, and process synthesis under uncertainty~\cite{pistikopoulos_theoretical_2012,oberdieck.diangelakis.ea2016MultiparametricProgrammingIts,Jenkins_1987_Using_Parametric}.
    The theoretical foundations of parametric optimization and postoptimality analysis have been developed over several decades; see, e.g., \citeauthor{dinkelbach_sensitivitatsanalysen_1969} \cite{dinkelbach_sensitivitatsanalysen_1969}, \citeauthor{gal_multiparametric_1972} \cite{gal_multiparametric_1972}, \citeauthor{gal_degeneracy_1995} \cite{gal_degeneracy_1995}, and \citeauthor{gal_advances_1997} \cite{gal_advances_1997}. More recent surveys provide comprehensive overviews of the theory and algorithms for solving (multi-)parametric optimization problems \cite{pistikopoulos_theoretical_2012,nemesch_survey_2025}.

    In the following, we focus on single-parametric problems in which the parameter ranges over an interval $\Lambda \coloneq [\lambda_{\min}, \lambda_{\max}]$.
    The behavior of such a problem is largely characterized by its \emph{optimal cost curve}, which maps each parameter value~$\lambda \in \Lambda$ to the optimal objective value of the corresponding non-parametric problem. The points at which the slope of this curve changes are known as \emph{breakpoints}.
    Within this setting, we consider right-hand-side single-parametric linear programming problems (RHS-PLPPs), in which the parameter appears exclusively on the right-hand side of the constraints. The corresponding dual linear programming problem is a parametric linear programming problem in which the parameter appears in the objective function (OBF-PLPP). Strong duality establishes a close structural relationship between each RHS-PLPP and its dual OBF-PLPP. In particular, the two problems have the same continuous, piecewise-linear, and convex optimal cost curve~\cite{geoffrion_parametric_1977}.

    However, a fundamental distinction arises regarding the representation of optimal solution sets. For an OBF-PLPP, a fixed solution induces an affine objective-value function in the parameter and may therefore remain optimal over an entire parameter interval, so finitely many explicit solutions suffice to represent an optimal solution set. In an RHS-PLPP, by contrast, the objective value of a fixed solution does not change if the parameter value changes. Consequently, whenever the optimal cost curve contains a strictly increasing or strictly decreasing segment -- which is typically the case -- every optimal solution set contains infinitely many distinct solutions. To obtain a finite representation, optimal solutions of RHS-PLPPs are therefore typically represented by finitely many optimal bases together with the parameter regions on which these bases remain optimal, rather than by explicit solutions \cite{gal_multiparametric_1972,ward_approaches_1990,borrelli_geometric_2003}.
    
    Since an optimal basis of an RHS-PLPP and an optimal solution of an OBF-PLPP each induce an affine objective-value function in the parameter, a breakpoint of the optimal cost curve requires a change of the optimal basis or solution, respectively. Hence, the number of breakpoints provides an important measure of the structural complexity of such problems.
    Unfortunately, many standard problems can have super-polynomially many breakpoints in the worst case --- including the parametric variants of the shortest path, minimum cut, and assignment problems \cite{Carstensen_1983_The_Complexity,Nikolova_2006_Stochastic_Shortest,Gassner_2010_A_fast}. 
    This means that, even though the corresponding non-parametric problems are linear programming problems that can be solved efficiently in polynomial time, the exact solution of the corresponding parametric linear programming problem may have super-polynomial complexity~\cite{murty_computational_1980}.

    To overcome this intractability, approximation algorithms for parametric optimization have been studied in different settings.
    Approximation provides a concept to substantially reduce the number of required solutions while still obtaining provable solution quality. For minimization problems and a given multiplicative approximation factor~$\alpha>1$, the goal is to find an $\alpha$-approximation set, i.e., a set of solutions such that, for every parameter value~$\lambda \in \Lambda$, it contains a solution that is feasible for the corresponding non-parametric problem and whose objective value is at most $\alpha$ times its optimum.

    For general parametric optimization problems in which the parameter appears in the objective function, several approximation frameworks have been proposed, comprehensive summaries of which can be found in recent surveys (e.g., \citeauthor{nemesch_survey_2025} \cite{nemesch_survey_2025}). Approximation schemes have also been developed for other forms of parameter dependence. For instance, \citeauthor{halman_fptas_2018} \cite{halman_fptas_2018} constructed a fully polynomial-time approximation scheme (FPTAS) for the parametric knapsack problem, in which the parameter appears in the item weights on the left-hand side of the constraints. Related parametric knapsack problems with max--max, max--min, min--max, and min--min objectives are considered by \citeauthor{halman_knapsack_2023} \cite{halman_knapsack_2023}.    

    By contrast, the literature on approximating parametric problems with parameter-de\-pendent feasible sets is limited. For the right-hand-side parametric linear programming problems considered here, existing approaches predominantly rely on additive error guarantees (see, e.g., \citeauthor{filippi_algorithm_2004}~\cite{filippi_algorithm_2004}). While the basic algorithm of \citeauthor{filippi_algorithm_2004} is polynomial in the size of its output, no polynomial bound on the output size in terms of the input size and the prescribed error tolerance is established.
    Multiplicative approximation frameworks for RHS-PLPPs therefore remain largely unexplored. In this paper, we bridge this gap by establishing new approximation algorithms for RHS-PLPPs.

    \subsection{Our Contribution}

After establishing fundamental structural properties of RHS-PLPPs and analyzing their inherent intractability, we present two approximation approaches. To the best of our knowledge, these provide the first multiplicative approximation results for RHS-PLPPs.

Our first approach uses a greedy staircase construction to compute an $\alpha$-approximation set. Starting from a parameter value where the optimal cost curve attains its global minimum, the algorithm iteratively expands the approximated parameter interval in both directions. At each step, following an approach similar to that of \citeauthor{Diakonikolas2008} \cite{Diakonikolas2008}, an $\alpha$-covering oracle greedily adds a new step to the staircase approximation by returning a set of feasible solutions that provides an $\alpha$-approximation over the largest possible next parameter interval to the left or right. This systematically reduces the global parametric approximation problem to a sequence of local covering subproblems.

We show that, under some basic assumptions, this framework yields an FPTAS whenever such a po\-ly\-no\-mi\-al-time $\alpha$-covering oracle is available. We then identify classes of RHS-PLPPs for which the required covering steps can be performed efficiently. For monotone problems, in which the feasible sets are nested as the parameter varies, every covering-oracle call can be answered by a single optimal solution, yielding an FPTAS. Moreover, the greedy construction produces an $\alpha$-approximation set of minimum cardinality. We further consider non-monotone RHS-PLPPs. We use covering sets, which record solutions together with parameter intervals on which they provide the required approximation guarantee, with these intervals covering the parameter range under consideration. We provide a polynomial-time procedure that answers a covering-oracle call by a singleton covering set whenever this is possible. This yields an FPTAS for classes on which the procedure succeeds for every required call and every rational input approximation factor greater than one. Whenever all calls in a given execution succeed, the returned covering set contains at most one more element than a minimum-cardinality covering set. We also derive a sufficient condition based on the minimal duration of the instance that guarantees successful calls and show how this condition can be checked explicitly for a class of moving-band RHS-PLPPs.

Finally, we introduce a second approximation approach based on a dual notion of approximation. Rather than allowing a multiplicative loss in the objective value, this approach returns solutions whose objective values are no greater than the corresponding optima while allowing controlled multiplicative violations of the parameter-dependent right-hand sides of the constraints. We show that primal and dual approximation sets are equivalent up to a uniform scaling of the solutions and that this correspondence preserves minimum cardinality as well as the existence of an FPTAS. Under a positivity condition on the right-hand sides, we obtain a dual FPTAS without requiring a covering oracle.

	\section{Preliminaries}
	We distinguish between the natural numbers~$\mathbb{N}$ starting with~$1$ and~$\mathbb{N}_0$ starting with~$0$. 
	For $p\in\mathbb{N}_0$, let $[p]\coloneqq\{i\in\mathbb{N}:1\leq i\leq p\}$; in particular, $[0]=\emptyset$.
	For $p\in\mathbb{N}$, we denote the nonnegative orthant by $\mathbb{R}^p_{\geq 0} \coloneq \{y \in \mathbb{R}^p: \, y \geq 0 \}$, where $0 \in \mathbb{R}^p$ is the $p$-dimensional zero vector and $\geq$ is the weak component-wise order defined by
	\begin{align*}
		y \geq y' \text{ if and only if } y_i \geq y_i' \quad \forall i \in [p].
	\end{align*}
	For $p,n\in\mathbb{N}$, given a function $f\colon\mathbb{R}^p\rightarrow \mathbb{R}^n$ and a set $Z \subseteq \mathbb{R}^p$, we denote the image of~$Z$ under~$f$ by $f(Z)\coloneq \{f(z): \, z \in Z\} \subseteq \mathbb{R}^n$.
	
	\medskip
	
	For running time discussions, we use the well-known Landau notation.
	For a problem~$\mathsf{P}$, we let~$T_{\mathsf{P}}$ denote the running time of an oracle for~$\mathsf{P}$.
    We write $\operatorname{poly}(n)$ for an unspecified function bounded above by a polynomial in~$n$.
	An algorithm is \emph{(input-)polynomial} if, for every input of size~$I$, its running time is in~$\mathcal{O}(\operatorname{poly}(I))$.
	A problem is called \emph{intractable} if it cannot be solved by a polynomial-time algorithm. 
    An algorithm is \emph{output-polynomial} if
    its running time is polynomially bounded in
    the size of the input and the size of the
    output it produces; that is, if there exists a bivariate polynomial~$q$ such that, for every input of size~$I$ producing an output of size~$O$, its running time is in $\mathcal{O}(q(I,O))$.
	
	\medskip
	
	Next, we introduce the relevant concepts from \emph{linear programming}. For a more detailed introduction, we refer to \citeauthor{luenberger_linear_2008} \cite{luenberger_linear_2008}.
    For $m,n\in\mathbb{N}$, let $A\in\mathbb{Q}^{m\times n}$, $c\in\mathbb{Q}^{n}$, and $b\in\mathbb{Q}^m$. Then, a minimization linear programming problem (LPP) is a problem of the form
    \begin{align*}
    \begin{alignedat}{2}
        &\min \quad F(x) \;&\coloneq\;& c^\mathsf{T}x \\
        &\,\,\,\text{s.t.}\quad x \in X \;&\coloneq\;&
        \left\{x \in \mathbb{R}_{\geq 0}^n |\, Ax \geq b \right\}.
    \end{alignedat}
    \end{align*}
    Its corresponding dual problem is the maximization linear programming problem
    \begin{align*}
    \begin{alignedat}{2}
        &\max \quad G(y) \;&\coloneq\;& b^\mathsf{T}y \\
        &\,\,\,\text{s.t.}\quad y \in Y \;&\coloneq\;&
        \left\{y \in \mathbb{R}_{\geq 0}^m \mid A^\mathsf{T}y \leq c \right\}.
    \end{alignedat}
    \end{align*}
	As shorthands, $F\colon X\to \mathbb{R}$ and $G\colon Y\to \mathbb{R}$ denote the \emph{objective functions}, and $X\subseteq\mathbb{R}^n$ and $Y\subseteq\mathbb{R}^m$ denote the sets of \emph{feasible solutions}, respectively.
	For feasible solutions $x \in X$ and $y \in Y$, we call the points~$F(x)$ and $G(y)$ their \emph{objective values}, and denote the sets of all objective values by $F(X)$ and $G(Y)$.
	The \emph{strong duality theorem} (see, e.g., \cite{luenberger_linear_2008}) states that if the primal LP has an optimal solution~$x^*$, then the dual LP has an optimal solution~$y^*$ with $F(x^*)=G(y^*)$. 
	
	\subsection{Parametric Linear Programming Problems}\label{sec:defRHS}
    We consider single-parametric linear programming problems in which a real-valued parameter affects either the objective function or the right-hand side of the constraints. In the former case, we speak of a parametric linear programming problem with a parameter-dependent objective function (OBF-PLPP), whereas in the latter case, we speak of a right-hand-side parametric linear programming problem (RHS-PLPP).
    Throughout the paper, the parameter is assumed to range over a compact interval. More precisely, we use the following definition:
	
	\begin{definition}[Parameter Set]\label{def:ParamSet}
		A \emph{parameter set~$\Lambda\subseteq\mathbb{R}$} is a set of the form 
        \[\Lambda \coloneq \{\lambda \in \mathbb{R}: \lambda_{\min} \leq \lambda \leq \lambda_{\max}\}\] 
        for $\lambda_{\min},\lambda_{\max}\in\mathbb{Q}$ with $\lambda_{\min}\leq\lambda_{\max}$.
	\end{definition}

    More general parameter sets, such as non-convex sets, can also be considered in parametric optimization. However, following the standard setting in the parametric optimization literature, we consider only parameter sets of the form defined above. Our definition of parametric optimization is based on the definition of an OBF-PLPP given by \citeauthor{bazgan_approximation_2022} \cite{bazgan_approximation_2022}.
	
	\begin{definition}[RHS- and OBF-Parametric Linear Programming Problems]\label{Def:OBF-/RHS-ParamOptDef}
		For $n,m \in \mathbb{N}$, let $A \in \mathbb{Q}^{m\times n}$, $c \in \mathbb{Q}^n$, and $b,r \in \mathbb{Q}^m$.
		\begin{itemize}
			\item[1.] An instance of an \emph{\mbox{RHS-(single-)parametric} minimization linear programming problem} is a triple $\Pi_{\text{RHS}}=(\Lambda, X, F)$ consisting of 
			\begin{enumerate}[label=(\alph*)]
				\item a parameter set $\Lambda\subseteq\mathbb{R}$ as in Definition \ref{def:ParamSet},
				\item sets~$X_{\lambda} \coloneq \{x \in \mathbb{R}^n: \, Ax \geq b + \lambda r,\, x\geq 0\} \subseteq\mathbb{R}^n$ of feasible solutions for every~$\lambda \in \Lambda$, collectively referred to as $X \coloneq \bigcup_{\lambda \in \Lambda} X_{\lambda}$ and
				\item a function $F\colon\mathbb{R}^n\to\mathbb{R}$, $x\mapsto F(x)=c^\T x$.
			\end{enumerate}
			\item[2.] An instance of the corresponding dual parametric LP, which is an \emph{\mbox{OBF-(single-)}para\-metric maximization linear programming problem}, is a triple $\Pi_{\text{OBF}}=(\Lambda,Y, g)$, consisting of
			\begin{enumerate}[label=(\alph*)]
				\item  the parameter set $\Lambda\subseteq\mathbb{R}$,
				\item the corresponding set~$Y \coloneq \{y \in \mathbb{R}^m: \, A^\mathsf{T} y \leq c,\, y \geq 0\} \subseteq\mathbb{R}^m$ of feasible solutions, and
				\item the function $g\colon Y \to \mathbb{R}^{2}, \, y \mapsto g(y) \coloneq (g_0(y),g_1(y))^\mathsf{T} \coloneq (b^\mathsf{T}y,r^\mathsf{T}y)^\mathsf{T}$.
			\end{enumerate}
		\end{itemize}
	\end{definition}

    For each fixed parameter value~$\lambda\in\Lambda$, let
    $G_\lambda\colon Y\to\mathbb{R}$ be defined by
    $G_\lambda(y)\coloneq g_0(y)+\lambda\cdot g_1(y)$.
    The corresponding conventional non-parametric linear programming problems are
    \begin{align*}
    \min_{x\in X_\lambda}\;F(x)
    &\tag{$\Pi_{\text{RHS}}(\lambda)$}\\
    \max_{y\in Y}\;G_\lambda(y)
    &\tag{$\Pi_{\text{OBF}}(\lambda)$}.
    \end{align*}

    We use the following notion of optimal solution sets for the associated parametric instances.

    \begin{definition}[Optimal Solution Sets]
    \label{Def:OptimalSolutionSets}
    A set $S_{\text{RHS}}\subseteq X$ is called an \emph{optimal solution set for~$\Pi_{\text{RHS}}$} if, for every~$\lambda\in\Lambda$, it contains an optimal solution of~$\Pi_{\text{RHS}}(\lambda)$. Analogously, a set $S_{\text{OBF}}\subseteq Y$ is called an \emph{optimal solution set for~$\Pi_{\text{OBF}}$} if, for every~$\lambda\in\Lambda$, it contains an optimal solution of~$\Pi_{\text{OBF}}(\lambda)$. When no ambiguity arises, we refer to an optimal solution set simply as a solution set. The cardinality of a solution set~$S$ is denoted by~$|S|$.
\end{definition}

    Let $\overline{\mathbb{R}}\coloneq\mathbb{R}\cup\{-\infty,+\infty\}$. The optimal cost curve 
    \[F^*\colon\Lambda\to\overline{\mathbb{R}}, \lambda\mapsto F^*(\lambda)\coloneq\inf_{x\in X_\lambda}F(x),\]
    assigns to every~$\lambda\in\Lambda$ the extended optimal objective value of~$\Pi_{\text{RHS}}(\lambda)$, where $+\infty$ represents infeasibility and $-\infty$ unboundedness below. Analogously, the optimal cost curve 
    \[G^*\colon\Lambda\to\overline{\mathbb{R}},\lambda\mapsto G^*(\lambda)\coloneq\sup_{y\in Y}G_\lambda(y),\]
    assigns to every~$\lambda\in\Lambda$ the extended optimal objective value of~$\Pi_{\text{OBF}}(\lambda)$, where $-\infty$ represents infeasibility and $+\infty$ unboundedness above. The optimal cost curves~$F^*$ and~$G^*$ are also called the lower envelope of~$\Pi_{\mathrm{RHS}}$ and the upper envelope of~$\Pi_{\mathrm{OBF}}$, respectively.
    
	The \emph{critical region} $\Lambda(y)\subseteq\Lambda$ of a solution~$y\in Y$ is the set of all parameter values for which~$y$ is an optimal solution of~$\Pi_{\text{OBF}}(\lambda)$. For every~$y\in Y$, the critical region~$\Lambda(y)$ is either empty or a closed interval. A parameter value at which the slope of the optimal cost curve changes is called a \emph{breakpoint}.
    
    Furthermore, we define the non-parametric problem that finds the global optimum of the optimal cost curve of an RHS-PLPP.
	\begin{definition}[RHS-Parametric Search Problem]\label{Def:ParametricSearchProblem}
		Given an RHS-PLPP instance $\Pi_{\text{RHS}}=(\Lambda, X, F)$, the associated \emph{RHS-parametric search problem} is
		\begin{align}\tag{$\Pi_{\text{PS}}$}\label{Problem:ParametricSearch}
			\begin{array}{rl}
				\displaystyle\min_{x,\lambda}\, F(x) \quad \quad & \\
				\quad \text{s.t.}\quad Ax - \lambda r &\geq b \\
				\qquad \quad \,\, x &\geq 0 \\
				\qquad \quad \,\, \lambda &\in [\lambda_{\min},\lambda_{\max}].			
			\end{array}
		\end{align}
	\end{definition}

    Throughout the paper, unless explicitly stated otherwise, we restrict attention to RHS-PLPP instances satisfying the following assumption.

    \begin{assumption}
    \label{Assumption:Standing}
    The RHS-PLPP instance $\Pi_{\text{RHS}}$ satisfies the following conditions:
    \begin{enumerate}
    \item $X_{\lambda_{\min}}\neq\emptyset$ and $X_{\lambda_{\max}}\neq\emptyset$, and at least one of $\Pi_{\text{RHS}}(\lambda_{\min})$ and $\Pi_{\text{RHS}}(\lambda_{\max})$ has an optimal solution. \label{Assumption:Standing1}
    \item $\inf_{\lambda\in\Lambda}F^*(\lambda)>0$. 
    \end{enumerate}
    \end{assumption}
    
    As shown later in Proposition~\ref{Proposition:ExistenceOptimalSolutions} and Corollary~\ref{Cor:OptimalSolutionSet}, the first condition implies that $\Pi_{\text{RHS}}(\lambda)$ has an optimal solution for every~$\lambda\in\Lambda$ and, consequently, that $\Pi_{\text{RHS}}$ has an optimal solution set. Together with the first condition, the second condition ensures that the multiplicative approximation framework introduced in Section~\ref{Sec:ParametricApproximation} is well-defined.

    \medskip

    We use the following optimization oracles throughout the paper.

    \begin{definition}[Oracles]\label{Definition:OptimizationOracles}
    We let $\mathrm{ALG}_{\text{RHS}}$ and $\mathrm{ALG}_{\text{OBF}}$ denote oracles that take as input a parameter value~$\lambda\in\Lambda$ and return an optimal solution for $\Pi_{\text{RHS}}(\lambda)$ and $\Pi_{\text{OBF}}(\lambda)$, respectively. Moreover, an \emph{RHS-parametric search oracle} $\mathrm{ALG}_{\text{PS}}$ takes an RHS-PLPP instance $\Pi_{\text{RHS}}$ as input and returns an optimal solution $(x_{\text{PS}},\lambda_{\text{PS}})$ for~$\Pi_{\text{PS}}$.
    \end{definition}
    
    As shown later, the first condition in Assumption~\ref{Assumption:Standing} guarantees the existence of the optimal solutions returned by all oracles in Definition~\ref{Definition:OptimizationOracles}, which are therefore well-defined. The second condition additionally ensures that the optimal value of~$\Pi_{\text{PS}}$ is positive.
    
    For every parameter value~$\lambda\in\Lambda$ whose encoding length is polynomially bounded in that of~$\Pi_{\text{RHS}}$, the problems $\Pi_{\text{RHS}}(\lambda)$ and $\Pi_{\text{OBF}}(\lambda)$ are linear programming problems whose encoding lengths are polynomially bounded in that of~$\Pi_{\text{RHS}}$. The parametric search problem~$\Pi_{\text{PS}}$ is likewise a linear programming problem whose encoding length is polynomially bounded in that of~$\Pi_{\text{RHS}}$. Hence, if the corresponding oracles are implemented using a polynomial-time linear programming algorithm and the pointwise oracles are queried only at parameter values of polynomially bounded encoding length, their running times are polynomial in the encoding length of~$\Pi_{\text{RHS}}$ \cite{schrijver_theory_1998}.
   
    \subsection{Elementary Properties and Intractability}
	
	In this section, we provide an overview of some important properties of RHS-PLPPs, which are used later in the construction of our approximation algorithms and in the applications. 
	From now on, if not stated otherwise, we assume that an RHS-PLPP instance $\Pi_{\text{RHS}}=(\Lambda, X, F)$ and its corresponding dual OBF-PLPP instance $\Pi_{\text{OBF}}=(\Lambda,Y, g)$ are given.
	
    \begin{lemma}[Feasibility Property]\label{Lemma:FeasibilityProp}
    If $X_{\lambda_1}$ and $X_{\lambda_2}$ are nonempty for some $\lambda_1,\lambda_2\in\Lambda$ with $\lambda_1\leq\lambda_2$, then $X_\lambda$ is nonempty for every $\lambda\in[\lambda_1,\lambda_2]$.
    \end{lemma}
	\begin{proof}
        Fix an arbitrary~$\lambda \in [\lambda_1,\lambda_2]$ and let $x_1\in X_{\lambda_1}$ and $x_2\in X_{\lambda_2}$ be feasible solutions for~$\lambda_1$ and~$\lambda_2$, respectively. We can then write~$\lambda$ as a convex combination 
        \[\lambda=\mu \lambda_1 + (1-\mu) \lambda_2\quad\text{ with }\mu\in [0,1].\]
        The corresponding convex combination $\mu x_1 + (1-\mu) x_2$ of~$x_1$ and~$x_2$ then satisfies
		\begin{align*}
				A(\mu x_1 + (1-\mu) x_2) 
                &= \mu Ax_1 + (1-\mu) Ax_2 \\
                &\geq \mu (b+\lambda_1 r) + (1-\mu) (b+\lambda_2 r) \\
				&= b + (\mu \lambda_1 + (1-\mu) \lambda_2)r \\
                &= b + \lambda r,
		\end{align*}
    and $\mu x_1 + (1-\mu) x_2 \geq 0$, i.e., $\mu x_1 + (1-\mu) x_2 \in X_{\lambda}$. 
	\end{proof}

    The following proposition provides simple conditions under which feasibility and the existence of optimal solutions are guaranteed for the entire parameter interval.

    \begin{proposition}[Existence of Optimal Solutions]\label{Proposition:ExistenceOptimalSolutions}
    Suppose that $X_{\lambda_{\min}}$ and $X_{\lambda_{\max}}$ are nonempty. Then, $X_\lambda$ is nonempty for every~$\lambda\in\Lambda$. Moreover, if $\Pi_{\text{RHS}}(\bar\lambda)$ is bounded below for some~$\bar\lambda\in\Lambda$, then $\Pi_{\text{RHS}}(\lambda)$ has an optimal solution for every~$\lambda\in\Lambda$.
    \end{proposition}

    \begin{proof}
    The nonemptiness of~$X_\lambda$ for every~$\lambda\in\Lambda$ follows from Lemma~\ref{Lemma:FeasibilityProp}. Moreover, for every $\lambda\in\Lambda$, the recession cone of $X_\lambda$ is
    \[
    \operatorname{rec}(X_\lambda)
    =
    \{d\in\mathbb{R}^n_{\geq 0}:Ad\geq 0\},
    \]
    and is therefore independent of $\lambda$.
    Now, let $\bar\lambda\in\Lambda$ be such that $\Pi_{\text{RHS}}(\bar\lambda)$ is bounded below. 
    Suppose, for contradiction, that $\Pi_{\text{RHS}}(\lambda)$ is unbounded below for some~$\lambda\in\Lambda$. Since $X_\lambda$ is a nonempty polyhedron, this implies the existence of a direction $d\in\operatorname{rec}(X_\lambda)$ such that $c^\T d<0$. Since the recession cone is independent of $\lambda$, we also have $d\in\operatorname{rec}(X_{\bar\lambda})$. Thus, for any $\bar x\in X_{\bar\lambda}$ and $t\geq 0$, we have $\bar x+td\in X_{\bar\lambda}$ and
    \begin{align*}
        c^\T(\bar x+td)=c^\T\bar x+t\,c^\T d\to-\infty
    \qquad\text{as }t\to\infty,
    \end{align*}
    contradicting the assumption that $\Pi_{\text{RHS}}(\bar\lambda)$ is bounded below. Hence, $\Pi_{\text{RHS}}(\lambda)$ is bounded below for every $\lambda\in\Lambda$. Since each $X_\lambda$ is a nonempty polyhedron, the minimum is attained, so $\Pi_{\text{RHS}}(\lambda)$ has an optimal solution for every~$\lambda\in\Lambda$ as claimed.
    \end{proof}

    Since an optimal solution exists for every parameter value under the conditions of Proposition~\ref{Proposition:ExistenceOptimalSolutions}, an optimal solution set can be obtained by selecting one such solution for each parameter value.

    \begin{corollary}[Existence of an Optimal Solution Set]\label{Cor:OptimalSolutionSet}
    Suppose that $X_{\lambda_{\min}}$ and $X_{\lambda_{\max}}$ are nonempty and that either $\Pi_{\text{RHS}}(\lambda_{\min})$ or $\Pi_{\text{RHS}}(\lambda_{\max})$ has an optimal solution. Then, $\Pi_{\text{RHS}}$ has an optimal solution set.
    \end{corollary}

    Consequently, the conditions in Assumption~\ref{Assumption:Standing} guarantee the existence of an optimal solution set for $\Pi_{\mathrm{RHS}}$.
    The following lemma characterizes the set of parameter values for which a given solution~$x\in X$ is feasible.

	\begin{lemma}[Duration Property]\label{Lemma:DurationProp}
		For any~$x \in X$, the set $\Lambda(x) \subseteq \Lambda$ of parameter values for which~$x$ is feasible is a closed interval, i.e., $\Lambda(x)=[\lambda^{\mathrm{LB}},\lambda^{\mathrm{UB}}]$ for some $\lambda^{\mathrm{LB}}, \lambda^{\mathrm{UB}}\in\Lambda$ with $\lambda^{\mathrm{LB}}\leq \lambda^{\mathrm{UB}}$.
	\end{lemma}
    \begin{proof}
    Fix any solution $x\in X$. By the definition of~$X$, there exists some $\bar\lambda\in\Lambda$ such that $x\in X_{\bar\lambda}$. For each $i\in[m]$, define
    \[
    \Lambda_i(x)\coloneq
    \left\{\lambda\in\Lambda:
    a_i^\T x\geq b_i+\lambda r_i
    \right\}.
    \]
    Then, $\Lambda(x)=\bigcap_{i=1}^m\Lambda_i(x)$.
    We show that each $\Lambda_i(x)$ is a nonempty closed interval. If $r_i=0$, then the $i$-th constraint does not depend on~$\lambda$. Since it is satisfied at $\bar\lambda$, it is thus satisfied for every $\lambda\in\Lambda$, and hence $\Lambda_i(x)=\Lambda$.
    If $r_i>0$, the constraint is equivalent to $\lambda\leq\frac{a_i^\T x-b_i}{r_i}$,
    whereas, if $r_i<0$, it is equivalent to $\lambda\geq\frac{a_i^\T x-b_i}{r_i}$.
    Thus, in either case, $\Lambda_i(x)$ is the intersection of~$\Lambda$ with a closed half-line and is therefore a closed interval. Moreover, $\bar\lambda\in\Lambda_i(x)$, so $\Lambda_i(x)$ is nonempty.
    Consequently, $\Lambda(x)$ is the intersection of finitely many closed intervals and is itself a closed interval. Since $\bar\lambda\in\Lambda(x)$, this interval is nonempty.
    \end{proof}

    We refer to the length of the parameter interval over which a solution remains feasible as its \emph{duration}. This notion is relevant for constructing approximation sets, since instances in which solutions remain feasible only over short parameter intervals may require larger approximation sets.

    \begin{definition}[Duration]\label{Def::Duration}
        The length $d^x\coloneq\lambda^{\mathrm{UB}}-\lambda^{\mathrm{LB}}\in [0,\lambda_{\max}-\lambda_{\min}]$ of the parameter subinterval $\Lambda(x)=[\lambda^{\mathrm{LB}},\lambda^{\mathrm{UB}}] \subseteq \Lambda$ on which a solution~$x \in X$ is feasible is called the \emph{duration} of~$x$. 
        The \emph{maximal duration} of an instance~$\Pi_{\text{RHS}}$ is defined as $d_{\max}^{\Pi_{\text{RHS}}} \coloneq \max_{x\in X} d^x$ and the \emph{minimal duration} as $d_{\min}^{\Pi_{\text{RHS}}} \coloneq \min_{x\in X} d^x$.
    \end{definition}

    For \emph{monotone} RHS-PLPPs, that is, for $r\geq 0$ or $r\leq 0$, the feasible sets are nested as the parameter varies~\cite{geoffrion_parametric_1977}.
	
    \begin{lemma}[Monotonicity Property]\label{Lemma:MonotonicityProp}\mbox{}
	\begin{enumerate}
		\item If $r \geq 0$, then $X_{\lambda_2} \subseteq X_{\lambda_1}$ for all $\lambda_1,\lambda_2 \in \Lambda$ with $\lambda_1 \leq \lambda_2$, and $F^*$ is non-decreasing.
		\item If $r \leq 0$, then $X_{\lambda_1} \subseteq X_{\lambda_2}$ for all $\lambda_1,\lambda_2 \in \Lambda$ with $\lambda_1 \leq \lambda_2$, and $F^*$ is non-increasing.
	\end{enumerate}
	\end{lemma}
    \begin{proof}
    Assume first that $r\geq 0$ and let
    $x\in X_{\lambda_2}$. Then, $Ax\geq b+\lambda_2r\geq b+\lambda_1r$, and hence $x\in X_{\lambda_1}$. Thus, $X_{\lambda_2}\subseteq X_{\lambda_1}$. Since the objective function is independent of $\lambda$ and the problem is a minimization problem, this inclusion implies that
    $F^*(\lambda_1)\leq F^*(\lambda_2)$.
    
    If $r\leq 0$, then $b+\lambda_1r\geq b+\lambda_2r$, so
    the same argument yields $X_{\lambda_1}\subseteq X_{\lambda_2}$ and therefore
    $F^*(\lambda_1)\geq F^*(\lambda_2)$.
    \end{proof}
	
    We next derive several basic properties of RHS-PLPPs from their duality relationship with the corresponding OBF-PLPPs.
	
	\begin{lemma}[Strong Duality between RHS- and OBF-PLPPs]\label{Lemma:Duality}
		If either the primal RHS-PLPP instance $\Pi_{\text{RHS}}$ or its corresponding dual OBF-PLPP instance $\Pi_{\text{OBF}}$ has an optimal solution set, then so does the other, and their optimal cost curves coincide, i.e., $F^*=G^*$.
	\end{lemma}
    \begin{proof}
    For each $\lambda\in\Lambda$, apply the strong duality theorem to the primal--dual pair $\Pi_{\text{RHS}}(\lambda)$ and $\Pi_{\text{OBF}}(\lambda)$. This shows that one of the two instances has an optimal solution if and only if the other has one, and that their optimal objective values coincide in this case. Hence, $\Pi_{\text{RHS}}$ has an optimal solution set if and only if $\Pi_{\text{OBF}}$ does, and their optimal cost curves then satisfy $F^*=G^*$.
    \end{proof}
	
    The preceding strong-duality result allows us to transfer known properties of the optimal cost curve of an OBF-PLPP, established by \citeauthor{dinkelbach_sensitivitatsanalysen_1969} \cite{dinkelbach_sensitivitatsanalysen_1969}, to its dual RHS-PLPP.

	\begin{corollary}[Optimal Cost Curve]\label{Cor:OptCostCurve}
		The optimal cost curve~$F^*$ of an RHS-PLPP instance~$\Pi_{\text{RHS}}$ is continuous, piecewise linear, and convex.
	\end{corollary}
    \begin{proof}
    By Corollary~\ref{Cor:OptimalSolutionSet}, $\Pi_{\text{RHS}}$ has an optimal solution set. Hence, by Lemma~\ref{Lemma:Duality}, its optimal cost curve coincides with that of the corresponding dual OBF-PLPP, i.e., $F^*=G^*$. The claim therefore follows from the corresponding properties of~$G^*$ established in~\cite{dinkelbach_sensitivitatsanalysen_1969}.
    \end{proof}

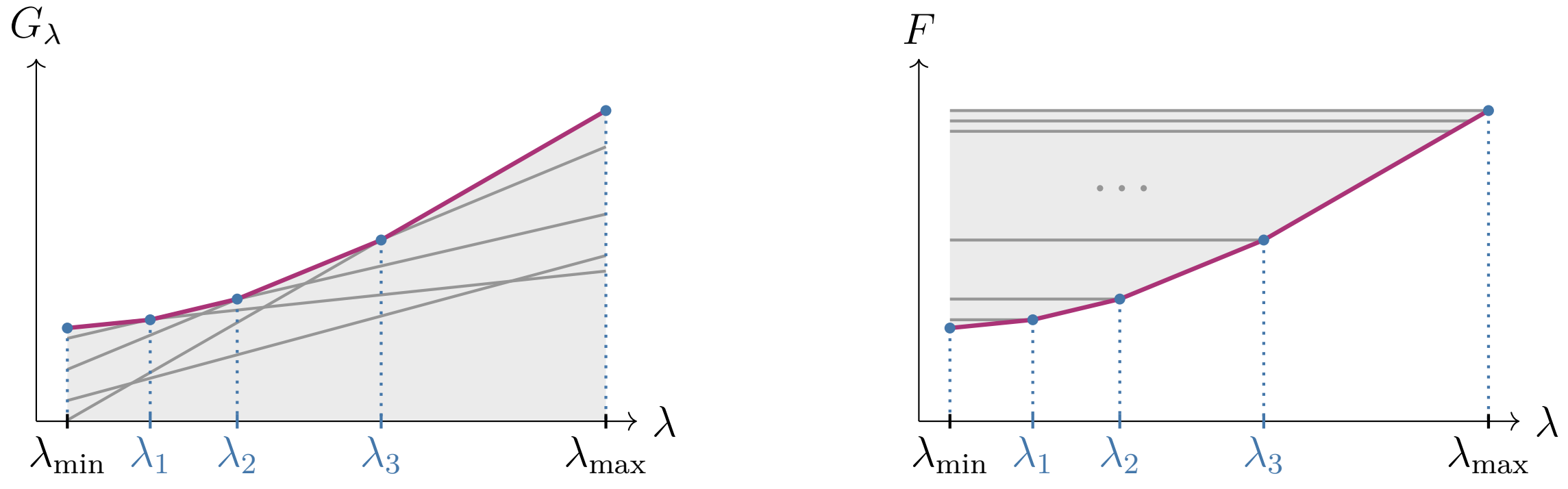
\begin{figure}[h]
		\centering
		\begin{minipage}[t]{.45\linewidth}
	\begin{tikzpicture}[scale=1]				
		\draw[->] (0,0) -- (5.8,0) node[right] {$\lambda$};
		\draw[->] (0,0) -- (0,3.5) node[above] {$G_\lambda$};
		
		\draw [brightgrey!80!black,thick] (0.3,0.9)--(5.5,1.45);
		\draw [brightgrey!80!black,thick] (0.3,0.8)--(5.5,2);
		\draw [brightgrey!80!black,thick] (0.3,0.5)--(5.5,2.65);
		\draw [brightgrey!80!black,thick] (0.3,0.01)--(5.5,3);
		\draw [brightgrey!80!black,thick] (0.3,0.2)--(5.5,1.6);
		
		\draw[very thick,brightpurple] (0.3,0.9) -- (1.1,0.98);
		\draw[very thick,brightpurple] (1.1,0.98) -- (1.94,1.18);
		\draw[very thick,brightpurple] (1.94,1.18) -- (3.33,1.75);
		\draw[very thick,brightpurple] (3.33,1.75) -- (5.5,3);
        \fill [brightblue](0.3,0.9) circle(1.5pt);
        \fill [brightblue](1.1,0.98) circle(1.5pt);
        \fill [brightblue](1.94,1.18) circle(1.5pt);
        \fill [brightblue](3.33,1.75) circle(1.5pt);
		\fill [brightblue](5.5,3) circle(1.5pt);

		\draw [brightblue,dotted,thick] (0.3,0.9)--(0.3,0);
		\draw [brightblue,dotted,thick] (1.1,0.98)--(1.1,0);
		\draw [brightblue,dotted,thick] (1.94,1.18)--(1.94,0);
		\draw [brightblue,dotted,thick] (3.33,1.75)--(3.33,0);
		\draw [brightblue,dotted,thick] (5.5,3)--(5.5,0);
		
		\draw [brightblue,thick] (1.1,0.07)--(1.1,-0.07);
		\draw [brightblue,thick] (1.94,0.07)--(1.94,-0.07);
		\draw [brightblue,thick] (3.33,0.07)--(3.33,-0.07);
		\draw [thick] (5.5,0.07)--(5.5,-0.07);
		\draw [thick] (0.3,0.07)--(0.3,-0.07);
		
		\begin{scope}[on background layer]
			\fill[brightgrey!30] (0.3,0)--(0.3,0.9)--(1.1,0.98)--(1.94,1.18)--(3.33,1.75)--(5.5,3)--(5.5,0);
		\end{scope}
		
		\tkzText[below,brightblue](1.1,0){$\lambda_1$}
		\tkzText[below,brightblue](1.94,0){$\lambda_2$}
		\tkzText[below,brightblue](3.33,0){$\lambda_3$}
		\tkzText[below](5.5,0){$\lambda_{\text{max}}$}
		\tkzText[below](0.3,0){$\lambda_{\text{min}}$}
	\end{tikzpicture}
	\\(a) OBF-PLPP
\end{minipage}
\hfill
\begin{minipage}[t]{.45\linewidth}
	\begin{tikzpicture}[scale=1]
		\draw[->] (0,0) -- (5.8,0) node[right] {$\lambda$};
		\draw[->] (0,0) -- (0,3.5) node[above] {$F$};
		
		\draw [brightgrey!80!black,thick] (0.3,0.98)--(1.1,0.98);
		\draw [brightgrey!80!black,thick] (0.3,1.18)--(1.94,1.18);
		\draw [brightgrey!80!black,thick] (0.3,1.75)--(3.33,1.75);
		\draw [brightgrey!80!black,thick] (0.3,3)--(5.5,3);
		\draw [brightgrey!80!black,thick] (0.3,2.9)--(5.33,2.9);
		\draw [brightgrey!80!black,thick] (0.3,2.8)--(5.15,2.8);
	
		\node[text=brightgrey!80!black] at (2, 2.25) {\textbf{\dots}};

		\draw[very thick,brightpurple] (0.3,0.9) -- (1.1,0.98);
		\draw[very thick,brightpurple] (1.1,0.98) -- (1.94,1.18);
		\draw[very thick,brightpurple] (1.94,1.18) -- (3.33,1.75);
		\draw[very thick,brightpurple] (3.33,1.75) -- (5.5,3);
        \fill [brightblue](0.3,0.9) circle(1.5pt);
        \fill [brightblue](1.1,0.98) circle(1.5pt);
        \fill [brightblue](1.94,1.18) circle(1.5pt);
        \fill [brightblue](3.33,1.75) circle(1.5pt);
		\fill [brightblue](5.5,3) circle(1.5pt);
		
		\draw [brightblue,dotted,thick] (0.3,0.9)--(0.3,0);
		\draw [brightblue,dotted,thick] (1.1,0.98)--(1.1,0);
		\draw [brightblue,dotted,thick] (1.94,1.18)--(1.94,0);
		\draw [brightblue,dotted,thick] (3.33,1.75)--(3.33,0);
		\draw [brightblue,dotted,thick] (5.5,3)--(5.5,0);
		
		\draw [brightblue,thick] (1.1,0.07)--(1.1,-0.07);
		\draw [brightblue,thick] (1.94,0.07)--(1.94,-0.07);
		\draw [brightblue,thick] (3.33,0.07)--(3.33,-0.07);
		\draw [thick] (5.5,0.07)--(5.5,-0.07);
		\draw [thick] (0.3,0.07)--(0.3,-0.07);
		
		\begin{scope}[on background layer]
			\fill[brightgrey!30] (0.3,3)--(0.3,0.9)--(1.1,0.98)--(1.94,1.18)--(3.33,1.75)--(5.5,3)--(5.5,3);
		\end{scope}

		\tkzText[below,brightblue](1.1,0){$\lambda_1$}
		\tkzText[below,brightblue](1.94,0){$\lambda_2$}
		\tkzText[below,brightblue](3.33,0){$\lambda_3$}
		\tkzText[below](5.5,0){$\lambda_{\text{max}}$}
		\tkzText[below](0.3,0){$\lambda_{\text{min}}$}
	\end{tikzpicture}
	\\(b) RHS-PLPP
\end{minipage}
        \caption{Illustration of the relationship between an OBF-PLPP and its dual RHS-PLPP in the monotone case~$r\geq0$. In the OBF-PLPP, the gray lines show how the objective values~$G_\lambda(y)$ of fixed dual solutions vary with~$\lambda$. In the RHS-PLPP, the horizontal gray segments show the constant objective values~$F(x)$ of fixed primal solutions over the parameter intervals on which they are feasible. The common optimal cost curve~$F^*=G^*$ is shown in purple. The shaded regions represent~$G_\lambda(Y)$ and~$F(X_\lambda)$, and the marked parameter values are the breakpoints.}\label{fig:OverallOBFRHS}
	\end{figure}

    The preceding results provide a structural interpretation of the relationship between an RHS-PLPP and its dual OBF-PLPP. For an OBF-PLPP, every fixed solution~$y\in Y$ induces an affine function $\lambda\mapsto G_\lambda(y)$, and the optimal cost curve~$G^*$ is the upper envelope of these functions. By contrast, every fixed solution~$x\in X$ of an RHS-PLPP has the constant objective value~$F(x)$, but is feasible only on the interval~$\Lambda(x)$. Hence, $F^*$ is the lower envelope of the corresponding horizontal line segments. By Lemma~\ref{Lemma:Duality}, these two envelopes coincide.
    Figure~\ref{fig:OverallOBFRHS} illustrates this relationship for the monotone case~$r\geq0$. By Lemma~\ref{Lemma:MonotonicityProp}, the feasible sets~$X_\lambda$ of the RHS-PLPP are nested, with $X_{\lambda_2}\subseteq X_{\lambda_1}$ whenever $\lambda_1\leq\lambda_2$, and the optimal cost curve~$F^*$ is non-decreasing. Together with the continuity, piecewise linearity, and convexity established in Corollary~\ref{Cor:OptCostCurve}, this yields the structure shown in the figure. The case~$r\leq0$ is analogous, with the inclusions and the direction of monotonicity reversed.
	
    The RHS-PLPP panel in Figure~\ref{fig:OverallOBFRHS} also highlights a fundamental difference in the possible sizes of optimal solution sets. \citeauthor{murty_computational_1980}~\cite{murty_computational_1980} showed that some OBF-PLPP instances require super-polynomially many solutions in every optimal solution set. For RHS-PLPPs, the situation can be even more extreme: since the objective value of a fixed primal solution is independent of~$\lambda$, even a simple strictly increasing optimal cost curve may require infinitely many distinct solutions. The following instance demonstrates this.

    \begin{lemma}[Infinite Optimal Solution Sets]
    \label{Lemma:Intractability}
    There exists an RHS-PLPP instance for which every optimal solution set has infinite cardinality.
    \end{lemma}
    
    \begin{proof}
    Consider the RHS-PLPP instance~$\Pi_{\mathrm{Inf}}$ given by:
    \begin{align}\tag{$\Pi_{\mathrm{Inf}}$}
    \left\{
    \begin{array}{ll}
    \displaystyle \min_x & x\\
    \,\,\text{s.t.} & x\geq\lambda,\\
                & x\in\mathbb{R}_{\geq 0}
    \end{array}
    \right\}_{\lambda\in[1,2]}.
    \end{align}
    For every~$\lambda\in[1,2]$, the unique optimal solution of~$\Pi_{\mathrm{Inf}}(\lambda)$ is $x_\lambda^*=\lambda$. Moreover, $x_\lambda^*$ is optimal only for the parameter value~$\lambda$. Hence, every optimal solution set must contain $x_\lambda^*=\lambda$ for every~$\lambda\in[1,2]$ and therefore contains the entire interval~$[1,2]$. Consequently, every optimal solution set has infinite cardinality.
    \end{proof}

    The corresponding dual OBF-PLPP instance behaves very differently:
    \begin{align*}
    \left\{
    \begin{array}{ll}
    \displaystyle \max_y & \lambda y\\
    \,\,\text{s.t.} & y\leq 1,\\
                & y\in\mathbb{R}_{\geq 0}
    \end{array}
    \right\}_{\lambda\in[1,2]}.
    \end{align*}
    For every~$\lambda\in[1,2]$, its unique optimal solution is~$y^*=1$. Hence, the singleton~$\{1\}$ is an optimal solution set for the dual OBF-PLPP instance, whereas every optimal solution set for the primal RHS-PLPP instance $\Pi_{\mathrm{Inf}}$ is infinite.
    
\begin{figure}[h]
    \centering
    \begin{tikzpicture}[thick] 
	
	\draw[brightgrey!80!black,->] (-1,0) -- (6,0) node[right, font=\large] {$x$}; 
	
	\draw[brightgrey!80!black] (0, -0.1) -- (0, 0.1) node[brightgrey!80!black,below=5mm, font=\large] {0};
	\draw[brightgrey!80!black] (2, -0.1) -- (2, 0.1) node[brightgrey!80!black,below=5mm, font=\large] {1};
	\draw[brightgrey!80!black] (4, -0.1) -- (4, 0.1) node[brightgrey!80!black,below=5mm, font=\large] {2};
	
	\node[left] at (2, -0.3) {\textcolor{black}{$x_1^*$}}; 
	\node[left] at (4, -0.3) {\textcolor{black}{$x_2^*$}}; 
	
	\draw[brightred, line width=3pt] (2,0.335) -- (2,-0.335);
	\draw[brightred, line width=3pt] (2,0.28) -- (2.2,0.28);
	\draw[brightred, line width=3pt] (2,-0.28) -- (2.2,-0.28); 
	\node[above=15pt, font=\large, brightred] at (2, 0) {$\lambda=1$}; 
	\draw[brightred, line width=3pt] (2,0) -- (5.7,0);
	\draw[black, line width=3pt, fill=black] (2,0) circle (0.05);
	
	\draw[brightblue, line width=3pt] (4,0.335) -- (4,-0.335);
	\draw[brightblue, line width=3pt] (4,0.28) -- (4.2,0.28);
	\draw[brightblue, line width=3pt] (4,-0.28) -- (4.2,-0.28); 
	\node[above=15pt, font=\large, brightblue] at (4, 0) {$\lambda=2$}; 
	\draw[brightblue, line width=3pt, dotted] (4,0) -- (5.7,0); 
	\draw[black, line width=3pt, fill=black] (4,0) circle (0.05);
	
	\draw[->, black, thick] (2.7,0) -- (3.3,0);
\end{tikzpicture}
    \caption{Movement of the feasible region of~$\Pi_{\mathrm{Inf}}$. The feasible region is shown in \textcolor{brightred}{red} for~$\lambda=1$ and in \textcolor{brightblue}{blue} for~$\lambda=2$. Although the unique optimal solution moves from~$x=1$ to~$x=2$ as~$\lambda$ varies over~$[1,2]$, it is generated by the same active parametric constraint $x\geq \lambda$ for every~$\lambda\in[1,2]$.}
    \label{fig:intracFeasReg}
    \end{figure}
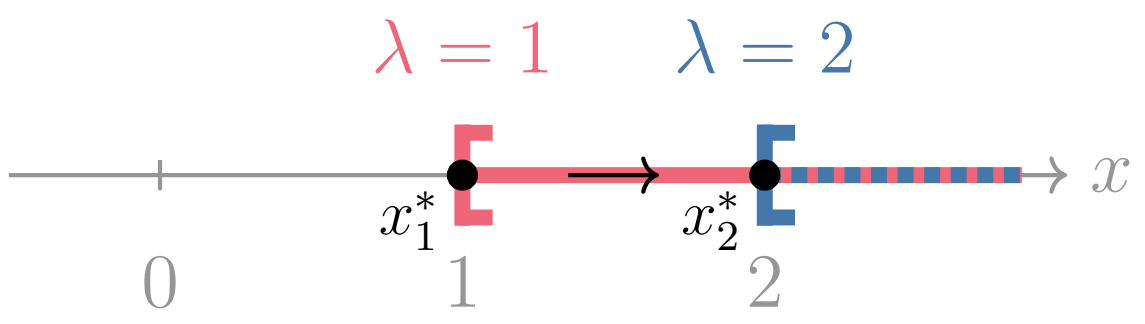

    Although infinitely many distinct optimal solution vectors are required for~$\Pi_{\mathrm{Inf}}$, their underlying basis structure is considerably simpler. As illustrated in Figure~\ref{fig:intracFeasReg}, the unique optimal solution moves with the parameter, but is generated by the same active parametric constraint for every~$\lambda\in[1,2]$.
    
    Motivated by this observation, we next represent solutions of an RHS-PLPP by their underlying bases. For the discussion of bases below, we transform each non-parametric problem~$\Pi_{\mathrm{RHS}}(\lambda)$ into standard equality form by introducing surplus variables and, for notational simplicity, continue to denote the resulting constraint matrix, variable vector, objective vector, and number of variables by~$A$, $x$, $c$, and~$n$, respectively.

    Let $B\subseteq[n]$ with $|B|=m$ be such that~$A_B$ is nonsingular. We call~$B$ a \emph{basis} and the index set $N\coloneq[n]\setminus B$ the corresponding \emph{nonbasis}. For a parameter value~$\lambda\in\Lambda$, the \emph{basic solution} associated with~$B$ is the vector $x^B(\lambda)\in\mathbb{R}^n$ defined by
    $x_N^B(\lambda)\coloneq 0$ and
    $x_B^B(\lambda)\coloneq A_B^{-1}(b+\lambda r)$.
    The set of parameter values for which this basic solution is feasible is denoted by
    \[
    \Lambda_B\coloneq
    \{\lambda\in\Lambda:x^B(\lambda)\in X_\lambda\}.
    \]
    A basis~$B$ is called \emph{optimal} for~$\lambda\in\Lambda$ if $A_B^{-1}(b+\lambda r)\geq0$ and $c_B^\T A_B^{-1}A\leq c^\T$. In this case, its associated basic solution~$x^B(\lambda)$ is an optimal solution of~$\Pi_{\mathrm{RHS}}(\lambda)$.
    A set~$S_B$ of bases is called a \emph{basic optimal solution set} for~$\Pi_{\mathrm{RHS}}$ if, for every~$\lambda\in\Lambda$, it contains a basis that is optimal for~$\lambda$.

    \begin{lemma}[Breakpoints]\label{Lemma:Breakpoints}
    Let $\Pi_{\mathrm{RHS}}$ be an RHS-PLPP instance and $\Pi_{\mathrm{OBF}}$ its corresponding dual OBF-PLPP instance. At every breakpoint of their common optimal cost curve~$F^*=G^*$, the basis used for~$\Pi_{\mathrm{RHS}}$ must change, in the sense that no basis of~$\Pi_{\mathrm{RHS}}$ can be optimal on parameter intervals on both sides of the breakpoint. Consequently, every basic optimal solution set for~$\Pi_{\mathrm{RHS}}$ contains at least as many bases as a minimum-cardinality optimal solution set for~$\Pi_{\mathrm{OBF}}$ contains solutions.
    \end{lemma}
   \begin{proof}
    Fix a basis~$B$ of~$\Pi_{\mathrm{RHS}}$. Its primal feasibility region is
    \[
    \Lambda_B=\bigl\{\lambda\in\Lambda:A_B^{-1}(b+\lambda r)\geq0\bigr\},
    \]
    which is defined by linear inequalities in~$\lambda$ and is therefore an interval. Moreover, the dual feasibility condition
    $
    c_B^\T A_B^{-1}A\leq c^\T
    $
    is independent of~$\lambda$. Consequently, the set of parameter values for which~$B$ is optimal is either~$\Lambda_B$, if~$B$ is dual feasible, or empty, and is therefore an interval.
    The objective value of the basic solution associated with~$B$ is affine in~$\lambda$ on~$\Lambda_B$, since
    \[
    \varphi_B(\lambda)\coloneq F\bigl(x^B(\lambda)\bigr)=c_B^\T A_B^{-1}(b+\lambda r).
    \]
    Suppose that the same basis~$B$ were optimal for parameter values on both sides of a breakpoint. Since its optimality region is an interval, it would then be optimal throughout a neighborhood of the breakpoint. Hence, the optimal cost curve would coincide with the same affine function~$\varphi_B$ on both sides of the breakpoint and could not change its slope there, contradicting the definition of a breakpoint. Thus, the basis used for~$\Pi_{\mathrm{RHS}}$ must change at every breakpoint.
    
    By Lemma~\ref{Lemma:Duality}, $\Pi_{\mathrm{RHS}}$ and $\Pi_{\mathrm{OBF}}$ have the same optimal cost curve. Let~$q$ denote the number of its breakpoints. A minimum-cardinality optimal solution set for~$\Pi_{\mathrm{OBF}}$ then contains exactly~$q+1$ solutions~\cite{bazgan_approximation_2022}. Similarly, a basic optimal solution set for~$\Pi_{\mathrm{RHS}}$ requires one basis before the first breakpoint and a different basis at each of the~$q$ breakpoints. Since the optimality region of every basis is an interval, no basis can be reused after it ceases to be optimal. Thus, every basic optimal solution set for~$\Pi_{\mathrm{RHS}}$ contains at least~$q+1$ bases, and the claimed cardinality bound follows.
    \end{proof}

    The breakpoints of an OBF-PLPP instance partition the parameter set into maximal linear segments of its optimal cost curve, and a single optimal solution suffices on each such segment. The analogous statement does not hold for basic optimal solution sets for the corresponding RHS-PLPP instance: even when there are no breakpoints—in fact, when the common optimal cost curve is constant on the entire parameter set—arbitrarily many distinct bases may be required. The proof is provided in the \hyperref[app:proof-basis-intractability]{Appendix}.

    \begin{lemma}[Multiple Bases Without Breakpoints]
    \label{Lemma:LinearSections}
    For every~$q\in\mathbb{N}$, there exists an RHS-PLPP instance whose optimal cost curve is constant on the entire parameter set and for which every basic optimal solution set contains at least~$q$ distinct bases.
    \end{lemma}

    The preceding lemma shows that the number of bases required by an RHS-PLPP instance is not determined by the number of breakpoints of its optimal cost curve. The following classical result from \citeauthor{murty_computational_1980} \cite{murty_computational_1980} shows that a basic optimal solution set may have super-polynomial cardinality in the encoding length of the instance: 

    \begin{lemma}[Super-polynomial Basis Complexity]
    \label{lemma:bases_intract}
    For every~$n\in\mathbb{N}$, there exists an RHS-PLPP instance of encoding length bounded by~$\operatorname{poly}(n)$ for which every basic optimal solution set contains at least~$2^{n-1}$ distinct bases.
    \end{lemma}

    The preceding results show that neither explicit solution vectors nor bases yield polynomial-size exact representations in general: an optimal solution set consisting of explicit vectors may necessarily be infinite, whereas a basic optimal solution set may require super-polynomially many bases. This motivates replacing optimal solution sets with approximation sets, as considered in the next section.

    \subsection{Parametric Approximation}
    \label{Sec:ParametricApproximation}
    
    We now formalize the two multiplicative approximation concepts used in the remainder of the paper. Throughout, we maintain the standing assumptions introduced at the end of Section~\ref{sec:defRHS}. By Lemma~\ref{Lemma:Duality}, these assumptions imply that $F^*(\lambda)=G^*(\lambda)>0$ for every~$\lambda\in\Lambda$, so multiplicative approximation guarantees are well-defined for both the RHS-PLPP and its corresponding dual OBF-PLPP.
    
    The following two definitions formalize our first approximation concept, which preserves feasibility while allowing a multiplicative loss in objective quality. 
    Our definitions follow the approximation-set framework of \citeauthor{bazgan_approximation_2022}~\cite{bazgan_approximation_2022} for OBF-PLPPs and extend it to RHS-PLPPs. Parametric approximation itself has been studied earlier (see, e.g., \cite{katoh_parametric_1987}).
    
    \begin{definition}[Parametric Approximation]
    \label{Def:Approximation}
    Let $\alpha>1$ and $\lambda\in\Lambda$.
    \begin{enumerate}
        \item A solution $x\in X_\lambda$ is an \emph{$\alpha$-approximate solution} of~$\Pi_{\mathrm{RHS}}(\lambda)$ if $F(x)\leq\alpha F^*(\lambda)$.
    
        \item A solution $y\in Y$ is an \emph{$\alpha$-approximate solution} of~$\Pi_{\mathrm{OBF}}(\lambda)$ if $G_\lambda(y)\geq\frac{1}{\alpha}G^*(\lambda)$.
    \end{enumerate}
    
    For a parameter subset~$\Lambda'\subseteq\Lambda$, a set $S_\alpha^{\mathrm{RHS}}\subseteq\mathbb{R}_{\geq0}^n$ is called an \emph{$\alpha$-approximation set for~$\Pi_{\mathrm{RHS}}$ on~$\Lambda'$} if, for every~$\lambda\in\Lambda'$, there exists a solution $x\in S_\alpha^{\mathrm{RHS}}\cap X_\lambda$ satisfying $F(x)\leq\alpha F^*(\lambda)$.
    
    Analogously, a set $S_\alpha^{\mathrm{OBF}}\subseteq Y$ is called an \emph{$\alpha$-approximation set for~$\Pi_{\mathrm{OBF}}$ on~$\Lambda'$} if, for every~$\lambda\in\Lambda'$, there exists a solution $y\in S_\alpha^{\mathrm{OBF}}$ satisfying $G_\lambda(y)\geq\frac{1}{\alpha}G^*(\lambda)$.

    For brevity, we say that a solution $x\in\mathbb{R}_{\geq0}^n$ \emph{$\alpha$-approximates~$\Pi_{\mathrm{RHS}}$ on~$\Lambda'$} if $\{x\}$ is an $\alpha$-approximation set for~$\Pi_{\mathrm{RHS}}$ on~$\Lambda'$. Analogously, a solution $y\in Y$ \emph{$\alpha$-approximates~$\Pi_{\mathrm{OBF}}$ on~$\Lambda'$} if $\{y\}$ is an $\alpha$-approximation set for~$\Pi_{\mathrm{OBF}}$ on~$\Lambda'$.
    
    When $\Lambda'=\Lambda$, we simply speak of an $\alpha$-approximation set for the corresponding parametric instance. When the instance is clear from the context, we omit the superscripts $\mathrm{RHS}$ and $\mathrm{OBF}$.
    
    An $\alpha$-approximation set is called \emph{minimum-cardinality} if it has minimum cardinality among all $\alpha$-approximation sets for the same instance and parameter subset.
    \end{definition}
    
    \begin{definition}[Parametric Approximation Algorithms]
    \label{Def:ApproximationAlgorithms}
    For a fixed~$\alpha>1$, an algorithm is called an \emph{$\alpha$-approximation algorithm for RHS-PLPPs} if, for every RHS-PLPP instance~$\Pi_{\mathrm{RHS}}$, it computes an $\alpha$-approximation set for~$\Pi_{\mathrm{RHS}}$ in time polynomial in the encoding length of the instance. An \emph{$\alpha$-approximation algorithm for OBF-PLPPs} is defined analogously.

    A \emph{polynomial-time approximation scheme} (PTAS) for RHS-PLPPs or OBF-PLPPs is an algorithm that takes as input an instance $\Pi_{\mathrm{RHS}}$ or $\Pi_{\mathrm{OBF}}$ and an $\varepsilon>0$, and returns a $(1+\varepsilon)$-approximation set in time polynomial in the encoding length of the instance for every fixed value of $\varepsilon$. It is called a \emph{fully polynomial-time approximation scheme} (FPTAS) if the running time is polynomial in both the encoding length of the instance and~$\frac{1}{\varepsilon}$.
    \end{definition}
    
    We refer to the approximation concept formalized in Definition~\ref{Def:Approximation} as \emph{primal approximation}; Definition~\ref{Def:ApproximationAlgorithms} specifies the corresponding primal approximation algorithms. By contrast, our second concept, \emph{dual approximation}, requires the objective-value guarantee $F(x)\leq F^*(\lambda)$ but permits a controlled multiplicative relaxation of the parameter-dependent right-hand sides. Here, the term ``dual'' refers to this exchange between objective-value approximation and feasibility relaxation and is unrelated to the LP duality between~$\Pi_{\mathrm{RHS}}$ and~$\Pi_{\mathrm{OBF}}$.

    To ensure that division by~$\alpha$ weakens every parameter-dependent right-hand side, we restrict dual approximation to RHS-PLPP instances with strictly positive right-hand sides.
    
    \begin{definition}[Dual Parametric Approximation]
    \label{Def:DualApprox}
    Let $\Pi_{\mathrm{RHS}}$ be an RHS-PLPP instance satisfying
    \begin{equation}
    \label{eq:positive-right-hand-sides}
    b_i+\lambda r_i>0
    \qquad\text{for every }i\in[m]\text{ and }\lambda\in\Lambda.
    \end{equation}
    For an approximation factor~$\alpha>1$ and a parameter value~$\lambda\in\Lambda$, define the \emph{$\alpha$-relaxed feasible set}
    \[
    X_\lambda^\alpha
    \coloneq
    \left\{
    x\in\mathbb{R}_{\geq0}^n:
    Ax\geq\frac{1}{\alpha}(b+\lambda r)
    \right\}
    \supseteq X_\lambda.
    \]
    
    A solution $x\in X_\lambda^\alpha$ is called a \emph{dual $\alpha$-approximate solution} of~$\Pi_{\mathrm{RHS}}(\lambda)$ if $F(x)\leq F^*(\lambda)$.
    
    For a parameter subset~$\Lambda'\subseteq\Lambda$, a set $S_\alpha^{\mathrm{D}}\subseteq\mathbb{R}_{\geq0}^n$ is called a \emph{dual $\alpha$-approximation set for~$\Pi_{\mathrm{RHS}}$ on~$\Lambda'$} if, for every~$\lambda\in\Lambda'$, there exists a solution $x\in S_\alpha^{\mathrm{D}}\cap X_\lambda^\alpha$ satisfying $F(x)\leq F^*(\lambda)$.
    
    When $\Lambda'=\Lambda$, we simply call~$S_\alpha^{\mathrm{D}}$ a dual $\alpha$-approximation set for~$\Pi_{\mathrm{RHS}}$. A dual $\alpha$-approximation set is called \emph{minimum-cardinality} if it has minimum cardinality among all dual $\alpha$-approximation sets for the same instance and parameter subset.
    \end{definition}
    
    \begin{definition}[Dual Parametric Approximation Algorithms]
    \label{Def:DualApproximationAlgorithms}
    For a fixed~$\alpha>1$, an algorithm is called a \emph{dual $\alpha$-approximation algorithm for RHS-PLPPs} if, for every RHS-PLPP instance~$\Pi_{\mathrm{RHS}}$ satisfying~\eqref{eq:positive-right-hand-sides}, it computes a dual $\alpha$-approximation set for~$\Pi_{\mathrm{RHS}}$ in time polynomial in the encoding length of the instance.

    A \emph{dual polynomial-time approximation scheme} (dual PTAS) for RHS-PLPPs is an algorithm that takes as input an instance $\Pi_{\mathrm{RHS}}$ and an $\varepsilon>0$, and returns a dual $(1+\varepsilon)$-approximation set in time polynomial in the encoding length of the instance for every fixed value of $\varepsilon$. It is called a \emph{dual fully polynomial-time approximation scheme} (dual FPTAS) if the running time is polynomial in both the encoding length of the instance and~$\frac{1}{\varepsilon}$.
    \end{definition}

	\section{Primal Approximation Algorithm}\label{sec:primal_approx_algo}
    In this section, let an RHS-PLPP instance $\Pi_{\mathrm{RHS}}=(\Lambda,X,F)$ and an approximation factor $1+\varepsilon$ with $\varepsilon\in\mathbb{Q}_{>0}$ be given. We develop a (primal) $(1+\varepsilon)$-approximation algorithm for RHS-PLPPs satisfying certain sufficient conditions. The basic idea is to approximate the convex optimal cost curve $F^*$ by a staircase whose levels lie between~$F^*$ and $(1+\varepsilon)\cdot F^*$. Starting from the global minimizer~$\lambda_{\mathrm{PS}}$ returned by an RHS-parametric search oracle, the algorithm iteratively extends the covered parameter region in both directions as far as possible at each successive staircase level.

    Algorithm~\ref{algo::aprox1D} constructs approximating solutions together with parameter intervals on which they provide the required approximation guarantee. For convenience, we record this information explicitly by associating each selected solution with such an interval. Thus, Algorithm~\ref{algo::aprox1D} returns a set of triples $(x,\lambda^{\mathrm{LB}},\lambda^{\mathrm{UB}})$, while the corresponding approximation set in the sense of Definition~\ref{Def:Approximation} consists only of the selected solutions. 
    The additional information concerning the intervals is not essential to the resulting approximation set, but is typically part of the solution representation in exact parametric optimization~\cite{gal_multiparametric_1972,borrelli_geometric_2003}. Recording the intervals explicitly motivates the following notion of a covering set.

    \begin{definition}[$\bm{(1+\varepsilon)$}-Covering Set]\label{Def:CovSet}
        Let $\varepsilon>0$ and $\lambda_1,\lambda_2\in\Lambda$ with $\lambda_1\leq\lambda_2$.
        A finite set $S\subseteq X\times\Lambda\times\Lambda$ is called a \emph{$(1+\varepsilon)$-covering set for $[\lambda_1,\lambda_2]$} if, for every $(x,\lambda^{\mathrm{LB}},\lambda^{\mathrm{UB}})\in S$, we have $\lambda_1\leq\lambda^{\mathrm{LB}}\leq\lambda^{\mathrm{UB}}\leq\lambda_2$, and the solution $x$ $(1+\varepsilon)$-approximates $\Pi_{\mathrm{RHS}}$ on $[\lambda^{\mathrm{LB}},\lambda^{\mathrm{UB}}]$, and if $
        [\lambda_1,\lambda_2] = \bigcup_{(x,\lambda^{\mathrm{LB}},\lambda^{\mathrm{UB}})\in S}
        [\lambda^{\mathrm{LB}},\lambda^{\mathrm{UB}}]$.
	\end{definition}

    For a set $S\subseteq X\times\Lambda\times\Lambda$, let 
    \[\operatorname{sol}(S)\coloneqq\{x\in X:\exists\,\lambda^{\mathrm{LB}},\lambda^{\mathrm{UB}}\in\Lambda\text{ with }(x,\lambda^{\mathrm{LB}},\lambda^{\mathrm{UB}})\in S\}\]
    denote its \emph{solution projection}. If~$S$ is a $(1+\varepsilon)$-covering set for $[\lambda_1,\lambda_2]\subseteq\Lambda$,
    then $\operatorname{sol}(S)$ is a $(1+\varepsilon)$-approximation set for $\Pi_{\mathrm{RHS}}$ on $[\lambda_1,\lambda_2]$ in the sense of Definition~\ref{Def:Approximation}. Clearly, $|\operatorname{sol}(S)|\leq|S|$, with strict inequality possible if the same solution occurs in several triples associated with different parameter intervals.

    \medskip

    The following observation, which follows from Corollary~\ref{Cor:OptCostCurve}, is central to the construction. Let $[\lambda_1,\lambda_2]\subseteq\Lambda$ with $\lambda_{\mathrm{PS}}\notin(\lambda_1,\lambda_2)$. Since $F^*$ is convex and $\lambda_{\mathrm{PS}}$ is a global minimizer, $F^*$ is non-increasing to the left of~$\lambda_{\mathrm{PS}}$ and non-decreasing to its right. Hence, $\min_{\lambda\in[\lambda_1,\lambda_2]}F^*(\lambda)=\min\{F^*(\lambda_1),F^*(\lambda_2)\}$.

    To describe how far the staircase construction can extend the covered parameter region at the objective threshold $(1+\varepsilon)\cdot F^*(\lambda)$, we introduce the following sublevel-set boundary functions:
    \begin{equation}\label{eq:epsilonParameterBoundaries}
    \begin{aligned}
    L_\varepsilon(\lambda) &\coloneqq \min\left\{\bar{\lambda}\in\Lambda:F^*(\bar{\lambda})\leq(1+\varepsilon)\cdot F^*(\lambda)\right\} &&\text{for }\lambda\in[\lambda_{\min},\lambda_{\mathrm{PS}}],\\
    R_\varepsilon(\lambda) &\coloneqq \max\left\{\bar{\lambda}\in\Lambda:F^*(\bar{\lambda})\leq(1+\varepsilon)\cdot F^*(\lambda)\right\} &&\text{for }\lambda\in[\lambda_{\mathrm{PS}},\lambda_{\max}].
    \end{aligned}
    \end{equation}
    Thus, $[L_\varepsilon(\lambda),\lambda]$ is the largest interval ending at~$\lambda$ such that $F^*(\bar{\lambda})\leq(1+\varepsilon)\cdot F^*(\lambda)$ for every $\bar{\lambda}$ in this interval, while $[\lambda,R_\varepsilon(\lambda)]$ is the corresponding largest interval starting at~$\lambda$. The boundary functions~$L_\varepsilon$ and~$R_\varepsilon$ are well defined. Indeed, by continuity and convexity of~$F^*$, every nonempty sublevel set $\{\bar{\lambda}\in\Lambda:F^*(\bar{\lambda})\leq H\}$ is a closed interval. For $H=(1+\varepsilon)\cdot F^*(\lambda)$, this interval is nonempty since it contains~$\lambda$. Moreover, its endpoints can be computed by standard techniques from linear programming. If an endpoint lies in the interior of~$\Lambda$, continuity of~$F^*$ implies that $F^*$ attains the value~$H$ at this endpoint.
        
    By the monotonicity stated above, $F^*(\lambda)$ is the minimum of~$F^*$ on $[L_\varepsilon(\lambda),\lambda]$ for $\lambda\in[\lambda_{\min},\lambda_{\mathrm{PS}}]$ and on $[\lambda,R_\varepsilon(\lambda)]$ for $\lambda\in[\lambda_{\mathrm{PS}},\lambda_{\max}]$. This motivates the following definition of a covering oracle for these two types of intervals.

    \begin{definition}[$\bm{(1+\varepsilon)}$-Covering Oracle]\label{Def:ApproxIntervalOracle}
    A \emph{$(1+\varepsilon)$-covering oracle} for an RHS-PLPP instance $\Pi_{\mathrm{RHS}}=(\Lambda,X,F)$ takes as input three rational values $\lambda_1,\lambda_2,H$ of one of the following two forms: 
    \begin{itemize}
        \item either $\bigl(\lambda_1,\lambda_2,H\bigr)=\bigl(L_\varepsilon(\lambda),\lambda,(1+\varepsilon)\cdot F^*(\lambda)\bigr)$ for some $\lambda\in[\lambda_{\min},\lambda_{\mathrm{PS}}]$, 
        \item or $\bigl(\lambda_1,\lambda_2,H\bigr)=\bigl(\lambda,R_\varepsilon(\lambda),(1+\varepsilon)\cdot F^*(\lambda)\bigr)$ for some $\lambda\in[\lambda_{\mathrm{PS}},\lambda_{\max}]$. 
    \end{itemize}
    
    It returns a $(1+\varepsilon)$-covering set $S$ for $[\lambda_1,\lambda_2]$ such that $F(x)\leq H$ for every solution $x\in\operatorname{sol}(S)$, where $H$ is the common objective threshold.
    We denote the returned covering set by $(1+\varepsilon)\text{-}\operatorname{Cov}(\lambda_1,\lambda_2,H)$. The oracle is called \emph{polynomial-time} if its running time is polynomial in the encoding lengths of $\Pi_{\mathrm{RHS}}$, $\varepsilon$, $\lambda_1$, $\lambda_2$, and~$H$, and in~$1/\varepsilon$.
    \end{definition}
    
    The threshold $H$ need not be given explicitly, since it can be computed in polynomial time from $\lambda_1$ and $\lambda_2$ by solving the corresponding pointwise RHS-PLPP. For convenience, however, we include it in the oracle input.
    The construction of the covering set in Algorithm~\ref{algo::aprox1D} relies on Assumption~\ref{Assumption:Standing} and the availability of a polynomial-time covering oracle. We summarize these requirements in the following assumption.

    \begin{assumption}[Sufficient Approximation Conditions]\label{Def:SufficientCond}
    For a given approximation factor $1+\varepsilon$ with $\varepsilon\in\mathbb{Q}_{>0}$, an RHS-PLPP instance $\Pi_{\mathrm{RHS}}$ satisfies the following \emph{sufficient approximation conditions}:
    \begin{enumerate}
        \item Assumption~\ref{Assumption:Standing} holds.
        \item A polynomial-time $(1+\varepsilon)$-covering oracle for $\Pi_{\mathrm{RHS}}$ in the sense of Definition~\ref{Def:ApproxIntervalOracle} is available.
    \end{enumerate}
    \end{assumption}
	
\begin{figure}[h]
		\centering
		\begin{tikzpicture}
	
	\definecolor{mygreen}{RGB}{200,255,200}
	\definecolor{mygray}{RGB}{230,230,230}
	\definecolor{myblue}{RGB}{200,220,255}
	
	\coordinate (P0) at (0,4);
	\coordinate (P1) at (1,2.5);
	\coordinate (P2) at (2,1.5);
	\coordinate (P3) at (3,1);
	\coordinate (P4) at (4,1);
	\coordinate (P5) at (5,1.5);
	\coordinate (P6) at (6,2.5);
	\coordinate (P7) at (7,4);
	
	\coordinate (D0) at (0,6);
	\coordinate (D1) at (1,3.75);
	\coordinate (D2) at (2,2.25);
	\coordinate (D3) at (3,1.5);
	\coordinate (D4) at (4,1.5);
	\coordinate (D5) at (5,2.25);
	\coordinate (D6) at (6,3.75);
	\coordinate (D7) at (7,6);
	
	\begin{scope}[on background layer]
        \fill[brightgrey!40,opacity=0.4] (P0)--(P1)--(P2)--(3,1) -- (3.5,1)-- (3.5,6.5) -- (0,6.5);
		\fill[brightgrey!75,opacity=0.4]  (3.5,1)-- (4,1) --(P5)--(P6)--(P7) -- (7,6.5) -- (3.5,6.5);
	\end{scope}
	
	\draw[->] (-0.5,0.5) -- (7.5,0.5) node[right] {$\lambda$};
	
	\draw (0,0.5) node[below] {$\lambda_{\min}$};
	\draw [black,thick] (0,0.57)--(0,0.43);
	\draw (3.5,0.5) node[below,brightgreen] {$\lambda_{\text{PS}}$};
	\draw [brightgreen,thick] (3.5,0.57)--(3.5,0.43);
	\draw (7,0.5) node[below] {$\lambda_{\max}$};
	\draw [black,thick] (7,0.57)--(7,0.43);
	
	\draw[thick, black] (P0) -- (P1) -- (P2) -- (P3) -- (P4) -- (P5) -- (P6) -- (P7);
	\foreach \p in {P0,P1,P2,P3,P4,P5,P6,P7} {
		\fill[black] (\p) circle (2pt);
	}
	
	\draw[thick, brightpurple]
	(D0) -- (D1) -- (D2) -- (D3) -- (D4) -- (D5) -- (D6) -- (D7);
	\foreach \p in {D0,D1,D2,D3,D4,D5,D6,D7} {
		\fill[brightpurple] (\p) circle (2pt);
	}
	
	\draw[gray!50, thin] (0,4.5) -- (7,4.5);
	\draw[gray!50, thin] (0,4.7) -- (7,4.7);
	\draw[gray!50, thin] (0,4.9) -- (7,4.9);
	\draw[gray!50, thin] (1,2.5) -- (6,2.5);
	\draw[gray!50, thin] (0.85,2.7) -- (6.15,2.7);
	\draw[gray!50, thin] (0.72,2.9) -- (6.28,2.9);
	
	\coordinate (S0) at (3.5,1);
	\coordinate (S1) at (3.5,1.5);
	\coordinate (S2) at (2,1.5);
	\coordinate (S3) at (5,1.5);
	\coordinate (S4) at (5,2.25);
	\coordinate (S5) at (2,2.25);
	\coordinate (S6) at (1.27,2.25);
	\coordinate (S7) at (5.73,2.25);
	\coordinate (S8) at (1.27,3.37);
	\coordinate (S9) at (5.73,3.37);
	\coordinate (S10) at (0.42,3.37);
	\coordinate (S11) at (6.58,3.37);
	\coordinate (S12) at (0.42,5.07);
	\coordinate (S13) at (6.58,5.07);
	\coordinate (S14) at (0,5.07);
	\coordinate (S15) at (7,5.07);
	
	\draw[thick, brightgreen]
	(S0) -- (S1) -- (S2) -- (S5) -- (S6) -- (S8) -- (S10) -- (S12) -- (S14);
	\draw[thick, brightgreen]
	(S1) -- (S3) -- (S4) -- (S7) -- (S9) -- (S11) -- (S13) -- (S15);
	
	\foreach \p in {S0,S1,S2,S3,S4,S5,S6,S7,S8,S9,S10,S11,S12,S13,S14,S15} {
		\fill[brightgreen] (\p) circle (2pt);
	}
	\node at (3.5, 3.65) {\dots};
	\node at (7.4,4) {$F^*$};
	\node[brightpurple] at (8,6) {$(1+\varepsilon) F^*$};
\end{tikzpicture}
		\caption{Schematic illustration of the greedy staircase construction in Algorithm~\ref{algo::aprox1D}.}
        \label{fig:greedy-staircase}
	\end{figure}
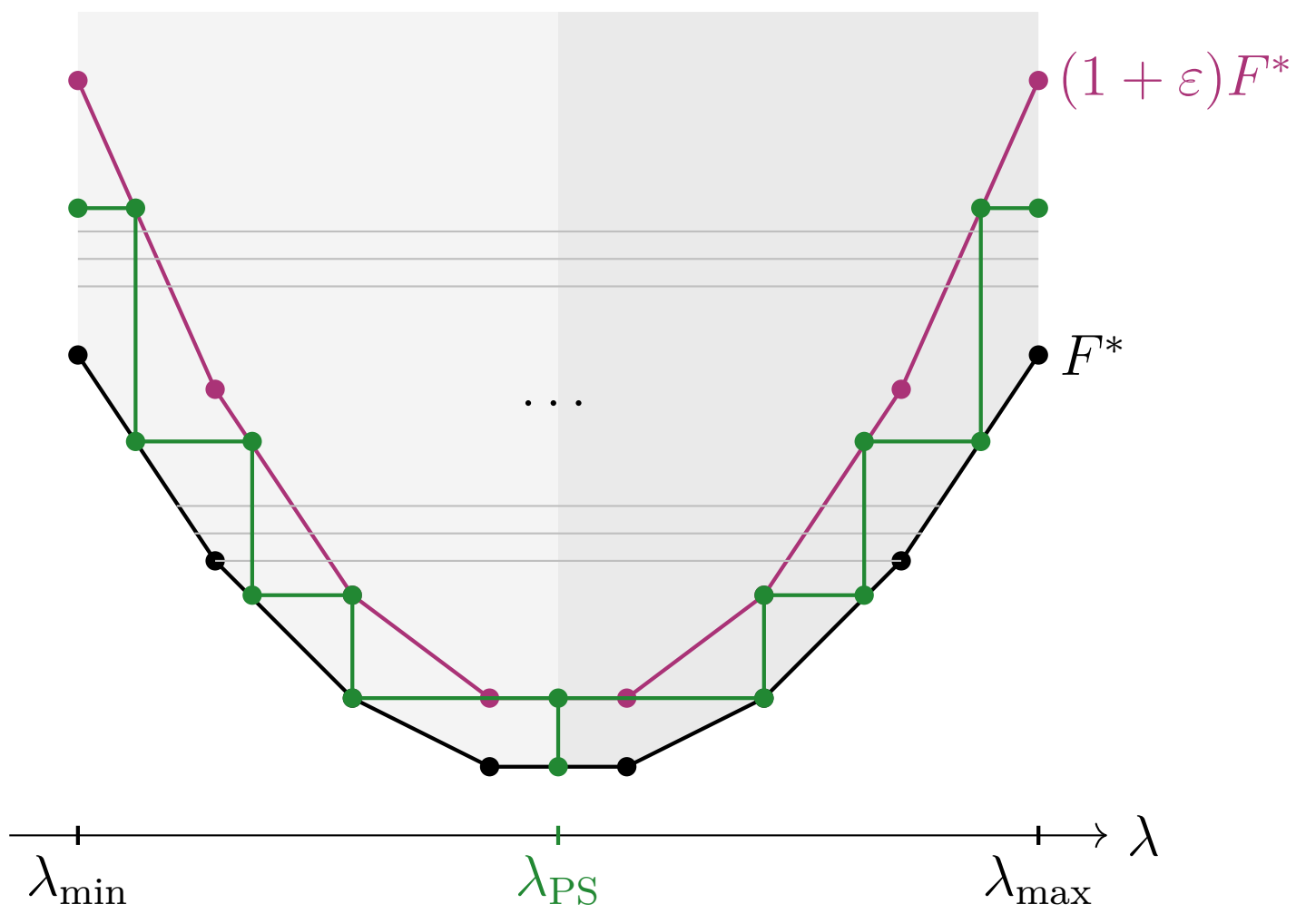

	Under these sufficient approximation conditions, Algorithm~\ref{algo::aprox1D} implements the staircase construction described above, which is illustrated schematically in Figure~\ref{fig:greedy-staircase}. We can now establish the correctness of the algorithm.
	
	\begin{algorithm}[H]
        \caption{A $(1+\varepsilon)$-Covering Algorithm for RHS-PLPPs}\label{algo::aprox1D}
		\begin{algorithmic}[1]
			\Require{}An instance $\Pi_{\text{RHS}}=(\Lambda,X,F)$ of an RHS-PLPP, an approximation factor $1+\varepsilon$ with $\varepsilon\in\mathbb{Q}_{>0}$.
			\Ensure{}A $(1+\varepsilon)$-covering set $S$ for $\Lambda$.
            \State{} Find $(x_\text{PS},\lambda_\text{PS})\gets \mathrm{ALG}_{\mathrm{PS}}(\Pi_{\text{RHS}})$.\label{algoline:aprox1D_initial_oracle_call}
            \If{$\lambda_{\min}=\lambda_{\max}$} \State{} \Return $S \gets \{(x_{\mathrm{PS}},\lambda_{\min},\lambda_{\max})\}$.
            \EndIf
            \State{} Initialize covering set $S \gets \emptyset$.
			\State{} Set $\lambda_{\text{old}}\gets \lambda_{\text{PS}}$, $H\gets (1+\varepsilon)\cdot F(x_{\text{PS}})$.
			\While{$\lambda_{\text{old}} > \lambda_{\min}$}\label{algoline:aprox1D_left-while}\Comment{Left branch}
			\State{} Find $\lambda_{\text{new}}\gets \min \{\lambda'\in\Lambda:\exists\, x\in X_{\lambda'}\text{ with }c^\T x\leq H\}$.\label{algoline:aprox1D_left-lambda_new}
            \State{} Find $S_{\mathrm{new}}\gets (1+\varepsilon)\text{-}\operatorname{Cov} (\lambda_{\mathrm{new}},\lambda_{\mathrm{old}},H)$.
			\State{} Update covering set $S \gets S\cup S_{\text{new}}$.
			\State{} Set $\lambda_{\text{old}}\gets \lambda_{\text{new}}$, $H\gets (1+\varepsilon)\cdot H$.
			\EndWhile
			\State{} Set $\lambda_{\text{old}}\gets \lambda_{\text{PS}}$, $H\gets (1+\varepsilon)\cdot F(x_{\text{PS}})$.\label{algoline:aprox1D_reinit}
			\While{$\lambda_{\text{old}} < \lambda_{\max}$}\label{algoline:aprox1D_right-while}\Comment{Right branch}
			\State{} Find $\lambda_{\text{new}}\gets \max \{\lambda'\in\Lambda:\exists\, x\in X_{\lambda'}\text{ with }c^\T x\leq H\}$.\label{algoline:aprox1D_right-lambda_new}
            \State{} Find $S_{\mathrm{new}}\gets (1+\varepsilon)\text{-}\operatorname{Cov}(\lambda_{\mathrm{old}},\lambda_{\mathrm{new}},H)$.
			\State{} Update covering set $S \gets S\cup S_{\text{new}}$.
			\State{} Set $\lambda_{\text{old}}\gets \lambda_{\text{new}}$, $H\gets (1+\varepsilon)\cdot H$.
			\EndWhile
			\State{} \Return $S$.
		\end{algorithmic}
	\end{algorithm}

\begin{theorem}[Correctness and Finiteness of Algorithm~\ref{algo::aprox1D}]
\label{thm:AlgorithmCorrectness}
Let $\varepsilon\in\mathbb{Q}_{>0}$ and suppose that
$\Pi_{\text{RHS}}=(\Lambda,X,F)$ satisfies the sufficient approximation conditions from Assumption~\ref{Def:SufficientCond}. Then each while loop in Algorithm~\ref{algo::aprox1D} executes only finitely many iterations, and the returned set~$S$ is a $(1+\varepsilon)$-covering set for~$\Lambda$. Consequently, $\operatorname{sol}(S)$ is a $(1+\varepsilon)$-approximation set for~$\Pi_{\mathrm{RHS}}$.
\end{theorem}
\begin{proof}
If $\lambda_{\min}=\lambda_{\max}$, the solution $x_{\mathrm{PS}}$ is
optimal for the unique parameter value, and the claim follows
immediately. Hence, suppose that $\lambda_{\min}<\lambda_{\max}$.

We first show that, upon termination of the two while loops, the set~$S$ returned by Algorithm~\ref{algo::aprox1D} is a $(1+\varepsilon)$-covering set for~$\Lambda$. Afterwards, we show that both while loops indeed terminate after finitely many iterations.

We first consider the left while loop starting in Line~\ref{algoline:aprox1D_left-while}. If $\lambda_{\mathrm{PS}}=\lambda_{\min}$, this loop is not executed, and the interval $[\lambda_{\min},\lambda_{\mathrm{PS}}]$ consists only of the single point $\lambda_{\mathrm{PS}}$, which will be covered by the right while loop considered below. Hence, suppose that $\lambda_{\mathrm{PS}}>\lambda_{\min}$. We show that, upon termination, the set of triples added to~$S$ during the left while loop is a $(1+\varepsilon)$-covering set for $[\lambda_{\min},\lambda_{\mathrm{PS}}]$.
To establish this claim, we show by induction that, in every executed
iteration, we have
\[
H=(1+\varepsilon)\cdot F^*(\lambda_{\mathrm{old}})
\]
at the beginning of the iteration and that the parameter value~$\lambda_{\mathrm{new}}$ computed in this iteration satisfies $\lambda_{\mathrm{new}} = L_\varepsilon(\lambda_{\mathrm{old}})$. These two relations imply that the covering-oracle call made in such an iteration is precisely
\[
(1+\varepsilon)\text{-}\operatorname{Cov}
\bigl(
L_\varepsilon(\lambda_{\mathrm{old}}),
\lambda_{\mathrm{old}},
(1+\varepsilon)\cdot F^*(\lambda_{\mathrm{old}})
\bigr),
\]
which is available by Assumption~\ref{Def:SufficientCond} and returns a $(1+\varepsilon)$-covering set for $[\lambda_{\mathrm{new}},\lambda_{\mathrm{old}}]$.

Initially, $\lambda_{\mathrm{old}}=\lambda_{\mathrm{PS}}$ and $H=(1+\varepsilon)\cdot F(x_{\mathrm{PS}})
=(1+\varepsilon)\cdot F^*(\lambda_{\mathrm{PS}})$,
so the first relation holds for the first iteration. Hence, the parameter value~$\lambda_{\mathrm{new}}$ computed in the first iteration satisfies
\[
\begin{aligned}
\lambda_{\mathrm{new}}
&=
\min\left\{
\lambda'\in\Lambda:
\exists\,x\in X_{\lambda'}\text{ with }c^\T x\leq H
\right\}\\
&=
\min\left\{
\lambda'\in\Lambda:
F^*(\lambda')\leq(1+\varepsilon)
\cdot F^*(\lambda_{\mathrm{old}})
\right\}\\
&=
L_\varepsilon(\lambda_{\mathrm{old}}),
\end{aligned}
\]
as claimed. Now consider an arbitrary executed iteration for which the two claimed relations hold. If $\lambda_{\mathrm{new}}=\lambda_{\min}$, this is the final iteration of the left while loop. 
If $\lambda_{\mathrm{new}}>\lambda_{\min}$, then $\lambda_{\mathrm{new}}$ is the left endpoint of the sublevel set $\{\lambda'\in\Lambda:F^*(\lambda')\leq H\}$. Moreover, since $H=(1+\varepsilon)\cdot F^*(\lambda_{\mathrm{old}}) >F^*(\lambda_{\mathrm{old}})$ and $\lambda_{\mathrm{old}}>\lambda_{\min}$, continuity of~$F^*$ implies that the sublevel set contains parameter values strictly smaller than $\lambda_{\mathrm{old}}$. Hence, $\lambda_{\mathrm{new}}<\lambda_{\mathrm{old}}\leq\lambda_{\max}$. Together with $\lambda_{\mathrm{new}}>\lambda_{\min}$, this shows that $\lambda_{\mathrm{new}}$ lies in the interior of~$\Lambda$. By continuity of~$F^*$, we therefore have $F^*(\lambda_{\mathrm{new}})=H$.
After the updates $\lambda_{\mathrm{old}}\leftarrow\lambda_{\mathrm{new}}$ and $H\leftarrow(1+\varepsilon)\cdot H$, it follows that
\[
H=(1+\varepsilon)\cdot F^*(\lambda_{\mathrm{old}}).
\]
Thus, the first claimed relation holds at the beginning of the next iteration. Applying the same calculation as above in that iteration then yields $\lambda_{\mathrm{new}}=L_\varepsilon(\lambda_{\mathrm{old}})$, which establishes the second claimed relation and completes the induction.

Thus, in every executed iteration of the left while loop, the algorithm invokes the covering oracle for the interval $[\lambda_{\mathrm{new}},\lambda_{\mathrm{old}}]$, where $\lambda_{\mathrm{new}}=L_\varepsilon(\lambda_{\mathrm{old}})$, and the oracle returns a $(1+\varepsilon)$-covering set for this interval. The first such interval has right endpoint $\lambda_{\mathrm{PS}}$. Since $\lambda_{\mathrm{old}}$ is updated to the previously computed value $\lambda_{\mathrm{new}}$ after every non-final iteration, the intervals corresponding to successive iterations are contiguous. Upon termination, the final such interval has left endpoint $\lambda_{\min}$. Hence, the parameter intervals associated with the triples added to~$S$ during the left while loop jointly cover $[\lambda_{\min},\lambda_{\mathrm{PS}}]$. Moreover, every such triple satisfies the required $(1+\varepsilon)$-approximation guarantee on its associated interval. Consequently, the set of triples added to~$S$ during the left loop is a $(1+\varepsilon)$-covering set for $[\lambda_{\min},\lambda_{\mathrm{PS}}]$, as claimed.

\medskip

The corresponding statement for the right while loop starting in Line~\ref{algoline:aprox1D_right-while} follows by a symmetric argument.

Combining the two parts, upon termination of the while loops, the triples
added to~$S$ jointly form a $(1+\varepsilon)$-covering set for $[\lambda_{\min},\lambda_{\mathrm{PS}}] \cup [\lambda_{\mathrm{PS}},\lambda_{\max}] = \Lambda$.
Indeed, if one of the two intervals is a singleton, its unique point $\lambda_{\mathrm{PS}}$ is covered by the other while loop, as discussed above.

\medskip

It remains to show that both while loops terminate after finitely many iterations. Consider first the left while loop.
Since $H$ is initialized at $(1+\varepsilon)\cdot F^*(\lambda_{\mathrm{PS}})$ and $F^*(\lambda_{\mathrm{PS}})>0$ by Assumption~\ref{Assumption:Standing}, we have $H>0$ initially. Moreover, $F^*(\lambda_{\min})$ is finite and $H$ is multiplied by $1+\varepsilon>1$ after every non-final iteration. Hence, its value grows geometrically and eventually satisfies $H\geq F^*(\lambda_{\min})$. Then $\lambda_{\min}$ belongs to the sublevel set \[\{\lambda'\in\Lambda:F^*(\lambda')\leq H\},\] so the minimization defining $\lambda_{\mathrm{new}}$ yields $\lambda_{\mathrm{new}}=\lambda_{\min}$, and the left while loop terminates. 
The same argument applies to the right while loop: eventually $H\geq F^*(\lambda_{\max})$, so the maximization defining $\lambda_{\mathrm{new}}$ yields $\lambda_{\mathrm{new}}=\lambda_{\max}$. Thus, both while loops terminate after finitely many iterations, which finishes the proof. 
\end{proof}
	
Having established the correctness of Algorithm~\ref{algo::aprox1D}, we now investigate its running time.

\begin{theorem}[Running Time of Algorithm~\ref{algo::aprox1D}]
\label{thm:AlgorithmRuntime}
Let $\varepsilon\in\mathbb{Q}_{>0}$ and suppose that $\Pi_{\text{RHS}}=(\Lambda,X,F)$ satisfies the sufficient approximation conditions from Assumption~\ref{Def:SufficientCond}. Then
Algorithm~\ref{algo::aprox1D} runs in time polynomial in $\langle\Pi_{\mathrm{RHS}}\rangle$, $\langle\varepsilon\rangle$, and
$1/\varepsilon$.
\end{theorem}

\begin{proof}
As stated in the paragraph following Definition~\ref{Definition:OptimizationOracles}, the parametric search problem is a linear programming problem of encoding length polynomially bounded in $\langle\Pi_{\mathrm{RHS}}\rangle$. Consequently, the initial parametric-search call in Line~\ref{algoline:aprox1D_initial_oracle_call} can be evaluated in polynomial time, and the returned values $x_{\mathrm{PS}}$ and $\lambda_{\mathrm{PS}}$ have polynomial encoding length.
If $\lambda_{\min}=\lambda_{\max}$, the algorithm returns immediately
after this call, and the claim follows. Hence, suppose that
$\lambda_{\min}<\lambda_{\max}$.

Algorithm~\ref{algo::aprox1D} starts at the global minimizer $\lambda_{\mathrm{PS}}$ of~$F^*$, with $F(x_{\mathrm{PS}})=F^*(\lambda_{\mathrm{PS}})$. By Corollary~\ref{Cor:OptCostCurve}, $F^*$ is convex and therefore non-increasing on $[\lambda_{\min},\lambda_{\mathrm{PS}}]$ and non-decreasing on $[\lambda_{\mathrm{PS}},\lambda_{\max}]$.
Hence, the maximum of~$F^*$ over $\Lambda$ is attained at one of the interval endpoints~$\lambda_{\min}$ or~$\lambda_{\max}$.
Let
$f_{\min}\coloneqq F^*(\lambda_{\mathrm{PS}})$ and
$f_{\max}\coloneqq
\max\{F^*(\lambda_{\min}),F^*(\lambda_{\max})\}$.
By Assumption~\ref{Assumption:Standing} and Proposition~\ref{Proposition:ExistenceOptimalSolutions}, both values are positive and
finite. Moreover, they are optimal values of linear programming problems with rational input data whose encoding lengths are polynomially bounded in
$\langle\Pi_{\mathrm{RHS}}\rangle$. Standard bounds for linear
programming problems therefore imply that there exists an integer $M\in\mathbb{N}$ whose value is polynomially bounded in $\langle\Pi_{\mathrm{RHS}}\rangle$ such that
\[
2^{-M}\leq f_{\min}\leq f_{\max}\leq2^M.
\]

For each $j\geq1$, let~$H_j$ denote the value of~$H$ at the beginning of the $j$-th iteration of either while loop. By the initialization and the update of~$H$ in Algorithm~\ref{algo::aprox1D}, we have $H_j=(1+\varepsilon)^j\cdot f_{\min}$. Once $H_j\geq f_{\max}$, the construction of~$f_{\max}$ implies that $F^*(\lambda)\leq H_j$ for all $\lambda\in\Lambda$. Hence, the minimization defining $\lambda_{\mathrm{new}}$ in Line~\ref{algoline:aprox1D_left-lambda_new} during the left while loop yields $\lambda_{\min}$, whereas the maximization defining $\lambda_{\mathrm{new}}$ in Line~\ref{algoline:aprox1D_right-lambda_new} during the right while loop yields $\lambda_{\max}$, and the loops terminate. Therefore, each while loop terminates after at most $N\coloneqq\min\{k\in\mathbb{N}_0:
2^M\leq(1+\varepsilon)^k\cdot2^{-M}\}$ iterations.
By the minimality of~$N$, we have $(1+\varepsilon)^{N-1}\cdot 2^{-M}<2^M$. 
Since $1/\ln(1+\varepsilon)\leq1+1/\varepsilon$ for every $\varepsilon>0$, taking logarithms yields
\[
N
<
1+\frac{2M\cdot\ln(2)}{\ln(1+\varepsilon)}
\leq
1+2M\cdot\ln(2)\cdot\left(1+\frac1\varepsilon\right).
\]
Thus, the two while loops together execute at most $2N$ iterations, which is polynomially bounded in $\langle\Pi_{\mathrm{RHS}}\rangle$ and $1/\varepsilon$.

It remains to bound the work performed in each iteration. For the $j$-th iteration, define
\[
P_j\coloneqq
\left\{
(\lambda,x)\in\Lambda\times\mathbb{R}_{\geq0}^n:
Ax\geq b+\lambda r,\;
c^\T x\leq H_j
\right\}.
\]
The minimization and maximization used to compute
$\lambda_{\mathrm{new}}$ in the left and right while loops,
respectively, are the linear programming problems
\[
\min_{(\lambda,x)\in P_j}\lambda
\qquad\text{and}\qquad
\max_{(\lambda,x)\in P_j}\lambda.
\]
Since $\langle H_j\rangle = \langle(1+\varepsilon)^j\cdot f_{\min}\rangle \in \mathcal{O}\bigl(\langle f_{\min}\rangle +j\cdot\langle1+\varepsilon\rangle \bigr)$ and $j\leq N$, the encoding length of~$H_j$, and hence that of either linear programming problem, is polynomial in $\langle\Pi_{\mathrm{RHS}}\rangle$, $\langle\varepsilon\rangle$, and $1/\varepsilon$. Standard bit-complexity bounds for linear programming therefore imply that $\lambda_{\mathrm{new}}$ has polynomial encoding length and can be computed in polynomial time.

Each iteration of one of the while loops additionally performs one call to the covering oracle. By the invariant established in the proof of Theorem~\ref{thm:AlgorithmCorrectness}, every such call is of one of the two forms specified in Definition~\ref{Def:ApproxIntervalOracle}. The encoding lengths of $\lambda_{\mathrm{new}}$ and $H_j$ were shown above to be polynomially bounded in the parameters stated in the theorem. Moreover, $\lambda_{\mathrm{old}}$ is either the polynomially encoded value $\lambda_{\mathrm{PS}}$ returned by the initial parametric-search call in Line~\ref{algoline:aprox1D_initial_oracle_call} or a previously computed value $\lambda_{\mathrm{new}}$. Hence, all inputs of every covering-oracle call have polynomial encoding length. Assumption~\ref{Def:SufficientCond} therefore implies that every such call can be evaluated within the polynomial bound claimed in the theorem. The same holds for all remaining arithmetic and set-update operations, since all manipulated numbers and covering-oracle outputs have polynomial encoding length. Together with the polynomial bound on the number of iterations, this shows that the total running time of Algorithm~\ref{algo::aprox1D} is polynomial in $\langle\Pi_{\mathrm{RHS}}\rangle$, $\langle\varepsilon\rangle$, and $1/\varepsilon$, as claimed.
\end{proof}

Combining the two theorems and suitably normalizing the approximation factor, we obtain an FPTAS:

\begin{corollary}\label{cor:GeneralFPTAS}
Let $\mathcal{P}$ be a class of RHS-PLPP instances such that, for every
$\Pi_{\mathrm{RHS}}\in\mathcal{P}$ and every $\eta\in\mathbb{Q}_{>0}$,
the instance satisfies Assumption~\ref{Def:SufficientCond} with $\eta$
in place of $\varepsilon$.
Given an instance $\Pi_{\mathrm{RHS}}\in\mathcal{P}$ and a desired approximation factor $1+\varepsilon$ with $\varepsilon>0$, set
$\hat{\varepsilon}\coloneqq 1/\lceil 1/\varepsilon\rceil$ and run
Algorithm~\ref{algo::aprox1D} on $\Pi_{\mathrm{RHS}}$ using
$1+\hat{\varepsilon}$ as the input approximation factor.
Then the resulting procedure, with output given by the solution projection of the computed covering set, is an FPTAS for~$\mathcal{P}$.
\end{corollary}
\begin{proof}
Fix $\Pi_{\mathrm{RHS}}\in\mathcal{P}$ and $\varepsilon>0$, and let $\hat{\varepsilon}$ be as defined in the statement. Then $\hat{\varepsilon}\in\mathbb{Q}_{>0}$ and $\hat{\varepsilon}\leq\varepsilon$. By the hypothesis of the corollary, applied with $\eta=\hat{\varepsilon}$, Assumption~\ref{Def:SufficientCond} is satisfied. Hence, Theorem~\ref{thm:AlgorithmCorrectness} shows that Algorithm~\ref{algo::aprox1D}, using $\hat{\varepsilon}$ in place of $\varepsilon$, returns a $(1+\hat{\varepsilon})$-covering set~$S$ whose solution projection~$\operatorname{sol}(S)$ is a $(1+\hat{\varepsilon})$-approximation set and therefore also a $(1+\varepsilon)$-approximation set.

By Theorem~\ref{thm:AlgorithmRuntime}, the running time of Algorithm~\ref{algo::aprox1D} is polynomial in $\langle\Pi_{\mathrm{RHS}}\rangle$, $\langle\hat{\varepsilon}\rangle$, and $1/\hat{\varepsilon}$. Moreover, $1/\hat{\varepsilon}=\lceil1/\varepsilon\rceil$ is a positive integer bounded above by $1+1/\varepsilon$. Consequently, \[\langle\hat{\varepsilon}\rangle\in O(\log(2+1/\varepsilon)).\]
Thus, both $1/\hat{\varepsilon}$ and $\langle\hat{\varepsilon}\rangle$ are polynomially bounded in $1/\varepsilon$. The running time is therefore polynomial in $\langle\Pi_{\mathrm{RHS}}\rangle$ and $1/\varepsilon$. Computing the solution projection takes polynomial time in the encoding length of the covering set and therefore preserves this bound. Hence, the procedure that returns this solution projection is an FPTAS for~$\mathcal{P}$.
\end{proof}

In the next section, we first show that monotone RHS-PLPPs admit polynomial-time $(1+\varepsilon)$-covering oracles that return a covering set consisting of a single triple on every call. We then consider certain non-monotone RHS-PLPPs and derive conditions under which the individual covering steps in Algorithm~\ref{algo::aprox1D} can likewise be evaluated in polynomial time using a single triple. A covering-oracle call returning a single triple should not, however, be confused with a singleton approximation set for the entire parameter interval: Algorithm~\ref{algo::aprox1D} may require several such calls, and several solutions may in fact be unavoidable. Moreover, the required number of solutions can depend strongly on the chosen approximation factor, as illustrated by the following example.

\begin{example}
Consider the following RHS-PLPP instance:
\begin{align*}
\left\{\begin{array}{rl}
\displaystyle\min_{x_1,x_2}\, x_1+x_2&\\
\text{s.t.}\,\,\,\,\quad x_1 \geq&\lambda,\\
x_2 \geq&5-\lambda,\\
x_1,x_2 \geq&0
\end{array}\right\}_{\lambda\in[1,5]}.
\end{align*}
For each $\lambda\in[1,5]$, the unique optimal solution is $x^*(\lambda)=(\lambda,5-\lambda)^\mathsf{T}$ with objective value~$5$. Consequently, every optimal solution set contains $S_{\mathrm{RHS}}\coloneqq\bigcup_{\lambda\in[1,5]}\{(\lambda,5-\lambda)^\mathsf{T}\}$ and is therefore infinite.

The solution $\bar{x}=(5,4)^\mathsf{T}$ is feasible for every $\lambda\in[1,5]$ and has objective value~$9$. Hence, $S_{9/5}=\{(5,4)^\mathsf{T}\}$ is a valid $9/5$-approximation set of cardinality one.

More generally, let $(1+\varepsilon)\in(1,2)$ and let $x\coloneqq(a,b)^\mathsf{T}$ be a $(1+\varepsilon)$-approximate solution for at least one parameter value in~$[1,5]$. Since $F^*(\lambda)=5$ throughout~$[1,5]$, we have $a+b\leq5(1+\varepsilon)$. Moreover, the parameter values for which $x$ is feasible form the interval $[5-b,a]\cap[1,5]$, whose length is at most $a+b-5\leq5\varepsilon$. Thus, each such solution can cover a parameter interval of length at most $5\varepsilon$. Since $[1,5]$ has length~$4$, every $(1+\varepsilon)$-approximation set contains at least $\lceil4/(5\varepsilon)\rceil$ solutions. In particular, more than one solution is necessary whenever $\varepsilon<4/5$, and the required cardinality grows without bound as $\varepsilon$ approaches zero.
\end{example}

\section{Efficient Covering Oracles}\label{sec:cov_oracle}

In this section, we identify classes of RHS-PLPP instances for which the covering-oracle calls required by Algorithm~\ref{algo::aprox1D} can be evaluated efficiently.
In Section~\ref{sec:monotone_problems}, we first consider monotone RHS-PLPPs and show that each covering-oracle call required by Algorithm~\ref{algo::aprox1D} can be evaluated in polynomial time by returning a covering set consisting of a single triple. Using this oracle in Algorithm~\ref{algo::aprox1D} yields an FPTAS and, for each rational input approximation factor, an approximation set of minimum cardinality. 
Section~\ref{sec:min_dur_problems} then considers non-monotone RHS-PLPPs. We first provide a polynomial-time test for whether an individual covering-oracle call can be answered by a covering set consisting of a single triple and then derive a sufficient condition, based on the minimal duration of the instance, under which this is guaranteed for every covering-oracle call generated by Algorithm~\ref{algo::aprox1D}. For an exemplary class of problems we refer to as \emph{moving-band RHS-PLPPs}, the minimal duration can be determined explicitly from the problem data, yielding a directly verifiable subclass for which the condition holds. More generally, the procedure yields an FPTAS for classes on which every required covering-oracle call succeeds for every rational input approximation factor greater than one. Whenever all calls in a given execution succeed, the returned covering set contains at most one more triple than a minimum-cardinality covering set.

\subsection{Monotone RHS-PLPPs}\label{sec:monotone_problems}

Following the monotonicity properties established in Lemma~\ref{Lemma:MonotonicityProp}, we call an RHS-PLPP instance \emph{monotonically non-decreasing} if its vector~$r$ satisfies $r\geq0$ and \emph{monotonically non-increasing} if $r\leq0$.

\medskip

For a monotonically non-decreasing instance, $\lambda_{\min}$ is a global minimizer of~$F^*$, whereas, for a monotonically non-increasing instance, $\lambda_{\max}$ is a global minimizer. Accordingly, when applying Algorithm~\ref{algo::aprox1D} to a monotone instance, we implement the initial parametric-search step in Line~\ref{algoline:aprox1D_initial_oracle_call} by setting $\lambda_{\mathrm{PS}}=\lambda_{\min}$ in the former case and $\lambda_{\mathrm{PS}}=\lambda_{\max}$ in the latter case, and by computing $x_{\mathrm{PS}}$ as an optimal solution of $\Pi_{\mathrm{RHS}}(\lambda_{\mathrm{PS}})$. If $r=0$, either choice is valid.

\medskip

A simple example of a monotonically non-decreasing instance is a parametric variant of the classical diet problem. The objective is to determine the least expensive diet satisfying prescribed minimum nutritional requirements. Let $x\in\mathbb{R}_{\geq0}^n$ denote the quantities of $n\in\mathbb{N}$ available food types, $c\in\mathbb{Q}_{>0}^n$ their respective unit prices, and $A\in\mathbb{Q}_{\geq0}^{m\times n}$ the content matrix of $m\in\mathbb{N}$ basic nutritional ingredients, where $a_{ij}$ denotes the amount of nutrient~$i\in[m]$ contained in one unit of food type~$j\in[n]$. Furthermore, let $b\in\mathbb{Q}^m$ denote the vector of minimum daily nutritional requirements.

To model changing dietary requirements, for instance due to different patient groups in a hospital or evolving nutritional recommendations, we introduce a parameter $\lambda\in\Lambda$ and consider the RHS-PLPP instance of the form
    \[
    \min_x\ c^\T x
    \qquad
    \text{s.t.}\quad
    Ax\geq b+\lambda r,\quad x\geq0,
    \]
where $r\in\mathbb{Q}_{\geq0}^m$ describes the direction and magnitude of the changes in the nutritional requirements. For instance, $\lambda$ may represent a person's activity level, with larger values corresponding to higher nutritional demands. Increasing $\lambda$ then tightens the nutritional constraints, resulting in a smaller feasible region and potentially higher minimum costs.

\medskip

For monotone RHS-PLPPs, the covering-oracle calls made by Algorithm~\ref{algo::aprox1D} can be answered using covering sets constructed from a single optimal solution:

\begin{lemma}[$\bm{(1+\varepsilon)}$-Covering Sets for Monotone RHS-PLPPs]
\label{Lemma:AlphaCovCorrectness}
Let $\Pi_{\mathrm{RHS}}=(\Lambda,X,F)$ be a monotone RHS-PLPP instance satisfying Assumption~\ref{Assumption:Standing}, let $\varepsilon\in\mathbb{Q}_{>0}$, and let $L_\varepsilon$ and $R_\varepsilon$ be as defined in \eqref{eq:epsilonParameterBoundaries}.
\begin{enumerate}
    \item Suppose that $r\geq0$.  
    For $\lambda\in\Lambda$, let $x^{\mathrm{R}}$ be an optimal solution of $\Pi_{\mathrm{RHS}}(R_\varepsilon(\lambda))$. Then, $\left\{\bigl(x^{\mathrm{R}},\lambda,R_\varepsilon(\lambda)\bigr)\right\}$ is a $(1+\varepsilon)$-covering set for $[\lambda,R_\varepsilon(\lambda)]$.
    \item Suppose that $r\leq0$.
    For $\lambda\in\Lambda$, let $x^{\mathrm{L}}$ be an optimal solution of $\Pi_{\mathrm{RHS}}(L_\varepsilon(\lambda))$.\, Then, $\left\{\bigl(x^{\mathrm{L}},L_\varepsilon(\lambda),\lambda\bigr)\right\}$ is a $(1+\varepsilon)$-covering set for $[L_\varepsilon(\lambda),\lambda]$.
\end{enumerate}
In both cases, the solution contained in the returned covering set has objective value at most $(1+\varepsilon)\cdot F^*(\lambda)$ and therefore satisfies the common objective threshold of the corresponding covering-oracle call. Furthermore, if $r\geq0$, then $R_\varepsilon(\lambda)$ is the largest parameter value $\mu\in[\lambda,\lambda_{\max}]$ for which a single solution is $(1+\varepsilon)$-approximate throughout $[\lambda,\mu]$; if $r\leq0$, then $L_\varepsilon(\lambda)$ is the smallest parameter value $\mu\in[\lambda_{\min},\lambda]$ for which a single solution is $(1+\varepsilon)$-approximate throughout $[\mu,\lambda]$.
\end{lemma}

\begin{proof}
We first consider the monotonically non-decreasing case $r\geq0$. By Lemma~\ref{Lemma:MonotonicityProp}, $F^*$ is non-decreasing and $X_{\lambda_2}\subseteq X_{\lambda_1}$ whenever $\lambda_1\leq\lambda_2$. Let $\bar{\lambda}\in[\lambda,R_\varepsilon(\lambda)]$. Since $x^{\mathrm{R}}\in X_{R_\varepsilon(\lambda)}$ and $\bar{\lambda}\leq R_\varepsilon(\lambda)$, we have $x^{\mathrm{R}}\in X_{\bar{\lambda}}$. Moreover,
\[
F(x^{\mathrm{R}})
=
F^*(R_\varepsilon(\lambda))
\leq
(1+\varepsilon)\cdot F^*(\lambda).
\]
Since $F^*$ is non-decreasing, $F^*(\lambda)\leq F^*(\bar{\lambda})$, and hence $F(x^{\mathrm{R}})\leq(1+\varepsilon)\cdot F^*(\bar{\lambda})$. Since $\bar{\lambda}\in[\lambda,R_\varepsilon(\lambda)]$ was arbitrary, $x^{\mathrm{R}}$ is therefore feasible and $(1+\varepsilon)$-approximate for every parameter value in this interval.

It remains to show that $R_\varepsilon(\lambda)$ is the largest parameter value with the stated property. Suppose, for contradiction, that there exist a parameter value $\mu\in(R_\varepsilon(\lambda),\lambda_{\max}]$ and a single solution~$x$ that is $(1+\varepsilon)$-approximate throughout $[\lambda,\mu]$. In particular, $x$ is $(1+\varepsilon)$-approximate at~$\lambda$ and feasible at~$\mu$. Therefore, 
\[
F^*(\mu)
\leq
F(x)
\leq
(1+\varepsilon)\cdot F^*(\lambda),
\]
contradicting the definition of $R_\varepsilon(\lambda)$. Together with the first part of the proof, this shows that $R_\varepsilon(\lambda)$ is the largest parameter value with the stated property.

The monotonically non-increasing case follows analogously, with the direction of the parameter interval reversed and $L_\varepsilon(\lambda)$ in place of $R_\varepsilon(\lambda)$.
\end{proof}

Lemma~\ref{Lemma:AlphaCovCorrectness} shows how to answer, using a covering set consisting of a single triple, all covering-oracle calls required by Algorithm~\ref{algo::aprox1D}. More precisely, if $r\geq 0$, then $\lambda_{\mathrm{PS}}=\lambda_{\min}$ and, in the non-degenerate case, only the right while loop is executed. Its covering-oracle calls are of the form
$(1+\varepsilon)\text{-}\operatorname{Cov}\bigl(\lambda,R_\varepsilon(\lambda),(1+\varepsilon)\cdot F^*(\lambda)\bigr)$
and are answered by
$\left\{\bigl(x^{\mathrm{R}},\lambda,R_\varepsilon(\lambda)\bigr)\right\}$.
Analogously, if $r\leq 0$, then $\lambda_{\mathrm{PS}}=\lambda_{\max}$ and, in the non-degenerate case, only the left while loop is executed. Its covering-oracle calls are of the form
$(1+\varepsilon)\text{-}\operatorname{Cov}\bigl(L_\varepsilon(\lambda),\lambda,(1+\varepsilon)\cdot F^*(\lambda)\bigr)$
and are answered by
$\left\{\bigl(x^{\mathrm{L}},L_\varepsilon(\lambda),\lambda\bigr)\right\}$.
Since the parameter value at which an optimal solution is required is part of the covering-oracle input, each such call can be answered by one call to the pointwise optimization oracle
$\operatorname{ALG}_{\mathrm{RHS}}$ from Definition~\ref{Definition:OptimizationOracles}.
Hence, all covering-oracle calls required by Algorithm~\ref{algo::aprox1D} can be evaluated in polynomial time.

In the degenerate case $\lambda_{\min}=\lambda_{\max}$, Algorithm~\ref{algo::aprox1D} directly returns the suitable covering set $\{(x_{\mathrm{PS}},\lambda_{\mathrm{PS}},\lambda_{\mathrm{PS}})\}$
without invoking the covering oracle. Consequently, for every monotone RHS-PLPP and every approximation factor $1+\varepsilon$ with $\varepsilon\in\mathbb{Q}_{>0}$, all covering steps required by Algorithm~\ref{algo::aprox1D} can be performed in polynomial time, which yields the following FPTAS result:

\begin{corollary}\label{Cor:MonotoneFPTAS}
Let $\mathcal{P}_{\mathrm{mon}}$ be the class of monotone RHS-PLPP instances satisfying Assumption~\ref{Assumption:Standing}. Given an instance $\Pi_{\mathrm{RHS}}\in\mathcal{P}_{\mathrm{mon}}$ and a desired approximation factor $1+\varepsilon$ with $\varepsilon>0$, set $\hat{\varepsilon}\coloneqq 1/\lceil 1/\varepsilon\rceil$ and run Algorithm~\ref{algo::aprox1D} on $\Pi_{\mathrm{RHS}}$ using $1+\hat{\varepsilon}$ as the input approximation factor, with the parametric-search call returning $\lambda_{\mathrm{PS}}=\lambda_{\min}$ if $r\geq 0$ and $\lambda_{\mathrm{PS}}=\lambda_{\max}$ otherwise, and each covering-oracle call answered by the corresponding singleton covering set from Lemma~\ref{Lemma:AlphaCovCorrectness} with $\hat{\varepsilon}$ in place of $\varepsilon$. Then the resulting procedure, with output given by the solution projection of the computed covering set, is an FPTAS for~$\mathcal{P}_{\mathrm{mon}}$.
\end{corollary}

\begin{proof}
For every $\Pi_{\mathrm{RHS}}\in\mathcal{P}_{\mathrm{mon}}$ and every $\eta\in\mathbb{Q}_{>0}$, Lemma~\ref{Lemma:AlphaCovCorrectness} and the subsequent construction provide a polynomial-time $(1+\eta)$-covering oracle in the sense of Definition~\ref{Def:ApproxIntervalOracle}. Here, admissible inputs with $\lambda_1=\lambda_2$ are answered in polynomial time by a singleton covering set containing an optimal solution at~$\lambda_1$. Hence, $\Pi_{\mathrm{RHS}}$ satisfies Assumption~\ref{Def:SufficientCond} with $\eta$ in place of $\varepsilon$. Moreover, the implementation of the parametric-search call described above returns an optimal solution of the parametric search problem and runs in polynomial time. The claim therefore follows from Corollary~\ref{cor:GeneralFPTAS}.
\end{proof}

By construction, the $(1+\varepsilon)$-covering oracle from Lemma~\ref{Lemma:AlphaCovCorrectness} extends each covering interval as far as possible from the current parameter value in the direction of the corresponding while loop. The following theorem shows that, as a consequence, Algorithm~\ref{algo::aprox1D} even yields an approximation set of minimum cardinality.

\begin{theorem}[Minimum-Cardinality Approximation Set]\label{thm:Algo1_minimum_cardinality}
Let $\Pi_{\mathrm{RHS}}=(\Lambda,X,F)$ be a monotone RHS-PLPP instance satisfying Assumption~\ref{Assumption:Standing}, and let $\varepsilon\in\mathbb{Q}_{>0}$. Let $S_{\mathrm{cov}}$ be the covering set returned by Algorithm~\ref{algo::aprox1D}, with the parametric-search call returning $\lambda_{\mathrm{PS}}=\lambda_{\min}$ if $r\geq 0$ and $\lambda_{\mathrm{PS}}=\lambda_{\max}$ otherwise, and each covering-oracle call answered by the corresponding singleton covering set from Lemma~\ref{Lemma:AlphaCovCorrectness}. Then $S_{\mathrm{ALG}}\coloneqq\operatorname{sol}(S_{\mathrm{cov}})$ is a $(1+\varepsilon)$-approximation set of minimum cardinality.
\end{theorem}
    \begin{proof}
    If $\lambda_{\min}=\lambda_{\max}$, Algorithm~\ref{algo::aprox1D} returns a covering set whose solution projection consists of a single optimal solution, so the claim follows. 
    
    Hence, suppose that $\lambda_{\min}<\lambda_{\max}$. We first consider the monotonically non-decreasing case, in which the parametric-search call returns $\lambda_{\mathrm{PS}}=\lambda_{\min}$ by the hypothesis of the theorem.
    Suppose that the right while loop of Algorithm~\ref{algo::aprox1D} performs $k$ iterations. Let $\lambda_0\coloneqq\lambda_{\min}$ and let $\lambda_i$, $i\in[k]$, denote the successive parameter boundaries generated by the algorithm. Thus, $\lambda_k=\lambda_{\max}$, and in iteration $i$ the oracle returns the singleton covering set $\{(x_i,\lambda_{i-1},\lambda_i)\}$, where $x_i$ is $(1+\varepsilon)$-approximate throughout $[\lambda_{i-1},\lambda_i]$. By the invariant established in the proof of Theorem~\ref{thm:AlgorithmCorrectness}, $\lambda_i=R_\varepsilon(\lambda_{i-1})$ for every $i\in[k]$. Hence, Lemma~\ref{Lemma:AlphaCovCorrectness} shows that $\lambda_i$ is the largest parameter value $\mu\in[\lambda_{i-1},\lambda_{\max}]$ for which a single solution is $(1+\varepsilon)$-approximate throughout $[\lambda_{i-1},\mu]$. 

Let $D=\{\hat{x}_1,\ldots,\hat{x}_m\}$ be a $(1+\varepsilon)$-approximation set of minimum cardinality. For every $\hat{x}_i\in D$, define
    \[
    \Lambda_{1+\varepsilon}(\hat{x}_i)
    \coloneqq
    \left\{
    \lambda\in\Lambda:
    \hat{x}_i\in X_\lambda,\;
    F(\hat{x}_i)\leq(1+\varepsilon)\cdot F^*(\lambda)
    \right\}.
    \]

By Lemma~\ref{Lemma:DurationProp}, the feasibility region of $\hat{x}_i$ is a closed interval. Since $F^*$ is non-decreasing and continuous, the set of parameter values satisfying $F(\hat{x}_i)\leq(1+\varepsilon)\cdot F^*(\lambda)$ is also a closed interval. The set $\Lambda_{1+\varepsilon}(\hat{x}_i)$ is the intersection of these two intervals and is nonempty, since otherwise $\hat{x}_i$ could be removed from $D$, contradicting its minimality. Consequently, $\Lambda_{1+\varepsilon}(\hat{x}_i)$ is a nonempty closed interval, whose endpoints we denote by $\hat{\lambda}_i^{\mathrm{LB}}$ and $\hat{\lambda}_i^{\mathrm{UB}}$, i.e., $\Lambda_{1+\varepsilon}(\hat{x}_i)=[\hat{\lambda}_i^{\mathrm{LB}},\hat{\lambda}_i^{\mathrm{UB}}]$.

Order the solutions $\hat{x}_1,\ldots,\hat{x}_m$ in $D$ so that $\hat{\lambda}_1^{\mathrm{LB}}<\cdots<\hat{\lambda}_m^{\mathrm{LB}}$. This can be done without loss of generality, since, if two intervals had the same left endpoint, one would be contained in the other and the corresponding solution could be removed, contradicting the minimality of $D$. Minimality further implies $\hat{\lambda}_1^{\mathrm{UB}}<\cdots<\hat{\lambda}_m^{\mathrm{UB}}$. Otherwise, if $\hat{\lambda}_{i+1}^{\mathrm{UB}}\leq\hat{\lambda}_i^{\mathrm{UB}}$ for some $i\in[m-1]$, then $\Lambda_{1+\varepsilon}(\hat{x}_{i+1})\subseteq\Lambda_{1+\varepsilon}(\hat{x}_i)$ and $\hat{x}_{i+1}$ could be removed. Moreover, since the intervals cover $\Lambda$, consecutive intervals overlap, i.e., $\hat{\lambda}_{i+1}^{\mathrm{LB}}\leq\hat{\lambda}_i^{\mathrm{UB}}$ for every $i\in[m-1]$. In particular, $\hat{\lambda}_1^{\mathrm{LB}}=\lambda_{\min}$ and $\hat{\lambda}_m^{\mathrm{UB}}=\lambda_{\max}$.

Since $S_{\mathrm{ALG}}$ is a $(1+\varepsilon)$-approximation set containing at most $k$ solutions, the minimality of $D$ implies $m\leq|S_{\mathrm{ALG}}|\leq k$. We show by induction that $\hat{\lambda}_i^{\mathrm{UB}}\leq\lambda_i$ for every $i\in[m]$.

For $i=1$, the solution $\hat{x}_1$ is $(1+\varepsilon)$-approximate throughout $[\lambda_0,\hat{\lambda}_1^{\mathrm{UB}}]$, since $\hat{\lambda}_1^{\mathrm{LB}}=\lambda_0=\lambda_{\min}$. By Lemma~\ref{Lemma:AlphaCovCorrectness}, $\hat{\lambda}_1^{\mathrm{UB}}\leq R_\varepsilon(\lambda_0)=\lambda_1$.
For the inductive step, assume that $\hat{\lambda}_i^{\mathrm{UB}}\leq\lambda_i$ for some $i<m$. Since consecutive intervals induced by solutions in~$D$ overlap, we have $\hat{\lambda}_{i+1}^{\mathrm{LB}} \leq \hat{\lambda}_i^{\mathrm{UB}} \leq \lambda_i$. Suppose, for contradiction, that $\hat{\lambda}_{i+1}^{\mathrm{UB}}>\lambda_{i+1}$. Then, 
\[[\lambda_i,\hat{\lambda}_{i+1}^{\mathrm{UB}}]\subseteq[\hat{\lambda}_{i+1}^{\mathrm{LB}},\hat{\lambda}_{i+1}^{\mathrm{UB}}]=\Lambda_{1+\varepsilon}(\hat{x}_{i+1}).\]
Thus, $\hat{x}_{i+1}$ is $(1+\varepsilon)$-approximate throughout this interval, contradicting the maximality of $\lambda_{i+1}=R_\varepsilon(\lambda_i)$ established in Lemma~\ref{Lemma:AlphaCovCorrectness}. Therefore, $\hat{\lambda}_{i+1}^{\mathrm{UB}}\leq\lambda_{i+1}$, which completes the induction. 

Since $D$ covers $\Lambda$, we have $\hat{\lambda}_m^{\mathrm{UB}}=\lambda_{\max}$. By the induction result, $\lambda_{\max}=\hat{\lambda}_m^{\mathrm{UB}}\leq\lambda_m\leq\lambda_{\max}$, so $\lambda_m=\lambda_{\max}$. Hence, Algorithm~\ref{algo::aprox1D} terminates after at most $m$ iterations, and therefore $k\leq m$. Together with $m\leq|S_{\mathrm{ALG}}|\leq k$, this yields $|S_{\mathrm{ALG}}|=m$, completing the proof for the monotonically non-decreasing case. The remaining case follows analogously, with $\lambda_{\mathrm{PS}}=\lambda_{\max}$ and the direction of the parameter interval reversed. 
\end{proof}

\subsection{Singleton Covering Sets Beyond Monotonicity}
\label{sec:min_dur_problems}

Monotonicity ensures that the feasible sets are nested, but this property is not required by Algorithm~\ref{algo::aprox1D}. In this subsection, we first give a polynomial-time test for whether a covering-oracle call can be answered by a singleton covering set. We then derive a sufficient condition based on the minimal duration of the instance that guarantees a successful answer for every covering-oracle call required by Algorithm~\ref{algo::aprox1D}, even for non-monotone RHS-PLPP instances. More generally, we obtain an FPTAS for classes on which all required calls succeed for every rational input approximation factor greater than one, and a cardinality guarantee whenever all calls in a given execution succeed.

To answer a covering-oracle call with input $(\lambda_1,\lambda_2,H)$ by a singleton covering set, it suffices to find a solution~$x$ that satisfies $F(x)\leq H$ and is feasible at both endpoints of $[\lambda_1,\lambda_2]$. By Lemma~\ref{Lemma:DurationProp}, feasibility at both endpoints implies feasibility throughout the interval. Moreover, the monotonicity of~$F^*$ on either side of $\lambda_{\mathrm{PS}}$ and Definition~\ref{Def:ApproxIntervalOracle} imply that $H=(1+\varepsilon)\min_{\lambda\in[\lambda_1,\lambda_2]}F^*(\lambda)$, so this solution is $(1+\varepsilon)$-approximate throughout the interval. To test whether such a solution exists, we minimize the left endpoint of a feasibility interval ending at $\lambda_2$ for a right-side covering-oracle call, or maximize the right endpoint of a feasibility interval starting at $\lambda_1$ for a left-side call, subject to the objective threshold~$H$. Both problems can be formulated as linear programming problems.

For a right-side covering-oracle call, where $\lambda_2=R_\varepsilon(\lambda_1)$, define
\begin{equation}
\tag{$\mathrm{EXT}^{\mathrm{R}}(\lambda_2,H)$}
\label{Eq:RightExtensionLP}
\begin{aligned}
\min_{x,L}\quad &L\\
\text{s.t.}\quad
&x\geq0,\\
&Ax\geq b+\lambda_2r,\\
&Ax\geq b+Lr,\\
&c^\T x\leq H,\\
&\lambda_{\mathrm{PS}}\leq L\leq\lambda_2.
\end{aligned}
\end{equation}
Let $(x_{\mathrm{ext}}^{\mathrm{R}},L^{\mathrm{R}})$ denote an optimal solution. Analogously, for a left-side covering-oracle call, where $\lambda_1=L_\varepsilon(\lambda_2)$, define
\begin{equation}
\tag{$\mathrm{EXT}^{\mathrm{L}}(\lambda_1,H)$}
\label{Eq:LeftExtensionLP}
\begin{aligned}
\max_{x,U}\quad &U\\
\text{s.t.}\,\quad
&x\geq0,\\
&Ax\geq b+\lambda_1r,\\
&Ax\geq b+Ur,\\
&c^\T x\leq H,\\
&\lambda_1\leq U\leq\lambda_{\mathrm{PS}}.
\end{aligned}
\end{equation}
Let $(x_{\mathrm{ext}}^{\mathrm{L}},U^{\mathrm{L}})$ denote an optimal solution.

Each extension problem is feasible for every covering-oracle input of the corresponding form. Indeed, if $\lambda_2=R_\varepsilon(\lambda_1)$, then $F^*(\lambda_2)\leq H$ by Definition~\ref{Def:ApproxIntervalOracle} and the definition of $R_\varepsilon$, so an optimal solution at $\lambda_2$ together with $L=\lambda_2$ is feasible for Problem~\eqref{Eq:RightExtensionLP}. The argument for Problem~\eqref{Eq:LeftExtensionLP} is symmetric, using an optimal solution at $\lambda_1$ together with $U=\lambda_1$. Moreover, the objective variables satisfy $\lambda_{\mathrm{PS}}\leq L\leq\lambda_2$ and $\lambda_1\leq U\leq\lambda_{\mathrm{PS}}$, respectively. Since their feasible regions are nonempty polyhedra and their objectives are bounded, both problems attain finite optimal values.

The following lemma characterizes the existence of a singleton covering set satisfying the common objective threshold $H$ in terms of the optimal value of the corresponding extension problem and shows that this condition can be tested in polynomial time.

\begin{lemma}[Polynomial Endpoint-Extension Test for Singleton Covering Sets]
\label{Lemma:EndpointExtensionTest}
Let $(\lambda_1,\lambda_2,H)$ be an input of a $(1+\varepsilon)$-covering oracle in the sense of Definition~\ref{Def:ApproxIntervalOracle}.
\begin{enumerate}[label=\textnormal{(\roman*)}]
\item If $\lambda_2=R_\varepsilon(\lambda_1)$, then there exists a singleton $(1+\varepsilon)$-covering set $\{(x,\lambda_1,\lambda_2)\}$ for $[\lambda_1,\lambda_2]$ satisfying $F(x)\leq H$ if and only if the optimal value $L^{\mathrm{R}}$ of Problem~\eqref{Eq:RightExtensionLP} satisfies $L^{\mathrm{R}}\leq\lambda_1$.
In this case, any optimal solution $(x_{\mathrm{ext}}^{\mathrm{R}},L^{\mathrm{R}})$ of Problem~\eqref{Eq:RightExtensionLP} yields such a covering set $\{(x_{\mathrm{ext}}^{\mathrm{R}},\lambda_1,\lambda_2)\}$. Regardless of whether $L^{\mathrm{R}}\leq\lambda_1$, no solution $\bar{x}\in X$ with $F(\bar{x})\leq H$ is feasible at any parameter value $\mu\in(\lambda_2,\lambda_{\max}]$.\label{lemma:endpoint_ext_case1}

\item If $\lambda_1=L_\varepsilon(\lambda_2)$, then there exists a singleton $(1+\varepsilon)$-covering set $\{(x,\lambda_1,\lambda_2)\}$ for $[\lambda_1,\lambda_2]$ satisfying $F(x)\leq H$ if and only if the optimal value $U^{\mathrm{L}}$ of Problem~\eqref{Eq:LeftExtensionLP} satisfies $U^{\mathrm{L}}\geq\lambda_2$. In this case, any optimal solution $(x_{\mathrm{ext}}^{\mathrm{L}},U^{\mathrm{L}})$ of Problem~\eqref{Eq:LeftExtensionLP} yields such a covering set $\{(x_{\mathrm{ext}}^{\mathrm{L}},\lambda_1,\lambda_2)\}$. Regardless of whether $U^{\mathrm{L}}\geq\lambda_2$, no solution $\bar{x}\in X$ with $F(\bar{x})\leq H$ is feasible at any parameter value $\mu\in[\lambda_{\min},\lambda_1)$.\label{lemma:endpoint_ext_case2}
\end{enumerate}
Both conditions can be tested in time polynomial in $\langle\Pi_{\mathrm{RHS}}\rangle$, $\langle\lambda_{\mathrm{PS}}\rangle$, $\langle\lambda_1\rangle$, $\langle\lambda_2\rangle$, and $\langle H\rangle$.
\end{lemma}

\begin{proof}
Suppose first that $\lambda_2=R_\varepsilon(\lambda_1)$ and $L^{\mathrm{R}}\leq\lambda_1$. By the definition of Problem \eqref{Eq:RightExtensionLP}, $x_{\mathrm{ext}}^{\mathrm{R}}$ is feasible at $L^{\mathrm{R}}$ and $\lambda_2$. Since its feasibility region is an interval by Lemma~\ref{Lemma:DurationProp} and $L^{\mathrm{R}}\leq\lambda_1\leq\lambda_2$, it is feasible throughout $[\lambda_1,\lambda_2]$. Moreover, 
\[F(x_{\mathrm{ext}}^{\mathrm{R}})\leq H=(1+\varepsilon)\cdot F^*(\lambda_1).\] 
Since $F^*$ is non-decreasing on $[\lambda_{\mathrm{PS}},\lambda_{\max}]$, we obtain $F(x_{\mathrm{ext}}^{\mathrm{R}})\leq(1+\varepsilon)\cdot F^*(\mu)$ for every $\mu\in[\lambda_1,\lambda_2]$. Hence, $\{(x_{\mathrm{ext}}^{\mathrm{R}},\lambda_1,\lambda_2)\}$ is a $(1+\varepsilon)$-covering set satisfying the common objective threshold $H$.

Conversely, suppose that $\{(x,\lambda_1,\lambda_2)\}$ is a $(1+\varepsilon)$-covering set for $[\lambda_1,\lambda_2]$ satisfying $F(x)\leq H$. Then $x$ is feasible at both $\lambda_1$ and $\lambda_2$ and satisfies $F(x)\leq H$, so $(x,\lambda_1)$ is feasible for Problem~\eqref{Eq:RightExtensionLP}. By optimality of $(x_{\mathrm{ext}}^{\mathrm{R}},L^{\mathrm{R}})$ for Problem~\eqref{Eq:RightExtensionLP}, we obtain $L^{\mathrm{R}}\leq\lambda_1$.

Now let $\bar{x}\in X$ satisfy $F(\bar{x})\leq H$ and be feasible at some $\mu\in\Lambda$. Then, \[F^*(\mu)\leq F(\bar{x})\leq H=(1+\varepsilon)\cdot F^*(\lambda_1).\] 
By the definition of $R_\varepsilon(\lambda_1)$, this implies $\mu\leq R_\varepsilon(\lambda_1)=\lambda_2$, proving the final assertion in~\ref{lemma:endpoint_ext_case1}. The proof of statement \ref{lemma:endpoint_ext_case2} follows by a symmetric argument. 

Finally, Problems~\eqref{Eq:RightExtensionLP} and~\eqref{Eq:LeftExtensionLP} are rational linear programming problems whose encoding lengths are polynomial in the quantities stated in the lemma. Thus, both conditions can be tested in polynomial time.
\end{proof}

The preceding lemma provides an exact test for whether a covering-oracle input can be answered by a singleton covering set satisfying the common objective threshold~$H$. Algorithm~\ref{algo:VerifiedSingletonOracle} applies this test directly. When invoked by Algorithm~\ref{algo::aprox1D}, the global minimizer $\lambda_{\mathrm{PS}}$ has already been computed and is therefore available to the procedure.

\begin{algorithm}[t]
\caption{A Singleton $(1+\varepsilon)$-Covering-Oracle Procedure}
\label{algo:VerifiedSingletonOracle}
\begin{algorithmic}[1]
\Require An RHS-PLPP instance $\Pi_{\mathrm{RHS}}=(\Lambda,X,F)$, a global minimizer $\lambda_{\mathrm{PS}}$, an approximation factor $1+\varepsilon$ with $\varepsilon\in\mathbb{Q}_{>0}$, and a covering-oracle input $(\lambda_1,\lambda_2,H)$ in the sense of Definition~\ref{Def:ApproxIntervalOracle}.
\Ensure A singleton $(1+\varepsilon)$-covering set $S$ for $[\lambda_1,\lambda_2]$ satisfying the common objective threshold $H$, or \textsc{false} if no such singleton covering set exists.
\If{$\lambda_1=\lambda_2$}
    \State Find $x\leftarrow\mathrm{ALG}_{\mathrm{RHS}}(\lambda_1)$.
    \State \Return $S\leftarrow\{(x,\lambda_1,\lambda_2)\}$.
\ElsIf{$\lambda_2\leq\lambda_{\mathrm{PS}}$}
    \State Find an optimal solution $(x,U)\leftarrow\mathrm{EXT}^{\mathrm{L}}(\lambda_1,H)$.
    \If{$U<\lambda_2$}
        \State \Return \textsc{false}.
    \EndIf
    \State \Return $S\leftarrow\{(x,\lambda_1,\lambda_2)\}$.
\Else
    \State Find an optimal solution $(x,L)\leftarrow\mathrm{EXT}^{\mathrm{R}}(\lambda_2,H)$.
    \If{$L>\lambda_1$}
        \State \Return \textsc{false}.
    \EndIf
    \State \Return $S\leftarrow\{(x,\lambda_1,\lambda_2)\}$.
\EndIf
\end{algorithmic}
\end{algorithm}

We call Algorithm~\ref{algo:VerifiedSingletonOracle} a covering-oracle \emph{procedure} because it attempts to answer a covering-oracle input by a singleton covering set satisfying the common objective threshold~$H$, but may return \textsc{false} if no such set exists. It is therefore not necessarily a covering oracle in the sense of Definition~\ref{Def:ApproxIntervalOracle}, which must return a covering set for every admissible input. The following lemma establishes the correctness and running time of the procedure.

\begin{lemma}[Correctness and Running Time of Algorithm~\ref{algo:VerifiedSingletonOracle}]
\label{Lemma:VerifiedSingletonOracle}
Let $(\lambda_1,\lambda_2,H)$ be a covering-oracle input in the sense of Definition~\ref{Def:ApproxIntervalOracle}. Then Algorithm~\ref{algo:VerifiedSingletonOracle} runs in time polynomial in the encoding lengths of $\Pi_{\mathrm{RHS}}$, $\lambda_{\mathrm{PS}}$, $\lambda_1$, $\lambda_2$, and $H$. If the algorithm returns a set~$S$, then $S$ is a singleton $(1+\varepsilon)$-covering set for $[\lambda_1,\lambda_2]$, and the solution contained in its sole triple satisfies the common objective threshold~$H$. Moreover, the algorithm returns \textnormal{\textsc{false}} if and only if no singleton $(1+\varepsilon)$-covering set for $[\lambda_1,\lambda_2]$ satisfying the common objective threshold~$H$ exists.
\end{lemma}

\begin{proof}
Suppose that $\lambda_1=\lambda_2$, and let $x$ be the optimal solution returned by $\mathrm{ALG}_{\mathrm{RHS}}(\lambda_1)$. By Definition~\ref{Def:ApproxIntervalOracle}, $H\geq F^*(\lambda_1)$. Hence, $x$ is feasible at $\lambda_1$ and satisfies $F(x)=F^*(\lambda_1)\leq H$. Thus, a singleton $(1+\varepsilon)$-covering set satisfying the common objective threshold exists, and the algorithm correctly returns $\{(x,\lambda_1,\lambda_2)\}$.

Now suppose that $\lambda_1<\lambda_2$. If $\lambda_2\leq\lambda_{\mathrm{PS}}$, the algorithm computes an optimal solution $(x,U)$ of Problem~\eqref{Eq:LeftExtensionLP}. By Lemma~\ref{Lemma:EndpointExtensionTest}~\ref{lemma:endpoint_ext_case2}, $U\geq\lambda_2$ holds if and only if there exists a singleton $(1+\varepsilon)$-covering set for $[\lambda_1,\lambda_2]$ satisfying the common objective threshold $H$. If $U\geq\lambda_2$, the algorithm returns such a covering set $\{(x,\lambda_1,\lambda_2)\}$; otherwise, it correctly returns \textsc{false}. If $\lambda_2>\lambda_{\mathrm{PS}}$, the admissible input forms in Definition~\ref{Def:ApproxIntervalOracle} imply $\lambda_1\geq\lambda_{\mathrm{PS}}$. The corresponding conclusions then follow from Lemma~\ref{Lemma:EndpointExtensionTest}~\ref{lemma:endpoint_ext_case1}, using an optimal solution $(x,L)$ of Problem~\eqref{Eq:RightExtensionLP} and the test $L\leq\lambda_1$.

Finally, Problems~\eqref{Eq:LeftExtensionLP} and~\eqref{Eq:RightExtensionLP}, as well as the pointwise problem solved when $\lambda_1=\lambda_2$, are rational linear programming problems whose encoding lengths are polynomial in the quantities stated in the lemma. Each execution solves one of these problems and performs the corresponding comparisons and output operations. Consequently, Algorithm~\ref{algo:VerifiedSingletonOracle} runs within the claimed polynomial bound.
\end{proof}

Note that a \textsc{false} output of Algorithm~\ref{algo:VerifiedSingletonOracle} only means that the requested interval cannot be covered by a single solution satisfying the common objective threshold. It does not exclude the existence of a covering set containing several triples for this interval or of another polynomial-time covering oracle.

\medskip

We next derive a sufficient condition based on the minimal duration of the instance under which Algorithm~\ref{algo:VerifiedSingletonOracle} answers every covering-oracle call required by Algorithm~\ref{algo::aprox1D} by a singleton covering set. Let $f_{\min}\coloneqq F^*(\lambda_{\mathrm{PS}})$ and $H_k\coloneqq(1+\varepsilon)^k\cdot f_{\min}$ for $k\in\mathbb{N}_0$. 

Set $\lambda_0^{\mathrm{L}}\coloneqq\lambda_0^{\mathrm{R}}\coloneqq\lambda_{\mathrm{PS}}$. If $\lambda_{\mathrm{PS}}>\lambda_{\min}$, let
\begin{align*}
    \begin{aligned}N_{\mathrm{L}}&\coloneqq\min\{k\in\mathbb{N}:H_k\geq F^*(\lambda_{\min})\} \quad\text{ and define} \\\lambda_k^{\mathrm{L}}&\coloneqq\min\{\lambda\in\Lambda:F^*(\lambda)\leq H_k\}
    \end{aligned}
\end{align*}
for $k\in[N_{\mathrm{L}}]$. Otherwise, set $N_{\mathrm{L}}\coloneqq0$. Similarly, if $\lambda_{\mathrm{PS}}<\lambda_{\max}$, let
\begin{align*}
    \begin{aligned}N_{\mathrm{R}}&\coloneqq\min\{k\in\mathbb{N}:H_k\geq F^*(\lambda_{\max})\} \quad\text{ and define} \\\lambda_k^{\mathrm{R}}&\coloneqq\max\{\lambda\in\Lambda:F^*(\lambda)\leq H_k\}
    \end{aligned}
\end{align*}
for $k\in[N_{\mathrm{R}}]$. Otherwise, set $N_{\mathrm{R}}\coloneqq0$. 

By the proof of Theorem~\ref{thm:AlgorithmCorrectness}, these sequences prescribe the successive parameter boundaries in Algorithm~\ref{algo::aprox1D}. They are determined independently of the outputs of the covering-oracle calls. The arguments in the proof of Theorem~\ref{thm:AlgorithmRuntime} show that their number and encoding lengths are polynomially bounded in $\langle\Pi_{\mathrm{RHS}}\rangle$, $\langle\varepsilon\rangle$, and $1/\varepsilon$, and that all boundaries can be computed within the same polynomial bound. Define the maximum staircase width by
\[
\Delta_\varepsilon^{\mathrm{grid}}
\coloneqq
\max\left(
\{0\}
\cup
\left\{
\lambda_{k-1}^{\mathrm{L}}-\lambda_k^{\mathrm{L}}:
k\in[N_{\mathrm{L}}]
\right\}
\cup
\left\{
\lambda_k^{\mathrm{R}}-\lambda_{k-1}^{\mathrm{R}}:
k\in[N_{\mathrm{R}}]
\right\}
\right).
\]
With this notation, the sufficient condition can be stated as follows.

\begin{lemma}[Minimal-Duration Condition]
\label{Lemma:MinimalDurationCertificate}
If $d_{\min}^{\Pi_{\mathrm{RHS}}} \geq \Delta_\varepsilon^{\mathrm{grid}}$, then Algorithm~\ref{algo:VerifiedSingletonOracle} returns a singleton $(1+\varepsilon)$-covering set for every covering-oracle call generated by Algorithm~\ref{algo::aprox1D}.
\end{lemma}

\begin{proof}
We first consider the covering-oracle calls generated by the right while loop of Algorithm~\ref{algo::aprox1D}. Fix $k\in[N_{\mathrm{R}}]$ and let $x$ be an optimal solution at $\lambda_k^{\mathrm{R}}$. For every $\mu\in\Lambda(x)$, we have $F^*(\mu)\leq F(x)=F^*(\lambda_k^{\mathrm{R}})\leq H_k$. By the definition of $\lambda_k^{\mathrm{R}}$, this implies $\mu\leq\lambda_k^{\mathrm{R}}$. Since $x$ is feasible at $\lambda_k^{\mathrm{R}}$, the right endpoint of $\Lambda(x)$ is therefore $\lambda_k^{\mathrm{R}}$. Writing $\Lambda(x)=[\lambda^{\mathrm{LB}},\lambda_k^{\mathrm{R}}]$, the hypothesis of the lemma and the definitions of minimal duration and maximum staircase width yield
\[
\lambda_k^{\mathrm{R}}-\lambda^{\mathrm{LB}}
=d_x
\geq d_{\min}^{\Pi_{\mathrm{RHS}}}
\geq\Delta_\varepsilon^{\mathrm{grid}}
\geq\lambda_k^{\mathrm{R}}-\lambda_{k-1}^{\mathrm{R}}.
\]
Consequently, $\lambda^{\mathrm{LB}}\leq\lambda_{k-1}^{\mathrm{R}}$. Hence, $x$ is feasible throughout $[\lambda_{k-1}^{\mathrm{R}},\lambda_k^{\mathrm{R}}]$ and satisfies $F(x)\leq H_k$. Since $H_k=(1+\varepsilon)\cdot F^*(\lambda_{k-1}^{\mathrm{R}})$ and $F^*$ is non-decreasing on this interval, $\{(x,\lambda_{k-1}^{\mathrm{R}},\lambda_k^{\mathrm{R}})\}$ is a $(1+\varepsilon)$-covering set satisfying the common objective threshold~$H_k$. By Lemma~\ref{Lemma:VerifiedSingletonOracle}, Algorithm~\ref{algo:VerifiedSingletonOracle} therefore returns a singleton covering set satisfying that threshold for the corresponding call.

The covering-oracle calls generated by the left while loop are handled analogously: every optimal solution at $\lambda_k^{\mathrm{L}}$ has $\lambda_k^{\mathrm{L}}$ as the left endpoint of its feasibility interval, and the duration bound ensures feasibility throughout $[\lambda_k^{\mathrm{L}},\lambda_{k-1}^{\mathrm{L}}]$. Each such solution also has objective value at most~$H_k$, so Lemma~\ref{Lemma:VerifiedSingletonOracle} again guarantees that Algorithm~\ref{algo:VerifiedSingletonOracle} returns a singleton $(1+\varepsilon)$-covering set for the corresponding call.
\end{proof}

We now consider a class of RHS-PLPP instances that allows arbitrary parameter-inde\-pend\-ent rational polyhedral constraints and rational linear objectives, while restricting the parameter dependence to a moving-band condition. For this class, the minimal duration can be expressed explicitly in terms of the model data, allowing the sufficient condition from Lemma~\ref{Lemma:MinimalDurationCertificate} to be checked directly.

\begin{example}[Moving-Band RHS-PLPPs]
\label{Ex:MovingBandClass}
Let $P\subseteq\mathbb{R}_{\geq0}^n$ be a nonempty parameter-independent rational polyhedron, let $Q\colon\mathbb{R}^n\rightarrow\mathbb{R}$ be an affine function with rational coefficients satisfying $Q(P)=[q_{\min},q_{\max}]\subseteq\Lambda=[\lambda_{\min},\lambda_{\max}]$, and let $\delta^-,\delta^+\in\mathbb{Q}_{>0}$.
Consider RHS-PLPPs with
\[
X_\lambda
\coloneqq
\left\{
x\in P:
Q(x)-\delta^-\leq\lambda\leq Q(x)+\delta^+
\right\},
\qquad
\lambda\in\Lambda,
\]
an arbitrary linear objective~$F(x)=c^\T x$ with $c\in\mathbb{Q}^n$, and suppose that Assumption~\ref{Assumption:Standing} is satisfied.

This structure arises naturally, for example, in blending problems (see \cite{Symonds_2002_Linear_Programming}). The components of $x$ may represent the quantities of the available materials, while $P$ contains parameter-independent blending constraints such as bounds on the component quantities or a prescribed total amount of the blend. The affine function $Q(x)$ represents a quality characteristic determined by the composition, and $\lambda$ denotes its desired target value. The constraint $Q(x)-\delta^-\leq\lambda\leq Q(x)+\delta^+$, which we call the \emph{moving-band condition}, is equivalent to $Q(x)\in[\lambda-\delta^+,\lambda+\delta^-]$. Thus, $\delta^+$ describes the admissible deviation below the target and $\delta^-$ the admissible deviation above it. Increasing $\lambda$ therefore tightens the lower quality requirement $Q(x)\geq\lambda-\delta^+$ while relaxing the upper requirement $-Q(x)\geq-\lambda-\delta^-$, and the feasible sets are generally not nested.

Since $Q(x)\in\Lambda$ for every $x\in P$, every $x\in P$ is feasible at $\lambda=Q(x)$, and hence $X=\bigcup_{\lambda\in\Lambda}X_\lambda=P$. For $x\in P$, the feasibility interval is
\[
\Lambda(x)
=
\left[
\max\{\lambda_{\min},Q(x)-\delta^-\},
\min\{\lambda_{\max},Q(x)+\delta^+\}
\right].
\]
Defining the function~$d\colon\Lambda\rightarrow\mathbb{R}_{\geq0}$ by
\[
d(q)
\coloneqq
\min\{q-\lambda_{\min},\delta^-\}
+
\min\{\lambda_{\max}-q,\delta^+\},
\]
we obtain $d_x=d(Q(x))$ for every $x\in P$. The function $d$ is concave and piecewise linear on $\Lambda$, since each summand is the pointwise minimum of two affine functions. Its minimum over $Q(P)=[q_{\min},q_{\max}]$ is therefore attained at an endpoint, and hence
\[
d_{\min}^{\Pi_{\mathrm{RHS}}}
=
\min\{d(q_{\min}),d(q_{\max})\}.
\]
Thus, the minimal duration can be determined directly from the range of $Q$, the parameter interval, and the two tolerance widths. If $q_{\min}$ and $q_{\max}$ are not given explicitly, they can be computed by minimizing and maximizing~$Q$ over~$P$, respectively, using two rational linear programming problems. In particular, if $Q(P)=\Lambda$, then
\[
d_{\min}^{\Pi_{\mathrm{RHS}}}
=
\min\{
\delta^-,
\delta^+,
\lambda_{\max}-\lambda_{\min}
\}.
\]
If additionally $\delta^-=\delta^+=\delta$ and $\delta\leq\lambda_{\max}-\lambda_{\min}$, then $d_{\min}^{\Pi_{\mathrm{RHS}}}=\delta$.

For a given approximation factor~$1+\varepsilon$, the sufficient condition in Lemma~\ref{Lemma:MinimalDurationCertificate} can now be checked by comparing $d_{\min}^{\Pi_{\mathrm{RHS}}}=\min\{d(q_{\min}),d(q_{\max})\}$ with $\Delta_\varepsilon^{\mathrm{grid}}$. A simpler, though more conservative, condition can be obtained without determining the individual staircase widths. Since every staircase interval lies entirely in one of the two branches,
\[
\Delta_\varepsilon^{\mathrm{grid}}
\leq
\max\{
\lambda_{\mathrm{PS}}-\lambda_{\min},
\lambda_{\max}-\lambda_{\mathrm{PS}}
\}
\]
for every $\varepsilon\in\mathbb{Q}_{>0}$. Consequently, every moving-band instance for which the global minimizer~$\lambda_{\mathrm{PS}}$ returned by the initial parametric-search call satisfies
\begin{equation}
\label{Eq:MovingBandDurationCondition}
\min\{d(q_{\min}),d(q_{\max})\}
\geq
\max\{
\lambda_{\mathrm{PS}}-\lambda_{\min},
\lambda_{\max}-\lambda_{\mathrm{PS}}
\}
\end{equation}
satisfies the minimal-duration condition of Lemma~\ref{Lemma:MinimalDurationCertificate} for every $\varepsilon\in\mathbb{Q}_{>0}$. Thus, \eqref{Eq:MovingBandDurationCondition} gives an $\varepsilon$-independent sufficient condition under which Algorithm~\ref{algo:VerifiedSingletonOracle} returns a singleton $(1+\varepsilon)$-covering set for every covering-oracle call generated by Algorithm~\ref{algo::aprox1D}.
\end{example}

The preceding results yield an FPTAS for classes on which Algorithm~\ref{algo:VerifiedSingletonOracle} successfully answers every covering-oracle call required by Algorithm~\ref{algo::aprox1D} for every rational input approximation factor greater than one.

\begin{corollary}
\label{Cor:VerifiedNonMonotoneFPTAS}
Let $\mathcal{C}$ be a class of RHS-PLPP instances satisfying Assumption~\ref{Assumption:Standing} such that, for every $\Pi_{\mathrm{RHS}}\in\mathcal{C}$ and every $\eta\in\mathbb{Q}_{>0}$, Algorithm~\ref{algo:VerifiedSingletonOracle} does not return \textnormal{\textsc{false}} when used to answer the covering-oracle calls in Algorithm~\ref{algo::aprox1D} with input approximation factor~$1+\eta$. Given an instance~$\Pi_{\mathrm{RHS}}\in\mathcal{C}$ and a desired approximation factor~$1+\varepsilon$ with $\varepsilon>0$, set $\hat{\varepsilon}\coloneqq1/\lceil1/\varepsilon\rceil$ and run Algorithm~\ref{algo::aprox1D} on $\Pi_{\mathrm{RHS}}$ using $1+\hat{\varepsilon}$ as the input approximation factor and Algorithm~\ref{algo:VerifiedSingletonOracle} for each covering-oracle call. Then the resulting procedure, with output given by the solution projection of the computed covering set, is an FPTAS for~$\mathcal{C}$.
\end{corollary}

\begin{proof}
Fix an instance~$\Pi_{\mathrm{RHS}}\in\mathcal{C}$ and $\varepsilon>0$, and let $\hat{\varepsilon}$ be as defined in the statement. Then $\hat{\varepsilon}\in\mathbb{Q}_{>0}$ and $\hat{\varepsilon}\leq\varepsilon$. By the hypothesis of the corollary, Algorithm~\ref{algo:VerifiedSingletonOracle} does not return \textsc{false} on any covering-oracle call generated by Algorithm~\ref{algo::aprox1D} with input approximation factor~$1+\hat{\varepsilon}$. By Lemma~\ref{Lemma:VerifiedSingletonOracle}, every such call therefore returns a singleton $(1+\hat{\varepsilon})$-covering set whose solution satisfies the common objective threshold. Hence, the argument in the proof of Theorem~\ref{thm:AlgorithmCorrectness} applies to this execution: Algorithm~\ref{algo::aprox1D} terminates and returns a $(1+\hat{\varepsilon})$-covering set whose solution projection is a $(1+\hat{\varepsilon})$-approximation set. Since $\hat{\varepsilon}\leq\varepsilon$, this projection is also a $(1+\varepsilon)$-approximation set.

The bounds established in the proof of Theorem~\ref{thm:AlgorithmRuntime}, together with the running-time guarantee in Lemma~\ref{Lemma:VerifiedSingletonOracle}, show that this execution runs in time polynomial in $\langle\Pi_{\mathrm{RHS}}\rangle$, $\langle\hat{\varepsilon}\rangle$, and $1/\hat{\varepsilon}$. Moreover, $1/\hat{\varepsilon}=\lceil1/\varepsilon\rceil\leq1+1/\varepsilon$ and $\langle\hat{\varepsilon}\rangle\in O(\log(2+1/\varepsilon))$. The running time is therefore polynomial in $\langle\Pi_{\mathrm{RHS}}\rangle$ and $1/\varepsilon$. Computing the solution projection takes polynomial time in the encoding length of the covering set and therefore preserves this bound. Hence, the procedure that returns this solution projection is an FPTAS for~$\mathcal{C}$.
\end{proof}

By Lemma~\ref{Lemma:MinimalDurationCertificate}, the hypothesis of Corollary~\ref{Cor:VerifiedNonMonotoneFPTAS} is satisfied whenever every instance~$\Pi_{\mathrm{RHS}}\in\mathcal{C}$ satisfies $d_{\min}^{\Pi_{\mathrm{RHS}}}\geq\Delta_\eta^{\mathrm{grid}}$ for every $\eta\in\mathbb{Q}_{>0}$. In particular, this applies to every class of moving-band RHS-PLPP instances satisfying \eqref{Eq:MovingBandDurationCondition} for the global minimizer~$\lambda_{\mathrm{PS}}$ returned by the initial parametric-search call.

\medskip

Finally, we establish a cardinality guarantee whenever Algorithm~\ref{algo:VerifiedSingletonOracle} successfully answers every covering-oracle call required by Algorithm~\ref{algo::aprox1D}, independently of whether this success follows from the minimal-duration condition. The argument uses the characterization and non-extension assertions in Lemma~\ref{Lemma:EndpointExtensionTest}. In the following result, the cardinality of a covering set refers to its number of triples; thus, two triples containing the same solution but different parameter intervals are counted as distinct elements.

\begin{theorem}[Near-Minimum Cardinality of the Covering Set]
\label{Thm:VerifiedCoverCardinality}
Let $\Pi_{\mathrm{RHS}}=(\Lambda,X,F)$ be an RHS-PLPP instance satisfying Assumption~\ref{Assumption:Standing}, and let $\varepsilon\in\mathbb{Q}_{>0}$. Suppose that Algorithm~\ref{algo::aprox1D} is executed on $\Pi_{\mathrm{RHS}}$ with input approximation factor~$1+\varepsilon$, using Algorithm~\ref{algo:VerifiedSingletonOracle} for every covering-oracle call, and that no call returns \textnormal{\textsc{false}}.
Let $S_{\mathrm{cov}}$ be the resulting $(1+\varepsilon)$-covering set for $\Lambda$, and let $S^*$ be a $(1+\varepsilon)$-covering set of minimum cardinality for $\Lambda$. Then,
\[
|S_{\mathrm{cov}}|
\leq
|S^*|+1.
\]
If $\lambda_{\mathrm{PS}}\in\{\lambda_{\min},\lambda_{\max}\}$, then $S_{\mathrm{cov}}$ has minimum cardinality.
\end{theorem}

\begin{proof}
If $\lambda_{\min}=\lambda_{\max}$, Algorithm~\ref{algo::aprox1D} returns a covering set consisting of a single triple, which has minimum cardinality. Hence, suppose that $\lambda_{\min}<\lambda_{\max}$.
Let $\Lambda_{\mathrm{L}}\coloneqq[\lambda_{\min},\lambda_{\mathrm{PS}}]$ and $\Lambda_{\mathrm{R}}\coloneqq[\lambda_{\mathrm{PS}},\lambda_{\max}]$, and let $S_{\mathrm{L}}$ and $S_{\mathrm{R}}$ denote the sets of triples generated by the left and right while loops, respectively. Then $S_{\mathrm{cov}}=S_{\mathrm{L}}\cup S_{\mathrm{R}}$ and $|S_{\mathrm{cov}}|\leq|S_{\mathrm{L}}|+|S_{\mathrm{R}}|$. 
Within each of the two parameter intervals~$\Lambda_{\mathrm{L}}$ and~$\Lambda_{\mathrm{R}}$, the parameter values at which a fixed solution~$x\in X$ is $(1+\varepsilon)$-approximate form an empty or closed interval. Indeed, its feasibility region is a closed interval by Lemma~\ref{Lemma:DurationProp}, while $F^*$ is monotone and continuous on each of the two parameter intervals.  
Moreover, at the beginning of each executed iteration, the common objective threshold is $H=(1+\varepsilon)\cdot F^*(\lambda_{\mathrm{old}})$ by the invariant established in the proof of Theorem~\ref{thm:AlgorithmCorrectness}.
Every solution~$\bar{x}$ that is $(1+\varepsilon)$-approximate at~$\lambda_{\mathrm{old}}$ therefore satisfies $F(\bar{x})\leq H$, so the non-extension assertions in Lemma~\ref{Lemma:EndpointExtensionTest} apply to every such solution. Since every covering-oracle call succeeds, each iteration reaches the largest possible right endpoint in the right while loop or the smallest possible left endpoint in the left while loop.
Consequently, the induction argument from the proof of Theorem~\ref{thm:Algo1_minimum_cardinality} shows that $S_{\mathrm{L}}$ has minimum cardinality for $\Lambda_{\mathrm{L}}$ whenever $\lambda_{\min}<\lambda_{\mathrm{PS}}$, and that $S_{\mathrm{R}}$ has minimum cardinality for $\Lambda_{\mathrm{R}}$ whenever $\lambda_{\mathrm{PS}}<\lambda_{\max}$.

If $\lambda_{\mathrm{PS}}\in\{\lambda_{\min},\lambda_{\max}\}$, only one while loop is executed, its corresponding parameter interval is~$\Lambda$, and it generates the set~$S_{\mathrm{cov}}$. The minimum-cardinality property established above for this interval therefore implies that $S_{\mathrm{cov}}$ has minimum cardinality among all $(1+\varepsilon)$-covering sets for~$\Lambda$ as claimed. Hence, suppose that $\lambda_{\min}<\lambda_{\mathrm{PS}}<\lambda_{\max}$. For every triple~$s=(x,\lambda_s^{\mathrm{LB}},\lambda_s^{\mathrm{UB}})\in S^*$, let $I(s)\coloneqq[\lambda_s^{\mathrm{LB}},\lambda_s^{\mathrm{UB}}]$ and define $C\coloneqq\{s\in S^*:\lambda_{\mathrm{PS}}\in I(s)\}$. Since the intervals associated with~$S^*$ cover~$\Lambda$, we have $C\neq\emptyset$. We first show that $|C|\leq2$. Suppose that $|C|\geq3$. Choose $s_{\mathrm{L}}\in C$ whose interval has the smallest left endpoint and $s_{\mathrm{R}}\in C$ whose interval has the largest right endpoint. Since both intervals contain $\lambda_{\mathrm{PS}}$, their union contains $I(s)$ for every $s\in C$. Hence, every element of $C\setminus\{s_{\mathrm{L}},s_{\mathrm{R}}\}$ could be removed without destroying the covering property, contradicting the minimum cardinality of~$S^*$. Thus, $1\leq|C|\leq2$.

Every triple in $S^*\setminus C$ has its associated interval entirely contained in one of the two parameter intervals~$\Lambda_{\mathrm{L}}$ and~$\Lambda_{\mathrm{R}}$. If $|C|=1$, let $s_0=(x_0,a_0,b_0)$ be its unique element. Replace it by the two restricted triples $(x_0,a_0,\lambda_{\mathrm{PS}})$ and $(x_0,\lambda_{\mathrm{PS}},b_0)$. Let $S_{\mathrm{L}}^*$ consist of the restricted triple for $\Lambda_{\mathrm{L}}$ and all triples~$s\in S^*\setminus C$ with $I(s)\subseteq\Lambda_{\mathrm{L}}$. Similarly, let $S_{\mathrm{R}}^*$ consist of the restricted triple for $\Lambda_{\mathrm{R}}$ and all triples~$s\in S^*\setminus C$ with $I(s)\subseteq\Lambda_{\mathrm{R}}$. Then $S_{\mathrm{L}}^*$ is a $(1+\varepsilon)$-covering set for~$\Lambda_{\mathrm{L}}$, and $S_{\mathrm{R}}^*$ is a $(1+\varepsilon)$-covering set for~$\Lambda_{\mathrm{R}}$. Since each triple in $S^*\setminus C$ belongs to exactly one of these sets and $s_0$ contributes one restricted triple to each, we have $|S_{\mathrm{L}}^*|+|S_{\mathrm{R}}^*|=|S^*|+1$.

If $|C|=2$, let $s_{\mathrm{L}}$ be the element whose interval has the smaller left endpoint and let $s_{\mathrm{R}}$ be the element whose interval has the larger right endpoint. Minimality of~$S^*$ implies that these are distinct elements. Since both intervals contain $\lambda_{\mathrm{PS}}$, we have $I(s_{\mathrm{R}})\cap\Lambda_{\mathrm{L}}\subseteq I(s_{\mathrm{L}})\cap\Lambda_{\mathrm{L}}$ and $I(s_{\mathrm{L}})\cap\Lambda_{\mathrm{R}}\subseteq I(s_{\mathrm{R}})\cap\Lambda_{\mathrm{R}}$. Restrict the interval of~$s_{\mathrm{L}}$ to $\Lambda_{\mathrm{L}}$ and that of~$s_{\mathrm{R}}$ to $\Lambda_{\mathrm{R}}$. Together with the triples in $S^*\setminus C$ whose associated intervals are contained entirely in~$\Lambda_{\mathrm{L}}$ or in~$\Lambda_{\mathrm{R}}$, the resulting restricted triples yield covering sets~$S_{\mathrm{L}}^*$ and~$S_{\mathrm{R}}^*$ satisfying $|S_{\mathrm{L}}^*|+|S_{\mathrm{R}}^*|=|S^*|$.
Thus, in either case, there exist covering sets for the two parameter intervals such that $|S_{\mathrm{L}}^*|+|S_{\mathrm{R}}^*|\leq|S^*|+1$. By the minimum-cardinality properties established above,
\[
|S_{\mathrm{cov}}|
\leq
|S_{\mathrm{L}}|+|S_{\mathrm{R}}|
\leq
|S_{\mathrm{L}}^*|+|S_{\mathrm{R}}^*|
\leq
|S^*|+1,
\]
which proves the claim.
\end{proof}

Note that the cardinality guarantee in Theorem~\ref{Thm:VerifiedCoverCardinality} concerns covering sets rather than approximation sets. Since, for non-monotone instances, a solution may be $(1+\varepsilon)$-approximate on two disjoint parameter intervals but not between them, this guarantee does not directly imply the analogous bound for approximation sets.

The additive term in Theorem~\ref{Thm:VerifiedCoverCardinality} reflects the separate construction of covering sets for $[\lambda_{\min},\lambda_{\mathrm{PS}}]$ and $[\lambda_{\mathrm{PS}},\lambda_{\max}]$. A one-directional greedy construction could instead proceed from $\lambda_{\min}$ to $\lambda_{\max}$, choosing at each step an interval that starts at the current parameter boundary and has the largest right endpoint for which a single solution is $(1+\varepsilon)$-approximate throughout the interval. The construction would then continue from that right endpoint. If each such maximum can be determined exactly, the standard interval-covering argument yields a covering set of minimum cardinality. The preceding extension tests do not establish a polynomial-time implementation of these steps, and we do not analyze this variant here.

\section{Dual Approximation Algorithm}
\label{sec:DualApproximation}

The primal approximation approaches developed in the previous sections rely on covering the parameter range with polynomially many intervals, on each of which a selected solution remains feasible and provides the required approximation guarantee. Even when all right-hand sides are strictly positive, an optimal solution need not remain feasible when the parameter changes. For example, consider  
\begin{align*}
    \left\{
    \begin{array}{ll}
    \displaystyle \min_{q,z} & z\\
    \,\,\text{s.t.} & z+q\geq2+\lambda,\\
                & z-q\geq2-\lambda, \\
                & q,z\geq0\\
    \end{array}
    \right\}_{\lambda\in[0,1]}.
    \end{align*}
These constraints describe a variant of the moving-band structure from Example~\ref{Ex:MovingBandClass}, with $|q-\lambda|\leq z-2$ and a tolerance width that depends on the decision variable~$z$. The unique optimal solution at parameter value~$\lambda$ is $(q,z)=(\lambda,2)$ and is feasible only at that parameter value. For any $\varepsilon\in\mathbb{Q}_{>0}$, relaxing the right-hand sides by a factor of~$1/(1+\varepsilon)$ makes this solution feasible throughout a parameter interval of positive length while retaining its objective value~$2=F^*(\lambda)$. This motivates the dual approximation approach developed in this section, which preserves the objective-value bound while allowing a controlled relaxation of feasibility.

\medskip

Throughout this section, let $\Pi_{\mathrm{RHS}}=(\Lambda,X,F)$ be an RHS-PLPP instance satisfying the first condition in Assumption~\ref{Assumption:Standing} and the positivity condition~\eqref{eq:positive-right-hand-sides} from Definition~\ref{Def:DualApprox}, and let $\varepsilon\in\mathbb{Q}_{>0}$. Unlike in the primal approximation setting, we do not require the optimal values to be strictly positive. However, the stated assumptions imply $F^*(\lambda)\geq0$ for every $\lambda\in\Lambda$. Indeed, if a feasible solution~$x\in X_\lambda$ satisfied $F(x)<0$, then $sx$ would remain feasible for every $s\geq1$, since $A(sx)\geq s(b+\lambda r)\geq b+\lambda r$, while $F(sx)=sF(x)\to-\infty$ as $s\to\infty$. This would contradict the existence of an optimal solution at~$\lambda$. We use the primal and dual approximation definitions with their unchanged inequalities also when the optimal value is zero.

\medskip

Before developing our dual approximation algorithm, we relate the two approximation notions. Under the assumptions stated below, uniformly scaling the solutions converts a dual approximation set into a primal approximation set with the same approximation factor and cardinality, and vice versa.

\begin{proposition}[Relation between Primal and Dual Approximation Sets]
\label{Prop:PrimalDualApproximationSets}
Let $\Pi_{\mathrm{RHS}}=(\Lambda,X,F)$ be an RHS-PLPP instance satisfying the first condition in Assumption~\ref{Assumption:Standing} and the positivity condition~\eqref{eq:positive-right-hand-sides} from Definition~\ref{Def:DualApprox}. For $\varepsilon\in\mathbb{Q}_{>0}$, define $X^{1+\varepsilon}\coloneqq\bigcup_{\lambda\in\Lambda}X_\lambda^{1+\varepsilon}$ and the $(1+\varepsilon)$-relaxed instance $\Pi_{\mathrm{rel\text{-}RHS}}\coloneqq(\Lambda,X^{1+\varepsilon},F)$, with the same linear objective~$F(x)=c^\T x$. Then, for every parameter subset~$\Lambda'\subseteq\Lambda$ and every set~$S\subseteq\mathbb{R}_{\geq0}^n$, the following statements are equivalent:
\begin{enumerate}[label=\textnormal{(\roman*)}]
    \item $S$ is a dual $(1+\varepsilon)$-approximation set for $\Pi_{\mathrm{RHS}}$ on $\Lambda'$;
    \item $S$ is a primal $(1+\varepsilon)$-approximation set for $\Pi_{\mathrm{rel\text{-}RHS}}$ on $\Lambda'$;
    \item $(1+\varepsilon)S\coloneqq\{(1+\varepsilon)x:x\in S\}$ is a primal $(1+\varepsilon)$-approximation set for $\Pi_{\mathrm{RHS}}$ on~$\Lambda'$.
\end{enumerate}
These correspondences preserve minimum cardinality. Moreover, for every fixed rational approximation factor~$\alpha>1$, a class of RHS-PLPP instances satisfying the above assumptions admits a primal $\alpha$-approximation algorithm if and only if it admits a dual $\alpha$-approximation algorithm. Such a class admits a primal FPTAS if and only if it admits a dual FPTAS.
\end{proposition}

\begin{proof}
Let $F_{\mathrm{rel}}^*(\lambda)$ denote the optimal value of $\Pi_{\mathrm{rel\text{-}RHS}}(\lambda)$. Moreover, we define the bijection $P\colon\mathbb{R}_{\geq0}^n\to\mathbb{R}_{\geq0}^n$ by $P(x)\coloneqq(1+\varepsilon)x$. For every $\lambda\in\Lambda$, we have $x\in X_\lambda^{1+\varepsilon}$ if and only if $P(x)\in X_\lambda$, so $P$ restricts to a bijection between these feasible sets. Substituting $y=(1+\varepsilon)x$ and using linearity of~$F$ yields
\begin{align*}
F_{\mathrm{rel}}^*(\lambda)
\coloneqq{}&
\min\left\{
F(x):
Ax\geq\nicefrac{1}{1+\varepsilon}(b+\lambda r),\
x\geq0
\right\}
\notag\\
={}&
\min\left\{
F(x):
A((1+\varepsilon)x)\geq b+\lambda r,\
(1+\varepsilon)x\geq0
\right\}
\notag\\
={}&
\min\left\{
\nicefrac{1}{(1+\varepsilon)}F(y):
Ay\geq b+\lambda r,\
y\geq0
\right\}
\notag\\
={}&
\nicefrac{1}{(1+\varepsilon)}F^*(\lambda).
\end{align*}

Now, $S$ is a dual $(1+\varepsilon)$-approximation set for $\Pi_{\mathrm{RHS}}$ on $\Lambda'$ if and only if, for every $\lambda\in\Lambda'$, there exists some $x\in S\cap X_\lambda^{1+\varepsilon}$ with $F(x)\leq F^*(\lambda)$. Since $(1+\varepsilon)F_{\mathrm{rel}}^*(\lambda)=F^*(\lambda)$, this is equivalent to $S$ being a primal $(1+\varepsilon)$-approximation set for $\Pi_{\mathrm{rel\text{-}RHS}}$ on $\Lambda'$, proving the equivalence of \textnormal{(i)} and \textnormal{(ii)}.

Moreover, for every $\lambda\in\Lambda'$ and $x\in S$, we have $x\in X_\lambda^{1+\varepsilon}$ if and only if $(1+\varepsilon)x\in X_\lambda$, and $F(x)\leq F^*(\lambda)$ if and only if $F((1+\varepsilon)x)\leq(1+\varepsilon)F^*(\lambda)$. Thus, \textnormal{(i)} is equivalent to \textnormal{(iii)}.
Finally, the bijective mapping~$P$ preserves the cardinality. Together with the equivalence of \textnormal{(i)}--\textnormal{(iii)}, this implies that the minimum cardinalities of the three types of approximation sets coincide.

Multiplying or dividing every solution by $1+\varepsilon$ gives polynomial-time conversions between primal and dual $(1+\varepsilon)$-approximation sets. For every fixed rational approximation factor, these conversions establish the equivalence of the existence of primal and dual approximation algorithms. For the FPTAS equivalence, use the rationalized accuracy~$\hat{\varepsilon}\coloneqq1/\lceil1/\varepsilon\rceil$ as in Corollary~\ref{cor:GeneralFPTAS}; this preserves the requested approximation guarantee and ensures the required running-time bound for the conversions.
\end{proof}

Proposition~\ref{Prop:PrimalDualApproximationSets} shows that, under the stated assumptions, primal and dual approximation are equivalent both for fixed rational approximation factors and with respect to the existence of an FPTAS. The dual formulation nevertheless provides a useful basis for the algorithmic construction: after solving the non-parametric problem at a fixed parameter value, we can explicitly compute an interval on which the resulting optimal solution remains feasible for the relaxed problems, and we can use the monotonicity of~$F^*$ on each branch to ensure the required objective-value bound. The resulting dual $(1+\varepsilon)$-approximation set can then be converted into a primal $(1+\varepsilon)$-approximation set by multiplying its solutions by $1+\varepsilon$.

To record the parameter intervals associated with the selected solutions, we use the following dual counterpart of Definition~\ref{Def:CovSet}.

\begin{definition}[Dual $\bm{(1+\varepsilon)}$-Covering Set] \label{Def:DualCoveringSet}
Let $\varepsilon>0$ and $\lambda_1,\lambda_2\in\Lambda$ with $\lambda_1\leq\lambda_2$. A finite set $S\subseteq\mathbb{R}_{\geq0}^n\times\Lambda\times\Lambda$ is called a \emph{dual $(1+\varepsilon)$-covering set for $[\lambda_1,\lambda_2]$} if, for every $(x,\lambda^{\mathrm{LB}},\lambda^{\mathrm{UB}})\in S$, we have $\lambda_1\leq\lambda^{\mathrm{LB}}\leq \lambda^{\mathrm{UB}}\leq\lambda_2$ and $x$ is a dual $(1+\varepsilon)$-approximate solution for every $\lambda\in[\lambda^{\mathrm{LB}},\lambda^{\mathrm{UB}}]$, and if $[\lambda_1,\lambda_2] = \bigcup_{(x,\lambda^{\mathrm{LB}},\lambda^{\mathrm{UB}})\in S}[\lambda^{\mathrm{LB}},\lambda^{\mathrm{UB}}]$.

Extending the notation introduced after Definition~\ref{Def:CovSet}, define the \emph{solution projection of~$S$}~by 
\[\operatorname{sol}(S)\coloneqq\{x\in\mathbb{R}_{\geq0}^n:\exists\,\lambda^{\mathrm{LB}},\lambda^{\mathrm{UB}}\in\Lambda \text{ with } (x,\lambda^{\mathrm{LB}},\lambda^{\mathrm{UB}})\in S\}.\] 
If $S$ is a dual $(1+\varepsilon)$-covering set for $[\lambda_1,\lambda_2]$, then $\operatorname{sol}(S)$ is a dual $(1+\varepsilon)$-approximation set for $\Pi_{\mathrm{RHS}}$ on $[\lambda_1,\lambda_2]$.
\end{definition}

The positivity condition in~\eqref{eq:positive-right-hand-sides} can be verified in polynomial time. Indeed, for every $i\in[m]$, the function $\lambda\mapsto b_i+\lambda r_i$ is affine on $\Lambda=[\lambda_{\min},\lambda_{\max}]$ and therefore attains its minimum at one of the two endpoints. Consequently, condition~\eqref{eq:positive-right-hand-sides} holds if and only if $b_i+\lambda_{\min}r_i>0$ and $b_i+\lambda_{\max}r_i>0$ for every $i\in[m]$.

\medskip

We next normalize the parameter interval and derive explicit intervals on which every solution feasible at a given parameter value remains feasible for the relaxed problems.
If $\lambda_{\min}=\lambda_{\max}$, we compute an optimal solution~$x^*$ at the unique parameter value by solving the corresponding linear programming problem in polynomial time. Then, $\{(x^*,\lambda_{\min},\lambda_{\max})\}$ is a dual $(1+\varepsilon)$-covering set for $\Lambda$. Hence, in the following, we assume that $\lambda_{\min}<\lambda_{\max}$.

Since $\lambda_{\min}<\lambda_{\max}$, the original RHS-PLPP instance can be reformulated equivalently over the normalized parameter set~$[0,1]$.
Define the affine bijection $h\colon[0,1]\to[\lambda_{\min},\lambda_{\max}]$ by \[h(\mu)\coloneqq \lambda_{\min}+\mu\cdot(\lambda_{\max}-\lambda_{\min}),\] and let $\widetilde b\coloneqq b+\lambda_{\min}\cdot r$ and $\widetilde r\coloneqq(\lambda_{\max}-\lambda_{\min})\cdot r$. Then, for every $\mu\in[0,1]$, substituting $\lambda=h(\mu)$ into the right-hand side of the original instance yields 
\[b+\lambda r = b+h(\mu)\cdot r = \widetilde b+\mu\cdot \widetilde r.\]
Consequently, the normalized feasible sets~$\widetilde X_\mu$ and the normalized $(1+\varepsilon)$-relaxed feasible sets~$\widetilde X_\mu^{1+\varepsilon}$ satisfy
\begin{align}
    \begin{aligned}
        X_{\lambda}&=X_{h(\mu)}=\widetilde X_\mu\coloneq \left\{x\in\mathbb{R}_{\ge0}^n: Ax \geq \widetilde b+\mu\widetilde r \right\} \quad \text{ and} \\
        X_{\lambda}^{1+\varepsilon} &= X_{h(\mu)}^{1+\varepsilon} = \widetilde X_\mu^{1+\varepsilon} \coloneq \left\{x\in\mathbb{R}_{\ge0}^n: Ax \geq \nicefrac{1}{(1+\varepsilon)} \left(\widetilde b+\mu\widetilde r \right)\right\}.
    \end{aligned}
\end{align}
Writing $\widetilde F^*(\mu)\coloneqq \min_{x\in\widetilde X_\mu}F(x)$, we have $\widetilde F^*(\mu)=F^*(h(\mu))$.
For a set $\widetilde S\subseteq \mathbb{R}_{\geq0}^n\times[0,1]\times[0,1]$, we define
\[h(\widetilde S) \coloneqq \left\{\bigl(x,h(\mu^{\mathrm{LB}}),h(\mu^{\mathrm{UB}})\bigr): \bigl(x,\mu^{\mathrm{LB}},\mu^{\mathrm{UB}}\bigr) \in\widetilde S \right\}.\]
Since $h$ is strictly increasing and by the preceding identities for the feasible sets and optimal values, $\widetilde S$ is a dual $(1+\varepsilon)$-covering set for the normalized instance on $[0,1]$ if and only if $h(\widetilde S)$ is a dual $(1+\varepsilon)$-covering set for the original instance on~$\Lambda$. Moreover, $\operatorname{sol}(h(\widetilde S)) =\operatorname{sol}(\widetilde S)$. The normalized instance~$\widetilde{\Pi}_{\mathrm{RHS}}$, with parameter set~$[0,1]$ and data~$A,c,\widetilde b,\widetilde r$, can be computed in polynomial time and has encoding length polynomial in the encoding length~$\langle\Pi_{\mathrm{RHS}}\rangle$ of the original instance~$\Pi_{\mathrm{RHS}}$.

\medskip

The positivity condition~\eqref{eq:positive-right-hand-sides} at the parameter endpoints implies $\widetilde b=b+\lambda_{\min}\cdot r>0$ and $\widetilde b+\widetilde r=b+\lambda_{\max}\cdot r>0$. For notational simplicity, we henceforth work with the normalized instance and omit the tildes, so that $\Lambda=[0,1]$, $b>0$, and $b+r>0$.

\medskip

For $\lambda\in[0,1]$, we call a solution~$x\in\mathbb{R}_{\geq0}^n$ \emph{$(1+\varepsilon)$-relaxed feasible at~$\lambda$} if $x\in X_\lambda^{1+\varepsilon}$.
Moreover, for every $t\in[0,1]$, we let \[\lambda_{1+\varepsilon}^-(t)\coloneqq\min\{\lambda\in[0,t]:b+\lambda r\leq(1+\varepsilon)(b+tr)\}\] be the \emph{minimal $(1+\varepsilon)$-relaxed feasibility endpoint} for~$t$, and let \[\lambda_{1+\varepsilon}^+(t)\coloneqq\max\{\lambda\in[t,1]:b+\lambda r\leq(1+\varepsilon)(b+tr)\}\] be the corresponding \emph{maximal $(1+\varepsilon)$-relaxed feasibility endpoint} for~$t$. Both defining sets are closed subsets of compact intervals and nonempty, since they contain $t$. Thus, the two endpoints are well defined. The proof that they determine an interval of relaxed feasibility relies on the following elementary interpolation property.

\begin{lemma}[Dual Feasibility Property]
\label{lemma:dual_feas_prop}
Let $\lambda_1,\lambda_2\in[0,1]$ with
$\lambda_1\leq\lambda_2$, and let
$\lambda\in[\lambda_1,\lambda_2]$. Then,
\begin{enumerate}[label=\textnormal{(\alph*)}]
    \item
    $X_{\lambda_1}^{1+\varepsilon}
    \cap X_{\lambda_2}^{1+\varepsilon}
    \subseteq X_\lambda^{1+\varepsilon}$,
    \item
    $X_{\lambda_1}
    \cap X_{\lambda_2}^{1+\varepsilon}
    \subseteq X_\lambda^{1+\varepsilon}$, and
    \item
    $X_{\lambda_1}^{1+\varepsilon}
    \cap X_{\lambda_2}
    \subseteq X_\lambda^{1+\varepsilon}$.
\end{enumerate}
\end{lemma}
\begin{proof}
Let $\lambda\in[\lambda_1,\lambda_2]$. Then, there exists $\mu\in[0,1]$ such that
$\lambda=\mu\lambda_1+(1-\mu)\lambda_2$.
To show statement~\textnormal{(a)}, take any $x\in X_{\lambda_1}^{1+\varepsilon}\cap X_{\lambda_2}^{1+\varepsilon}$.
Now, $x$ satisfies both relaxed feasibility constraints, i.e., $Ax\geq\frac{1}{1+\varepsilon}(b+\lambda_1r)$ and $Ax\geq\frac{1}{1+\varepsilon}(b+\lambda_2r)$. Taking a convex combination of these inequalities yields
\[
\begin{aligned}
Ax =\mu Ax+(1-\mu)Ax \geq
\frac{1}{1+\varepsilon}
\bigl(\mu(b+\lambda_1r)+(1-\mu)(b+\lambda_2r)\bigr)= \frac{1}{1+\varepsilon}(b+\lambda r).
\end{aligned}
\]
Hence, $x\in X_\lambda^{1+\varepsilon}$, which shows statement~\textnormal{(a)}. Moreover, the positivity condition~\eqref{eq:positive-right-hand-sides} gives $\frac{1}{1+\varepsilon}(b+\lambda'r)\leq b+\lambda'r$ and hence $X_{\lambda'}\subseteq X_{\lambda'}^{1+\varepsilon}$ for every $\lambda'\in[0,1]$. Statements~\textnormal{(b)} and~\textnormal{(c)} therefore follow immediately from~\textnormal{(a)}.
\end{proof}

The following lemma gives explicit formulas for the relaxed feasibility endpoints and shows that, for every $t\in[0,1]$, each solution feasible at~$t$ remains $(1+\varepsilon)$-relaxed feasible throughout the interval between these endpoints.

For notational convenience, let
$I^+\coloneqq\{i\in[m]:r_i>0\}$ and
$I^-\coloneqq\{i\in[m]:r_i<0\}$.
For fixed~$t\in[0,1]$, the inequalities indexed by~$I^+$ in the condition $b+\lambda r\leq(1+\varepsilon)(b+tr)$ limit how far~$\lambda$ can be increased from~$t$, whereas those indexed by~$I^-$ limit how far it can be decreased. The quantities introduced below capture the most restrictive of these bounds in each direction and yield explicit formulas for the relaxed feasibility endpoints.
Using the usual conventions
$\min\emptyset=+\infty$ and $\max\emptyset=-\infty$, define
\begin{align*}
    \begin{aligned}
    \rho^+&\coloneqq\min_{i\in I^+} b_i/r_i \quad\text{ and} \\
    \rho^-&\coloneqq\max_{i\in I^-} b_i/r_i.
    \end{aligned}
\end{align*}
Since $b>0$, we have $\rho^+\in(0,+\infty]$.
Moreover, since $b+r>0$ and $r_i<0$ for every $i\in I^-$,
we have $b_i/r_i<-1$ and therefore
$\rho^-\in[-\infty,-1)$.

\begin{lemma}[Dual Feasibility Interval]
\label{lem:DualFeasibilityInterval}
For every $t\in[0,1]$, the two $(1+\varepsilon)$-relaxed feasibility endpoints are given by
\begin{align*}
    \begin{aligned}
&\lambda_{1+\varepsilon}^+(t)
=
\min\{1,(1+\varepsilon)t+\varepsilon\rho^+\} \quad \text{and} \\&\lambda_{1+\varepsilon}^-(t)
=
\max\{0,(1+\varepsilon)t+\varepsilon\rho^-\}.
\end{aligned}
\end{align*}
Here,
$\lambda_{1+\varepsilon}^-(t)<t$ whenever $t>0$ and
$\lambda_{1+\varepsilon}^+(t)>t$ whenever $t<1$.
Moreover, for every $x\in X_t$, we have
$x\in X_\lambda^{1+\varepsilon}$ for all
$\lambda\in[\lambda_{1+\varepsilon}^-(t),\lambda_{1+\varepsilon}^+(t)]$.
\end{lemma}

\begin{proof}
Fix $t\in[0,1]$. For every $i\in[m]$, the inequality
$b_i+\lambda r_i\leq(1+\varepsilon)(b_i+tr_i)$ is equivalent to
\begin{align*}
    \begin{aligned}
&\lambda
\leq
(1+\varepsilon)t+\varepsilon\frac{b_i}{r_i}
\quad\text{if }i\in I^+, \quad
       \text{and to} \\&
\lambda
\geq
(1+\varepsilon)t+\varepsilon\frac{b_i}{r_i}
\quad\text{if }i\in I^-.
\end{aligned}
\end{align*}
If $r_i=0$, the inequality
$b_i\leq(1+\varepsilon)b_i$ holds automatically since $b_i>0$.
Since $\lambda=t$ satisfies every componentwise inequality, all lower bounds are at most~$t$ and all upper bounds are at least~$t$. Thus, the lower bounds do not restrict the maximum over $[t,1]$, and the upper bounds do not restrict the minimum over $[0,t]$. Intersecting with these parameter intervals therefore yields
\begin{align*}
    \begin{aligned}
        &\lambda_{1+\varepsilon}^+(t) =
\min\left\{1,\,\min_{i\in I^+} \left((1+\varepsilon) t+\varepsilon\frac{b_i}{r_i}\right)\right\}=\min\{1,(1+\varepsilon)t+\varepsilon\rho^+\} \quad
       \text{and} \\&\lambda_{1+\varepsilon}^-(t) =
\max\left\{0,\,\max_{i\in I^-} \left((1+\varepsilon) t+\varepsilon\frac{b_i}{r_i}\right)\right\} = \max\{0,(1+\varepsilon)t+\varepsilon\rho^-\}.
    \end{aligned}
\end{align*}

Suppose first that $t<1$. If $I^+=\emptyset$, then $\lambda_{1+\varepsilon}^+(t)=1>t$. Otherwise, $\rho^+>0$ and $t\geq0$ imply $(1+\varepsilon)t+\varepsilon\rho^+=t+\varepsilon(t+\rho^+)>t$. Since also $1>t$, both arguments of the minimum defining $\lambda_{1+\varepsilon}^+(t)$ are greater than~$t$, and hence $\lambda_{1+\varepsilon}^+(t)>t$.

Similarly, suppose that $t>0$. If $I^-=\emptyset$, then $\lambda_{1+\varepsilon}^-(t)=0<t$. Otherwise, $\rho^-<-1$ and $t\leq1$ imply $t+\rho^-<0$, so $(1+\varepsilon)t+\varepsilon\rho^-=t+\varepsilon(t+\rho^-)<t$. Since also $0<t$, both arguments of the maximum defining $\lambda_{1+\varepsilon}^-(t)$ are smaller than~$t$, and hence $\lambda_{1+\varepsilon}^-(t)<t$.

\medskip

Finally, let $x\in X_t$. By the definition of $\lambda_{1+\varepsilon}^-(t)$, $b+\lambda_{1+\varepsilon}^-(t)r\leq (1+\varepsilon)(b+tr)$.
Since $Ax\geq b+tr$, it follows that $x\in X_{\lambda_{1+\varepsilon}^-(t)}^{1+\varepsilon}$. Together with $x\in X_t$, Lemma~\ref{lemma:dual_feas_prop}~\textnormal{(c)} therefore yields $x\in X_\lambda^{1+\varepsilon}$ for every $\lambda\in[\lambda_{1+\varepsilon}^-(t),t]$.
Analogously, the definition of $\lambda_{1+\varepsilon}^+(t)$ implies $x\in X_{\lambda_{1+\varepsilon}^+(t)}^{1+\varepsilon}$. Lemma~\ref{lemma:dual_feas_prop}~\textnormal{(b)} then gives $x\in X_\lambda^{1+\varepsilon}$ for every $\lambda\in[t,\lambda_{1+\varepsilon}^+(t)]$. Combining the two intervals proves the claim.
\end{proof}

Lemma~\ref{lem:DualFeasibilityInterval} gives explicit formulas for the relaxed feasibility endpoints associated with a parameter value~$t\in[0,1]$. In Algorithm~\ref{algo:DualApproximation}, successive intervals are determined by using the endpoint of one interval as the starting parameter value for the next. We therefore derive explicit expressions for the parameter sequences obtained by repeatedly computing the minimal or maximal relaxed feasibility endpoint of the preceding parameter value.

\begin{lemma}[Iterated Relaxed-Feasibility Endpoints]
\label{lem:IteratedRelaxedFeasibilityEndpoints}
For an arbitrary $t\in[0,1]$, define two sequences
$(t^-_k)_{k\in\mathbb{N}_0}$ by
$t_0^-\coloneqq t$,
$t_{k+1}^-\coloneqq\lambda_{1+\varepsilon}^-(t_k^-)$
and $(t^+_k)_{k\in\mathbb{N}_0}$ by
$t_0^+\coloneqq t$,
$t_{k+1}^+\coloneqq\lambda_{1+\varepsilon}^+(t_k^+)$.
Then, for every $k\in\mathbb{N}$,
\begin{align*}
    \begin{aligned}
&t_k^+
=
\min\left\{
1,
(1+\varepsilon)^k t
+
\bigl((1+\varepsilon)^k-1\bigr)\rho^+
\right\}
\quad\text{ and } \\
&t_k^-
=
\max\left\{
0,
(1+\varepsilon)^k t
+
\bigl((1+\varepsilon)^k-1\bigr)\rho^-
\right\}.
\end{aligned}
\end{align*}
In particular, $t_{k^-}^-=0$ and $t_{k^+}^+=1$ for some finite
$k^-,k^+\in\mathbb{N}_0$.
\end{lemma}
\begin{proof}
We first consider $(t_k^+)_{k\in\mathbb{N}_0}$.
We show by induction that
\[
t_k^+
=
\min\left\{
1,
(1+\varepsilon)^k t
+
\bigl((1+\varepsilon)^k-1\bigr)\rho^+
\right\} \text{ for every } k\in\mathbb{N}.
\]
For $k=1$, the claim follows directly from
Lemma~\ref{lem:DualFeasibilityInterval}.
Now suppose that the formula holds for some $k\in\mathbb{N}$ and set
$
a_k^+
\coloneqq
(1+\varepsilon)^k t
+
\bigl((1+\varepsilon)^k-1\bigr)\rho^+$.
If $a_k^+\geq1$, then $t_k^+=1$ and therefore
$t_{k+1}^+=\lambda_{1+\varepsilon}^+(1)=1$.
Moreover, since $\rho^+>0$, $(1+\varepsilon)a_k^++\varepsilon\rho^+>1$, so the claimed formula also holds for $k+1$.
Otherwise, $t_k^+=a_k^+$, and
Lemma~\ref{lem:DualFeasibilityInterval} yields
\begin{align*}
t_{k+1}^+ =
\min\{1,(1+\varepsilon)t_k^++\varepsilon\rho^+\}\ =
\min\left\{
1,
(1+\varepsilon)^{k+1}t
+
\bigl((1+\varepsilon)^{k+1}-1\bigr)\rho^+
\right\}.
\end{align*}
Since $t\geq 0$, $\rho^+>0$, and $1+\varepsilon>1$, the expression $(1+\varepsilon)^k t + \bigl((1+\varepsilon)^k-1\bigr)\rho^+$ tends to $+\infty$ as $k\to\infty$. Hence, it is at least $1$ for some finite $k^+\in\mathbb{N}$, and therefore $t_{k^+}^+=1$.

We now consider $(t_k^-)_{k\in\mathbb{N}_0}$.
Analogously, induction yields
\[
t_k^-
=
\max\left\{
0,
(1+\varepsilon)^k t
+
\bigl((1+\varepsilon)^k-1\bigr)\rho^-
\right\} \text{ for every } k\in\mathbb{N}. 
\]
Since $t\leq1$, $\rho^-<-1$, and $1+\varepsilon>1$, the expression $(1+\varepsilon)^k t + \bigl((1+\varepsilon)^k-1\bigr)\rho^-$ tends to~$-\infty$ as $k\to\infty$. It is therefore at most $0$ for some finite $k^-\in\mathbb{N}_0$, which implies $t_{k^-}^-=0$.
\end{proof}

The relaxed feasibility endpoints derived above provide the main ingredient for Algorithm~\ref{algo:DualApproximation}. Starting from a global minimizer of the optimal cost curve, the algorithm proceeds separately towards both endpoints of the parameter interval. At each step, it uses an optimal solution at the current parameter value and computes the corresponding relaxed feasibility endpoint to determine the next interval. The solution remains $(1+\varepsilon)$-relaxed feasible throughout this interval, while the monotonicity of~$F^*$ in the direction of traversal ensures the required objective-value bound. The formulas in Lemma~\ref{lem:DualFeasibilityInterval} allow each endpoint to be computed in time polynomial in the encoding lengths of the instance, $\varepsilon$, and the current parameter value. We subsequently establish the correctness of the algorithm and bound its overall running time.

\begin{algorithm}[t]
\caption{A Dual $(1+\varepsilon)$-Covering Algorithm for RHS-PLPPs}
\label{algo:DualApproximation}
\begin{algorithmic}[1]
\Require A normalized RHS-PLPP instance $\Pi_{\mathrm{RHS}}=([0,1],X,F)$ satisfying the first condition in Assumption~\ref{Assumption:Standing} as well as $b>0$ and $b+r>0$, and $\varepsilon\in\mathbb{Q}_{>0}$.
\Ensure A dual $(1+\varepsilon)$-covering set~$S$ for $[0,1]$.
\State Find
$(x_{\mathrm{PS}},\lambda_{\mathrm{PS}})
\leftarrow
\mathrm{ALG}_{\mathrm{PS}}(\Pi_{\mathrm{RHS}})$.
\State Initialize $S\leftarrow\emptyset$.
\State Set
$t_{\mathrm{old}}\leftarrow\lambda_{\mathrm{PS}}$ and
$x_{\mathrm{old}}\leftarrow x_{\mathrm{PS}}$.
\While{$t_{\mathrm{old}}>0$} \label{DualAlgo:leftwhile}
    \State Set
    $t_{\mathrm{new}}\leftarrow
    \lambda_{1+\varepsilon}^-(t_{\mathrm{old}})$.
    \State Set
    $S\leftarrow
    S\cup
    \{(x_{\mathrm{old}},t_{\mathrm{new}},t_{\mathrm{old}})\}$.
    \State Set $t_{\mathrm{old}}\leftarrow t_{\mathrm{new}}$.
    \If{$t_{\mathrm{old}}>0$}
        \State Find
        $x_{\mathrm{old}}
        \leftarrow
        \mathrm{ALG}_{\mathrm{RHS}}(t_{\mathrm{old}})$.
    \EndIf
\EndWhile
\State Set
$t_{\mathrm{old}}\leftarrow\lambda_{\mathrm{PS}}$ and
$x_{\mathrm{old}}\leftarrow x_{\mathrm{PS}}$.
\While{$t_{\mathrm{old}}<1$}\label{DualAlgo:rightwhile}
    \State Set
    $t_{\mathrm{new}}\leftarrow
    \lambda_{1+\varepsilon}^+(t_{\mathrm{old}})$.
    \State Set
    $S\leftarrow
    S\cup
    \{(x_{\mathrm{old}},t_{\mathrm{old}},t_{\mathrm{new}})\}$.
    \State Set $t_{\mathrm{old}}\leftarrow t_{\mathrm{new}}$.
    \If{$t_{\mathrm{old}}<1$}
        \State Find
        $x_{\mathrm{old}}
        \leftarrow
        \mathrm{ALG}_{\mathrm{RHS}}(t_{\mathrm{old}})$.
    \EndIf
\EndWhile
\State \Return $S$.
\end{algorithmic}
\end{algorithm}

\begin{theorem}[Correctness of
Algorithm~\ref{algo:DualApproximation}]
\label{thm:DualApproximationCorrectness}
Let $\varepsilon\in\mathbb{Q}_{>0}$ and suppose that $\Pi_{\mathrm{RHS}}=([0,1],X,F)$ satisfies the first condition in Assumption~\ref{Assumption:Standing} and the positivity condition~\eqref{eq:positive-right-hand-sides} from Definition~\ref{Def:DualApprox}.
Then each while loop in Algorithm~\ref{algo:DualApproximation}
executes only finitely many iterations, and the returned set $S$ is a dual
$(1+\varepsilon)$-covering set for $[0,1]$. Consequently,
$\operatorname{sol}(S)$ is a dual
$(1+\varepsilon)$-approximation set for the normalized RHS-PLPP instance $\Pi_{\mathrm{RHS}}$.
\end{theorem}

\begin{proof}
Since $\lambda_{\mathrm{PS}}$ is a global minimizer of the convex function $F^*$, the optimal cost curve is non-increasing on $[0,\lambda_{\mathrm{PS}}]$ and non-decreasing on $[\lambda_{\mathrm{PS}},1]$.
We first consider the left while loop beginning in Line~\ref{DualAlgo:leftwhile}. If $\lambda_{\mathrm{PS}}=0$, this loop is not executed, and the interval $[0,\lambda_{\mathrm{PS}}]$ consists only of the single point $0$, which will be covered by the right while loop considered below. Hence, suppose that $\lambda_{\mathrm{PS}}>0$. We show that, upon termination, the set of triples added to~$S$ during the left while loop is a dual $(1+\varepsilon)$-covering set for $[0,\lambda_{\mathrm{PS}}]$.
At the beginning of each executed iteration, $x_{\mathrm{old}}$ is an optimal solution at~$t_{\mathrm{old}}$. This holds initially by the definition of $(x_{\mathrm{PS}},\lambda_{\mathrm{PS}})$ and after every nonterminal iteration by the call to $\mathrm{ALG}_{\mathrm{RHS}}(t_{\mathrm{old}})$.

By Lemma~\ref{lem:DualFeasibilityInterval}, $x_{\mathrm{old}}$ is $(1+\varepsilon)$-relaxed feasible for every $\lambda\in[t_{\mathrm{new}},t_{\mathrm{old}}]$, where $t_{\mathrm{new}} =\lambda_{1+\varepsilon}^-(t_{\mathrm{old}})$. Moreover, for every such $\lambda$, we have $\lambda\leq t_{\mathrm{old}}\leq\lambda_{\mathrm{PS}}$ and therefore
\[
F(x_{\mathrm{old}})
=
F^*(t_{\mathrm{old}})
\leq
F^*(\lambda).
\]
Hence, $\{(x_{\mathrm{old}},t_{\mathrm{new}},t_{\mathrm{old}})\}$ is a dual $(1+\varepsilon)$-covering set for $[t_{\mathrm{new}},t_{\mathrm{old}}]$.
The first interval has right endpoint~$\lambda_{\mathrm{PS}}$. After each interval is added, $t_{\mathrm{old}}$ is set to the computed value~$t_{\mathrm{new}}$, so successive intervals share an endpoint. If the loop terminates, its condition $t_{\mathrm{old}}>0$ is false. Since $t_{\mathrm{old}}\in[0,1]$, this implies $t_{\mathrm{old}}=0$, so the final interval has left endpoint~$0$. Hence, these intervals jointly cover $[0,\lambda_{\mathrm{PS}}]$. Together with the dual approximation guarantee established above, this shows that, upon termination, the triples added during the left while loop form a dual $(1+\varepsilon)$-covering set for $[0,\lambda_{\mathrm{PS}}]$.

\medskip 

It remains to show that the left while loop terminates after finitely many iterations. Define $t_0^{\mathrm{L}}\coloneqq\lambda_{\mathrm{PS}}$ and $t_{k+1}^{\mathrm{L}}\coloneqq\lambda_{1+\varepsilon}^-(t_k^{\mathrm{L}})$ for $k\in\mathbb{N}_0$. These are precisely the successive values of~$t_{\mathrm{old}}$ while the loop is executed. By Lemma~\ref{lem:IteratedRelaxedFeasibilityEndpoints}, there exists a finite $K_{\mathrm{L}}\in\mathbb{N}_0$ such that $t_{K_{\mathrm{L}}}^{\mathrm{L}}=0$. Consequently, the loop terminates after at most $K_{\mathrm{L}}$ iterations.

\medskip

The argument for the right while loop beginning in Line~\ref{DualAlgo:rightwhile} is analogous. Combining the two parts, upon termination of the while loops, the triples added to~$S$ jointly form a dual $(1+\varepsilon)$-covering set for $[0,\lambda_{\mathrm{PS}}] \cup [\lambda_{\mathrm{PS}},1] = [0,1]$. Indeed, if one of the two intervals is a singleton, its unique point $\lambda_{\mathrm{PS}}$ is covered by the other while loop, as discussed above. The statement about $\operatorname{sol}(S)$ follows from Definition~\ref{Def:DualCoveringSet}.
\end{proof}

\begin{theorem}[Running Time of
Algorithm~\ref{algo:DualApproximation}]
\label{thm:DualApproximationRuntime}
Let $\varepsilon\in\mathbb{Q}_{>0}$ and suppose that $\Pi_{\mathrm{RHS}}=([0,1],X,F)$ satisfies the first condition in Assumption~\ref{Assumption:Standing} and the positivity condition~\eqref{eq:positive-right-hand-sides} from Definition~\ref{Def:DualApprox}. Then Algorithm~\ref{algo:DualApproximation} runs in time polynomial in
$\langle\Pi_{\mathrm{RHS}}\rangle$, $\langle\varepsilon\rangle$, and
$1/\varepsilon$. 
\end{theorem}

\begin{proof}
As stated in the paragraph following Definition~\ref{Definition:OptimizationOracles}, the parametric search problem is a linear programming problem of encoding length polynomially bounded in $\langle\Pi_{\mathrm{RHS}}\rangle$. Consequently, the initial parametric-search call can be evaluated in polynomial time, and the returned values $x_{\mathrm{PS}}$ and $\lambda_{\mathrm{PS}}$ have polynomial encoding length.
Moreover, the index sets~$I^+$ and~$I^-$ and every finite value among~$\rho^+$ and~$\rho^-$ can be computed in polynomial time. Each such finite value is one of the ratios~$b_i/r_i$ and therefore has encoding length polynomially bounded in $\langle\Pi_{\mathrm{RHS}}\rangle$.

We first bound the number of iterations of the right loop. If $\lambda_{\mathrm{PS}}=1$, the loop is not executed, and we set $K_{\mathrm R}\coloneqq0$. If $\lambda_{\mathrm{PS}}<1$ and $I^+=\emptyset$, it terminates after one iteration, and we set $K_{\mathrm R}\coloneqq1$. Hence, suppose that $\lambda_{\mathrm{PS}}<1$ and $I^+\neq\emptyset$.
By Lemma~\ref{lem:IteratedRelaxedFeasibilityEndpoints}, the parameter value after $k$ updates is
$
t_k^{\mathrm R}
=
\min\left\{
1,
(1+\varepsilon)^k(\lambda_{\mathrm{PS}}+\rho^+)-\rho^+
\right\}$.
Thus, the loop reaches $1$ once
$
(1+\varepsilon)^k
\geq
\frac{1+\rho^+}{\lambda_{\mathrm{PS}}+\rho^+}.
$
Consequently, using that $\log_2(1+\varepsilon)\geq\varepsilon/(1+\varepsilon)$ for $\varepsilon>0$, the number of iterations is at most
\[
K_{\mathrm R}
:=
\left\lceil
\frac{
\log_2\!\left(
\frac{1+\rho^+}{\lambda_{\mathrm{PS}}+\rho^+}
\right)}
{\log_2(1+\varepsilon)}
\right\rceil \leq \left\lceil \log_2\!\left(
\frac{1+\rho^+}{\lambda_{\mathrm{PS}}+\rho^+}
\right)\cdot \left(1+\frac{1}{\varepsilon}\right)\right\rceil.
\]

Analogously, if $\lambda_{\mathrm{PS}}=0$, the left loop is not executed, and we set $K_{\mathrm L}\coloneqq0$. If $\lambda_{\mathrm{PS}}>0$ and $I^-=\emptyset$, it terminates after one iteration, and we set $K_{\mathrm L}\coloneqq1$. Hence, suppose that $\lambda_{\mathrm{PS}}>0$ and $I^-\neq\emptyset$. Lemma~\ref{lem:IteratedRelaxedFeasibilityEndpoints} then yields
$
t_k^{\mathrm L}
=
\max\left\{
0,
(1+\varepsilon)^k(\lambda_{\mathrm{PS}}+\rho^-)-\rho^-
\right\}$.
Since $\rho^-<-1$, we have $-\rho^--\lambda_{\mathrm{PS}}>0$. The loop reaches $0$ once
$
(1+\varepsilon)^k
\geq
\frac{-\rho^-}{-\rho^--\lambda_{\mathrm{PS}}}$,
and therefore performs at most
\[
K_{\mathrm L}
:=
\left\lceil
\frac{
\log_2\!\left(
\frac{-\rho^-}{-\rho^--\lambda_{\mathrm{PS}}}
\right)}
{\log_2(1+\varepsilon)}
\right\rceil \leq \left\lceil \log_2\!\left(
\frac{-\rho^-}{-\rho^--\lambda_{\mathrm{PS}}}
\right)\cdot \left(1+\frac{1}{\varepsilon}\right)\right\rceil
\]
iterations. In each case involving a logarithm, its argument is a rational number greater than~$1$ whose encoding length is polynomially bounded in $\langle\Pi_{\mathrm{RHS}}\rangle$. Its logarithm is therefore polynomially bounded in $\langle\Pi_{\mathrm{RHS}}\rangle$ as well. Together with the bounds for the skipped and one-step loops, this shows that $K_{\mathrm L}$ and $K_{\mathrm R}$ are polynomially bounded in $\langle\Pi_{\mathrm{RHS}}\rangle$ and $1/\varepsilon$.

\medskip

It remains to bound the encoding lengths of the generated parameter values. The initial value~$\lambda_{\mathrm{PS}}$ already has polynomially bounded encoding length, and the terminal values~$0$ and~$1$ have constant encoding length. A loop whose corresponding index set is empty generates only its terminal endpoint, if it is executed. If $\rho^+$ is finite, every nonterminal parameter value generated by the right loop has the form $(1+\varepsilon)^k(\lambda_{\mathrm{PS}}+\rho^+)-\rho^+$. Similarly, if $\rho^-$ is finite, every nonterminal parameter value generated by the left loop has the form $(1+\varepsilon)^k(\lambda_{\mathrm{PS}}+\rho^-)-\rho^-$. In both cases, $k\leq\max\{K_{\mathrm L},K_{\mathrm R}\}$. We have $\langle(1+\varepsilon)^k\rangle\in O(1+k\langle\varepsilon\rangle)$, and $\lambda_{\mathrm{PS}}$ and the finite values~$\rho^\pm$ have encoding lengths polynomially bounded in $\langle\Pi_{\mathrm{RHS}}\rangle$. Together with the bound on~$k$, this shows that every generated parameter value has encoding length polynomially bounded in $\langle\Pi_{\mathrm{RHS}}\rangle$, $\langle\varepsilon\rangle$, and $1/\varepsilon$.

Consequently, every call to $\mathrm{ALG}_{\mathrm{RHS}}$ solves a linear programming problem whose encoding length is polynomially bounded in $\langle\Pi_{\mathrm{RHS}}\rangle$, $\langle\varepsilon\rangle$, and $1/\varepsilon$, and can be evaluated in time polynomial in these quantities. An optimal rational solution of polynomially bounded encoding length can be returned. The endpoint computations, comparisons, and updates of~$S$ also take polynomial time in the same quantities. Each iteration produces exactly one triple consisting of such a solution and two polynomially encoded parameter values. Since the total number of iterations is polynomially bounded, the returned dual $(1+\varepsilon)$-covering set has polynomially bounded cardinality and total encoding length, and the overall running time of Algorithm~\ref{algo:DualApproximation} is polynomial in $\langle\Pi_{\mathrm{RHS}}\rangle$, $\langle\varepsilon\rangle$, and $1/\varepsilon$ as claimed.
\end{proof}

Combining the preceding results with the normalization described earlier yields a dual FPTAS that also records the parameter intervals on which its solutions provide the required approximation guarantee.

\begin{corollary}
\label{cor:DualFPTAS}
Let $\mathcal{P}_{\mathrm{dual}}$ be the class of RHS-PLPP instances satisfying the first condition in Assumption~\ref{Assumption:Standing} and the positivity condition~\eqref{eq:positive-right-hand-sides} from Definition~\ref{Def:DualApprox}.
Given an instance~$\Pi_{\mathrm{RHS}}\in\mathcal{P}_{\mathrm{dual}}$ and a desired dual approximation factor~$1+\varepsilon$ with $\varepsilon>0$, set $\hat{\varepsilon}\coloneqq1/\lceil1/\varepsilon\rceil$. If $\lambda_{\min}=\lambda_{\max}$, compute an optimal solution at the unique parameter value and form the corresponding singleton dual covering set. Otherwise, normalize the parameter interval to $[0,1]$ by the change of parameter $\lambda=\lambda_{\min}+\mu(\lambda_{\max}-\lambda_{\min})$, run Algorithm~\ref{algo:DualApproximation} on the resulting normalized instance with input approximation factor~$1+\hat{\varepsilon}$, and transform the returned interval endpoints back to~$\Lambda$. Then the resulting procedure, with output given by the solution projection of the computed covering set, is a dual FPTAS for~$\mathcal{P}_{\mathrm{dual}}$.
\end{corollary}
\begin{proof}
Fix $\Pi_{\mathrm{RHS}}\in\mathcal{P}_{\mathrm{dual}}$ and $\varepsilon>0$, and let $\hat{\varepsilon}$ be as defined in the statement. Then $\hat{\varepsilon}\in\mathbb{Q}_{>0}$ and $\hat{\varepsilon}\leq\varepsilon$. If $\lambda_{\min}=\lambda_{\max}$, an optimal solution at the unique parameter value yields a singleton dual $(1+\varepsilon)$-covering set and can be computed in polynomial time. Hence, suppose that $\lambda_{\min}<\lambda_{\max}$. Let $\widetilde{\Pi}_{\mathrm{RHS}}$ be the normalized instance obtained under the change of parameter $\lambda=h(\mu)\coloneqq\lambda_{\min}+\mu(\lambda_{\max}-\lambda_{\min})$, and apply Algorithm~\ref{algo:DualApproximation} to this instance with input approximation factor~$1+\hat{\varepsilon}$. By Theorem~\ref{thm:DualApproximationCorrectness}, the algorithm returns a dual $(1+\hat{\varepsilon})$-covering set
$\widetilde S$ for $[0,1]$.

Since $\hat{\varepsilon}\leq\varepsilon$ and the normalized right-hand sides are strictly positive, for every $\mu\in[0,1]$, we have $\widetilde X_\mu^{1+\hat{\varepsilon}}\subseteq\widetilde X_\mu^{1+\varepsilon}$. Hence, $\widetilde S$ is also a dual $(1+\varepsilon)$-covering set for the normalized instance.
By the correspondence between covering sets under the normalization, $h(\widetilde S)$ is a dual $(1+\varepsilon)$-covering set for the original parameter interval $\Lambda$, and $\operatorname{sol}(h(\widetilde S))=\operatorname{sol}(\widetilde S)$ is a dual $(1+\varepsilon)$-approximation set for the original instance.

Moreover, $1/\hat{\varepsilon} =\lceil1/\varepsilon\rceil\leq1+1/\varepsilon$ and $\langle\hat{\varepsilon}\rangle \in O(\log(2+1/\varepsilon))$. The normalization, the application of~$h$ to the interval endpoints, and the running time from Theorem~\ref{thm:DualApproximationRuntime} are therefore polynomial in the encoding length of the original instance and in $1/\varepsilon$. Computing the solution projection takes polynomial time in the encoding length of the covering set and therefore preserves this bound. Hence, the resulting algorithm is a dual FPTAS. 
\end{proof}

The dual FPTAS complements the primal approximation framework developed in the preceding sections by preserving the objective-value bound while allowing a controlled relaxation of the parameter-dependent right-hand sides. The computed dual covering set also records the parameter intervals on which each selected solution provides the required guarantee. By Proposition~\ref{Prop:PrimalDualApproximationSets}, this construction additionally yields a primal FPTAS for~$\mathcal{P}_{\mathrm{dual}}$.

\section{Conclusion}

This paper develops multiplicative primal and dual approximation approaches for RHS-PLPPs, providing, to the best of our knowledge, the first results of this kind. Although the optimal cost curve is continuous, piecewise linear, and convex, an optimal solution set may require infinitely many distinct solutions, and a representation by optimal bases may require super-polynomially many bases. Under the conditions established in this paper, approximation sets of polynomially bounded cardinality can be computed efficiently and provide the required guarantees throughout the parameter interval.

\medskip

The connection between primal and dual approximation provides a useful perspective beyond the individual algorithms. Under the stated assumptions, uniform scaling transfers approximation sets between the two notions while preserving their approximation factor and cardinality. Consequently, algorithmic and cardinality results established for one notion can also inform the other. This correspondence makes the dual formulation useful even when the ultimate aim is to obtain solutions that satisfy the original constraints exactly.

For applications involving repeated queries at different parameter values, the explicit intervals associated with the selected solutions provide certificates of where their guarantees apply. Once such a representation has been computed, a solution with the required guarantee can be selected for a given parameter value without solving another linear programming problem. The distinction between covering sets and approximation sets also clarifies two aspects of representation size: the number of distinct solutions and the number of intervals needed to describe their use.

\medskip

Several directions for future research arise from these results. One is to identify broader classes of non-monotone RHS-PLPPs for which the covering subproblems can be solved efficiently, particularly when more than one solution is required to cover a parameter interval. Another is to establish cardinality guarantees for approximation sets in the non-monotone case and for the dual approximation construction. Extending the approximation and cardinality results to multi-parametric linear programming problems is a further direction. Finally, mixed-integer RHS-parametric problems provide a natural setting in which to investigate the limits of these approaches: fixed solutions still have interval-shaped feasibility regions in the single-parameter case, but the convexity of the optimal cost curve and the scaling correspondence need not be preserved.

\printbibliography[heading=bibintoc]{}

    \section*{Appendix}
\phantomsection
\label{app:proof-basis-intractability}

    This appendix provides a proof for Lemma~\ref{Lemma:LinearSections}, and hence shows that arbitrarily many distinct bases may be required in a basic optimal solution set, even when the optimal cost curve is constant. 
    
    \begin{proof}
    Fix~$q\in\mathbb{N}$. Define the sequence $(b_k)_{k\in\mathbb{N}_0}$ by $b_0=0.2$ and $b_k=2b_{k-1}+0.1k$ for~$k\in\mathbb{N}$. Equivalently,
    \[
    b_k=0.4\cdot2^k-0.1k-0.2.
    \]
    Consider the following RHS-PLPP instance:
    \begin{align*}
    \left\{
    \begin{array}{ll}
    \displaystyle \min_{x_1,x_2} & x_2\\
    \,\,\text{s.t.} & x_1\geq 2-b_k+2^k\lambda \qquad \text{for }k=0,\ldots,q-1,\\
                & x_2\geq1,\\
                & (x_1,x_2)\in\mathbb{R}_{\geq0}^2
    \end{array}
    \right\}_{\lambda\in[0,2]}.
    \end{align*}
    
    For every~$\lambda\in[0,2]$, the corresponding non-parametric problem is feasible. Since the objective is to minimize~$x_2$ and the only constraint involving~$x_2$ is $x_2\geq1$, its optimal objective value is~$1$. Hence, the optimal cost curve is constant on the entire parameter set.
    
    For~$k=0,\ldots,q-1$, let $g_k(\lambda)\coloneq2-b_k+2^k\lambda$.
    We show that each~$g_k$ uniquely determines the lower bound on~$x_1$ over a nonempty parameter interval. Define $\lambda_0\coloneq0$, $\lambda_q\coloneq2$, and
    \[
    \lambda_k\coloneq0.4-\frac{1}{10\cdot2^{k-1}}
    \qquad\text{for }k=1,\ldots,q-1.
    \]
    Then,
    \[
    0=\lambda_0<\lambda_1<\cdots<\lambda_{q-1}<0.4<\lambda_q.
    \]
    
    We claim that, for every~$k\in\{0,\ldots,q-1\}$ and every~$\bar\lambda\in(\lambda_k,\lambda_{k+1})$, one has $g_k(\bar\lambda)>g_j(\bar\lambda)$ for all~$j\neq k$. To prove this, fix~$k\in\{0,\ldots,q-1\}$ and~$\bar\lambda\in(\lambda_k,\lambda_{k+1})$.
    
    First, let~$j>k$, and let~$\widehat{\lambda}_{j,k}$ denote the unique parameter value satisfying $g_k(\widehat{\lambda}_{j,k})=g_j(\widehat{\lambda}_{j,k})$.
    Solving this equation yields
    \[
    \widehat{\lambda}_{j,k}
    =\frac{b_j-b_k}{2^j-2^k}
    =0.4-\frac{j-k}{10(2^j-2^k)}
    =0.4-\frac{j-k}{10\cdot2^k(2^{j-k}-1)}.
    \]
    Since $2^{j-k}-1\geq j-k$, we obtain
    \[
    \widehat{\lambda}_{j,k}
    \geq0.4-\frac{1}{10\cdot2^k}
    =\lambda_{k+1}
    >\bar\lambda.
    \]
    Moreover, $g_j$ has the larger slope because $2^j>2^k$. Since~$\bar\lambda$ lies strictly before the parameter value at which the two functions attain the same value, it follows that $g_k(\bar\lambda)>g_j(\bar\lambda)$.
    
    Now, let~$j<k$. The unique parameter value~$\widehat{\lambda}_{j,k}=\widehat{\lambda}_{k,j}$ satisfying $g_k(\widehat{\lambda}_{j,k})=g_j(\widehat{\lambda}_{j,k})$ is given by
    \[
    \widehat{\lambda}_{j,k}
    =0.4-\frac{k-j}{10\cdot2^j(2^{k-j}-1)}.
    \]
    Since $2^{k-j}-1\leq(k-j)2^{k-j-1}$, we obtain
    \[
    \widehat{\lambda}_{j,k}
    \leq0.4-\frac{1}{10\cdot2^{k-1}}
    =\lambda_k
    <\bar\lambda.
    \]
    Moreover, $g_k$ has the larger slope because $2^k>2^j$. Since~$\bar\lambda$ lies strictly after the parameter value at which the two functions attain the same value, it follows that $g_k(\bar\lambda)>g_j(\bar\lambda)$.
    
    Consequently, on each interval~$(\lambda_k,\lambda_{k+1})$, the inequality $x_1\geq g_k(\lambda)$ is the unique tightest lower-bound constraint on~$x_1$. Together with~$x_2\geq1$, it determines the unique optimal vertex~$\bigl(g_k(\lambda),1\bigr)^\T$.
    After introducing surplus variables, this vertex is nondegenerate: the surplus variables associated with the constraints $x_1\geq g_k(\lambda)$ and $x_2\geq1$ are nonbasic, while all other surplus variables are positive. 
    It therefore has a unique basis, which is dual feasible because the associated basic solution is nondegenerate and optimal. These bases are distinct for different values of~$k$.
    Hence, every basic optimal solution set must contain at least one distinct basis for each of the~$q$ intervals~$(\lambda_k,\lambda_{k+1})$ and therefore at least~$q$ bases.
    \end{proof}

\end{document}